\documentstyle[12pt]{article}
\newcommand{\preprint}[1]{#1}
\newcommand{\journal}[1]{}
\newcommand{\hide}[1]{}
\newtheorem{thm}{Theorem}[section]
\newtheorem{prop}[thm]{Proposition}

\newtheorem{conj}[thm]{Conjecture}
\newtheorem{thm-defi}[thm]{Theorem/Definition}
\newtheorem{example}[thm]{Example}
\newtheorem{cor}[thm]{Corollary}
\newtheorem{question}[thm]{Question}
\newtheorem{new-lemma}[thm]{Lemma}
\newtheorem{defi}[thm]{Definition}
\newtheorem{rem}[thm]{Remark}

\newtheorem{condition}[thm]{Condition}

\newcommand{\A}{{\cal A}}
\newcommand{\C}{{\cal C}}

\newcommand{\E}{{\cal E}}
\newcommand{\F}{{\cal F}}
\newcommand{\G}{{\cal G}}

\newcommand{\R}{{\cal R}}

\newcommand{\M}{{\cal M}}

\newcommand{\W}{{\cal W}}
\newcommand{\X}{{\cal X}}

\renewcommand{\P}{{\cal P}}
\newcommand{\PP}{{\Bbb P}}

\newcommand{\RealNumbers}{{\Bbb R}}
\newcommand{\Integers}{{\Bbb Z}}
\newcommand{\ComplexNumbers}{{\Bbb C}}
\newcommand{\RationalNumbers}{{\Bbb Q}}
\newcommand{\LieAlg}[1]{{\frak #1}}

\newcommand{\monrep}{{mon}}
\newcommand{\monclass}{{\overline{mon}}}

\newcommand{\LongIsomRightArrow}{\stackrel{\cong}{\longrightarrow}}
\newcommand{\LongIsomLeftArrow}{\stackrel{\cong}{\longleftarrow}}
\newcommand{\RightArrowOf}[1]{\stackrel{#1}{\rightarrow}}

\newcommand{\LongRightArrowOf}[1]{\stackrel{#1}{\longrightarrow}}
\newcommand{\LongIsomRightArrowOf}[1]{
\stackrel
{\stackrel{#1}{\cong}}
{\longrightarrow}
}
\newcommand{\LongIsomLefttArrowOf}[1]{
\stackrel
{\stackrel{#1}{\cong}}
{\longleftarrow}
}

\newcommand{\StructureSheaf}[1]{{\cal O}_{#1}}
\newcommand{\EndProof}{\nolinebreak\hfill  $\Box$}
\newcommand{\restricted}[2]{#1_{\mid_{#2}}}

\newcommand{\rank}{{\rm rank}}
\newcommand{\coker}{{\rm coker}}

\newcommand{\Sym}{{\rm Sym}}
\newcommand{\Ext}{{\rm Ext}}

\newcommand{\Hom}{{\rm Hom}}

\newcommand{\Abs}[1]{\mid\!#1\!\mid}

\newcommand{\Wedge}[1]{\stackrel{#1}{\wedge}}

\newcommand{\lcm}{{\rm lcm}}

\input amssymb.sty

\begin{document}
\begin{center}
\begin{large}
{\bf 
\noindent
Integral constraints on 
the monodromy group of the hyperk\"{a}hler 
resolution of a 
symmetric product of a $K3$ 
surface
}
\end{large}
\\
Eyal Markman\footnote{Department of mathematics and statistics, 
University of Massachusetts, Amherst, MA 01003, USA. 
E-mail: markman@math.umass.edu}
\end{center}

{\small 
{\bf Abstract:}
Let $S^{[n]}$ be the Hilbert scheme of length $n$ subschemes of a $K3$ 
surface $S$. 
$H^2(S^{[n]},\Integers)$ is endowed with the
Beauville-Bogomolov bilinear form. 
Denote by $Mon$ the subgroup of $GL[H^*(S^{[n]},\Integers)]$
generated by monodromy operators, and let $Mon^2$
be its image in $OH^2(S^{[n]},\Integers)$.
We prove that $Mon^2$ is the subgroup 
generated by reflections with respect to $+2$ and $-2$ classes
(Theorem \ref{thm-introduction-Mon-2-is-W}). 
Thus $Mon^2$ does not
surject onto $OH^2(S^{[n]},\Integers)/(\pm 1)$, when
$n-1$ is not a prime power.

As a consequence, we get 
counter examples to a version of the weight $2$ Torelli question
for hyperk\"{a}hler varieties $X$ deformation equivalent to $S^{[n]}$. 
The weight $2$ Hodge structure on $H^2(X,\Integers)$
does not determine the bimeromorphic class of $X$, whenever 
$n-1$ is not a prime power (the first case being $n=7$). 
There are at least $2^{\rho(n-1)-1}$ distinct bimeromorphic classes of $X$
with a given generic weight $2$ Hodge structure, where $\rho(n-1)$ is
the Euler number of $n-1$. 

The second main result 
states, that if a monodromy operator acts as the identity on 
$H^2(S^{[n]},\Integers)$, then it acts as the identity on 
$H^{k}(S^{[n]},\Integers)$, $0\!\leq\! k\! \leq\! n\!+\!2$ (Theorem 
\ref{thm-monodromy-operator-is-determined-by-its-weight-2-action}).
We conclude the injectivity of the restriction homomorphism 
$Mon\!\rightarrow\! Mon^2$,
if $n\!\equiv\! 0$, or $n\!\equiv\! 1$ modulo 4
(Corollary
\ref{cor-reduction-of-monodromy-conjecture-to-generators-conjecture}).
}

{\scriptsize 
\tableofcontents
} 

\section{Introduction}
\label{sec-introduction}
%

%
\subsection{The main results}
An {\em irreducible holomorphic symplectic manifold} is a simply connected
compact K\"{a}hler manifold $X$, such that $H^0(X,\Omega^2_X)$ is generated
by an everywhere non-degenerate holomorphic two-form.
An irreducible holomorphic symplectic manifold of real dimension $4n$ 
admits a Riemannian metric with holonomy $Sp(n)$ 
\cite{beauville-varieties-with-zero-c-1}. Such a metric is called 
{\em hyperk\"{a}hler}.

Let $S$ be a smooth complex projective $K3$ surface and
$S^{[n]}$ the Hilbert scheme of length $n$ zero dimensional subschemes of $S$.
The  Hilbert scheme $S^{[n]}$ is a smooth projective 
variety, which is also an irreducible holomorphic symplectic manifold
\cite{beauville-varieties-with-zero-c-1}.
We study in this paper the monodromy group of
irreducible holomorphic symplectic manifolds, which are  
deformation equivalent to $S^{[n]}$.

\begin{defi}
\label{def-monodromy}
{\rm
Let $X$ be an irreducible holomorphic symplectic manifold. 
An automorphism $g$ of the cohomology ring 
$H^*(X,\Integers)$ is called a {\em monodromy operator}, 
if there exists a 
family $\X \rightarrow B$ (which may depend on $g$) 
of irreducible holomorphic symplectic manifolds, having $X$ as a fiber
over a point $b_0\in B$, 
and such that $g$ belongs to the image of $\pi_1(B,b_0)$ under
the monodromy representation. 
The {\em monodromy group} $Mon(X)$ of $X$ is the subgroup 
of $GL(H^*(X,\Integers))$ generated by all the monodromy operators. 
}
\end{defi}

Let $X$ be an irreducible holomorphic symplectic manifold
deformation equivalent to $S^{[n]}$.
The second cohomology
$H^2(X,\Integers)$ admits the symmetric bilinear 
Beauville-Bogomolov pairing \cite{beauville-varieties-with-zero-c-1}.
$H^2(S^{[n]},\Integers)$, $n\geq 2$, for example, is the orthogonal direct sum 
$H^2(S,\Integers)\oplus \Integers\delta$, where $\delta$ is half the class of
the big diagonal, and the Beauville-Bogomolov pairing 
restricts to the standard one on $H^2(S,\Integers)$, while 
$(\delta,\delta)=2-2n$. Other examples are described in 
Theorem \ref{thm-irreducibility}.

Given an element $u$ of $H^2(X,\Integers)$ with
$(u,u)$ equal $2$ or $-2$, let $\rho_u$ be the isometry of
$H^2(X,\Integers)$ given by
$\rho_u(w)=\frac{-2}{(u,u)}w+(w,u)u$. Then $\rho_u$ is the
reflection in $u$, when $(u,u)=-2$, and $-\rho_u$ is the reflection in $u$,
when $(u,u)=2$. 
Set 
\begin{equation}
\label{eq-W}
\W \ \ \ := \ \ \ \langle
\rho_u \ : \ u \in H^2(X,\Integers) \ \
\mbox{and} \ \ (u,u) = \pm2
\rangle
\end{equation} 
to be the subgroup of $O(H^2(X,\Integers))$ generated by the 
elements $\rho_u$. Then $\W$ is a normal subgroup of finite index
in $O(H^2(X,\Integers))$. 

Let $Mon^2$ be the image of $Mon(X)$ in $GL[H^2(X,\Integers)]$. 
The inclusion $\W\subset Mon^2$ was proven in 
\cite{markman-monodromy-I}.
We prove the reversed inclusion and obtain:

\begin{thm}
\label{thm-introduction-Mon-2-is-W}
\ \ \ $Mon^2 \ \ \ = \ \ \ \W$.
\end{thm}

The main ideas of the proof are summarized in section
\ref{sec-introduction-the-Mukai-lattice-as-a-quotient}, where
the proof of 
Theorem \ref{thm-introduction-Mon-2-is-W} is reduced to that of 
Theorem \ref{thm-introduction-invariant-Q-4}. 

Let $K^2\subset Mon(X)$ be the subgroup, which acts trivially on 
$H^2(X,\Integers)$. We get the short exact sequence
\[
0 \rightarrow K^2 \rightarrow Mon(X) \rightarrow \W \rightarrow 0.
\]
The sequence admits a $Mon(X)$-equivariant splitting 
\[
\nu \ : \  \W \ \ \ \rightarrow \ \ \ Mon(X)
\] 
(Corollary 1.8 in \cite{markman-monodromy-I}).
When $X$ is a smooth and compact moduli space 
of stable sheaves on a K3 surface $S$, the monodromy operator
$\nu(g)$ is described explicitly in terms of the Chern classes of a 
universal sheaf over $S\times X$ (see 
equation (\ref{eq-monrep-g}) and Lemma \ref{lemma-mu-is-an-isomorphism}). 
The following statement was conjectured in \cite{markman-monodromy-I}.

\begin{conj} 
\label{conj-Mon-isomorphic-to-Mon-2}
Any monodromy operator, which acts as the identity on $H^2(X,\Integers)$,
is the identity operator.
\end{conj}

Geometric implications of the conjecture are described in section 
\ref{sec-reduction-of-monodromy-conjecture-to-generators-conjecture}. 
The conjecture is related to a Torelli question in 
Proposition \ref{prop-torelli-implies-K-2-is-trivial}.
We prove here the following three results.

\begin{prop}
\label{prop-mon-is-K2-times-W}
$Mon(X)$ is the direct product of $\nu(\W)$ and the subgroup 
$K^2\subset Mon(X)$. 
$K^2$ is a finite abelian group.
\end{prop}

$K^2$ is thus contained in the center of $Mon(X)$. $K^2$ consists
of ring and $Mon(X)$-module automorphisms of $H^*(X,\Integers)$.
The Proposition follows immediately from
Theorem \ref{thm-introduction-Mon-2-is-W} above and  Lemma
4.7 in \cite{markman-monodromy-I} 
(the proof of which relies also on results of Verbitsky \cite{verbitsky}).
Proposition \ref{prop-mon-is-K2-times-W} is proven in section 
\ref{sec-proof-of-thm-arithmetic-invariant}. 

\begin{thm}
\label{thm-monodromy-operator-is-determined-by-its-weight-2-action}
Any monodromy operator, which acts as the identity on $H^2(X,\Integers)$,
acts as the identity also on
$H^k(X,\Integers)$, $0\leq k \leq \frac{\dim_\RealNumbers(X)}{4}+2$. 
\end{thm}

The Theorem is proven in section \ref{sec-extension-classes}.
Following is the strategy of the proof.
$H^{2i}(X,\RationalNumbers)$ decomposes as a direct sum
of irreducible $Mon(X)$-submodules. $H^{2i}(X,\Integers)$ does not.
Let $(A_{2i-2})^{2i}$ be the subgroup of $H^{2i}(X,\Integers)$
of classes, which are polynomials in integral classes of degree $\leq 2i-2$.
$H^{2i}(X,\Integers)$ is thus an extension of $Mon(X)$-modules:\\
$0\rightarrow (A_{2i-2})^{2i}\rightarrow H^{2i}(X,\Integers)\rightarrow 
H^{2i}(X,\Integers)/(A_{2i-2})^{2i}\rightarrow 0$. 
We calculate the order of these extension classes.
Their non-triviality imposes constraints on ring automorphisms of the
$Mon(X)$-module $H^*(X,\RationalNumbers)$ to preserve the integral 
structure, and hence to belong to $K^2$.

Set $n:=\frac{\dim_\RealNumbers(X)}{4}$. 
Results of M. Lehn and C. Sorger state, that
$H^*(X,\RationalNumbers)$ is generated by $H^k(X,\RationalNumbers)$, 
$2\leq k \leq 2\lfloor\frac{n}{2}\rfloor+4$. 
See Theorem \ref{thm-vanishing-of-varphi-k} below
for a more detailed statement. 
``Half'' the cases of Conjecture \ref{conj-Mon-isomorphic-to-Mon-2}
follow from their results and 
Theorem \ref{thm-monodromy-operator-is-determined-by-its-weight-2-action}, 
and in the remaining cases the order of $K^2$ is at most $2$.
\hide{
Conjecture 
\ref{conj-Mon-isomorphic-to-Mon-2} follows from 
Theorem 
\ref{thm-monodromy-operator-is-determined-by-its-weight-2-action}
and Conjecture \ref{conj-vanishing-of-varphi-k}, 
if $n\equiv 0$ or $n\equiv 1$ modulo $4$
(Proposition 
\ref{prop-reduction-of-monodromy-conjecture-to-generators-conjecture}).
If $n\equiv 2$ or $n\equiv 3$ modulo $4$, then Theorem 
\ref{thm-monodromy-operator-is-determined-by-its-weight-2-action} 
and Conjecture
\ref{conj-vanishing-of-varphi-k} imply that $K^2$
has order at most $2$ (Proposition 
\ref{prop-reduction-of-monodromy-conjecture-to-generators-conjecture}). 
The case $n=3$ of Conjecture
\ref{conj-Mon-isomorphic-to-Mon-2} follows from
Theorem \ref{thm-monodromy-operator-is-determined-by-its-weight-2-action}
independently of Conjecture \ref{conj-vanishing-of-varphi-k}.
Conjecture \ref{conj-Mon-isomorphic-to-Mon-2}
is clear for $n=2$, since $H^*(X,\RationalNumbers)$
is generated by $H^2$ (Lemma \ref{lemma-10-in-markman-diagonal}). 
}

\begin{cor}
\label{cor-reduction-of-monodromy-conjecture-to-generators-conjecture}
The restriction homomorphism 
$Mon(S^{[n]})\rightarrow Mon^2(S^{[n]})$ is injective
(and Conjecture \ref{conj-Mon-isomorphic-to-Mon-2} holds for 
$S^{[n]}$), if $n\equiv 0$ or $n\equiv 1$  modulo $4$, or if $n\leq 3$. 
If $n\geq 6$ and $n\equiv 2$ or  $n\equiv 3$ modulo 4, then the kernel
$K^2$ of the restriction homomorphism has order $\leq 2$.
\end{cor}

The Corollary is proven in section 
\ref{sec-reduction-of-monodromy-conjecture-to-generators-conjecture}.
%
\subsection{Monodromy and the Generic Torelli question}
A negative answer to a 
common formulation of the Generic Torelli question
follows from 
Theorem \ref{thm-introduction-Mon-2-is-W}
(Corollary \ref{cor-introduction-counter-example-to-Torelli}).
Let $X$ be an irreducible holomorphic symplectic manifold, 
deformation equivalent to $S^{[n]}$,
with $n\geq 1$. 
A {\em marking} for $X$ is a choice of an isometry
$\phi: H^2(X,\Integers)\rightarrow \Lambda$ with a fixed lattice $\Lambda$.
The {\em period}, of the marked manifold $(X,\phi)$, is the
line $\phi[H^{2,0}(X)]$ considered as a point in the projective space
$\PP[\Lambda\otimes\ComplexNumbers]$. The period lies in the period domain 
\begin{equation}
\label{eq-period-domain}
\Omega \ := \ \{
\ell \ : \ (\ell,\ell)=0 \ \ \ \mbox{and} \ \ \ 
(\ell,\bar{\ell}) > 0
\}.
\end{equation}
$\Omega$ is an open subset, in the classical topology, of the quadric in 
$\PP[\Lambda\otimes\ComplexNumbers]$ of isotropic lines. 
The arithmetic group 
$\PP{O}\Lambda:=O\Lambda/(\pm 1)$ naturally acts on $\Omega$,
and the stabilizer of a generic point is trivial. 
$\PP{O}\Lambda$-orbits
in $\Omega$ are in bijection with {\em Hodge-isometry classes} of integral 
weight $2$ Hodge structures. The analogue
$\W(\Lambda)$, of the group given in (\ref{eq-W}),
maps injectively into $\PP{O}\Lambda$ and we
denote its image by $\PP\W(\Lambda)$. 
The index $[\PP\W(\Lambda):\PP{O}\Lambda]$  
of $\PP\W(\Lambda)$ in $\PP{O}\Lambda$ is calculated as follows.
The {\em Euler number} $\rho(m)$ of an integer $m\geq 1$,    
is the number $r$ of distinct primes $p_i>1$, in the prime decomposition
$m=p_1^{e_1}\cdots p_{r}^{e_r}$.
Then $[\PP\W(\Lambda):\PP{O}\Lambda]$ is equal to 
$2^{\rho(n-1)-1}$, for $n\geq 3$, and $\PP\W(\Lambda)=\PP{O}\Lambda$
for $1\leq n\leq 6$ (Lemma \ref{lemma-index-of-smaller-symmetry-group}).  

\begin{cor}
\label{cor-introduction-counter-example-to-Torelli}
Let $h\in \Omega$ be a point, whose stabilizer 
in $\PP{O}\Lambda$ is contained in $\PP\W(\Lambda)$.
There are at least $[\PP\W(\Lambda):\PP{O}\Lambda]$ isomorphism classes of
irreducible holomorphic symplectic manifolds,
deformation equivalent to $S^{[n]}$, with weight $2$ Hodge structures
Hodge-isometric to $h$.
\end{cor}

The short proof of Corollary 
\ref{cor-introduction-counter-example-to-Torelli} is included below.
A stronger version of the corollary
is proven in section \ref{sec-period-maps}  (Theorem 
\ref{thm-non-bimeromorphic-classes}).
There is a (non-Hausdorff) moduli space ${\frak M}_\Lambda$ of marked 
irreducible holomorphic symplectic manifolds, 
with a second integral cohomology 
lattice isometric to $\Lambda$ \cite{huybrects-basic-results}. 
The period map 
\begin{eqnarray}
\label{eq-period-map}
P \ : \ {\frak M}_\Lambda & \longrightarrow & \Omega,
\\
\nonumber
(X,\phi) & \mapsto & \phi[H^{2,0}(X)]
\end{eqnarray}
is a local isomorphism, by the Local Torelli Theorem
\cite{beauville-varieties-with-zero-c-1}.
The Surjectivity Theorem states, that the restriction of the period map
to every connected component of ${\frak M}_\Lambda$ is surjective 
\cite{huybrects-basic-results}.

\medskip
{\bf Proof of Corollary \ref{cor-introduction-counter-example-to-Torelli}:}
The group $O\Lambda$ acts on ${\frak M}_\Lambda$, by changing the marking,
and the period map is $O\Lambda$-equivariant. 
Let $(S^{[n]},\phi_0)$ be a marked Hilbert scheme in a connected
component ${\frak M}^0_\Lambda$ of ${\frak M}_\Lambda$.
The marking $\phi_0$ conjugates $Mon^2(S^{[n]})$ to the largest subgroup 
$Mon^2(\Lambda)$ of
$O\Lambda$, leaving invariant the connected component ${\frak M}^0_\Lambda$. 
The inclusion $Mon^2(\Lambda)\subset Mon^2(S^{[n]})^{\phi_0}$
follows from the definition of $Mon^2(S^{[n]})$ and the reverse inclusion
follows from the universal property of the moduli space ${\frak M}_\Lambda$. 
Note that {\em the moduli space ${\frak M}_\Lambda$ has at least 
$[Mon^2(\Lambda):O\Lambda]$
connected components}.

The quotient 
${\frak M}^0_\Lambda/Mon^2(\Lambda)$ is the moduli space of 
isomorphism classes of irreducible holomorphic symplectic manifolds, 
deformation equivalent to $S^{[n]}$. 
The period map induces a map from 
${\frak M}^0_\Lambda/Mon^2(\Lambda)$ to the moduli space
$\Omega/\PP{O}\Lambda$ of Hodge isometry classes, which factors through
a surjective map onto the intermediate quotient
\[
{\frak M}^0_\Lambda/Mon^2(\Lambda) \ \ \longrightarrow \ \ 
\Omega/Mon^2(\Lambda). 
\]
Now $\Omega/\PP{O}\Lambda$ is the quotient of
$\Omega/Mon^2(\Lambda)$ and the cardinality of the fiber, 
containing the orbit $Mon^2(\Lambda)\cdot h$, is equal to 
the index of the image of $Mon^2(\Lambda)$ in $\PP{O}\Lambda$, by our
assumption on $h$. 
Theorem \ref{thm-introduction-Mon-2-is-W}
implies the equality $Mon^2(\Lambda)=\W(\Lambda)$, and
so the relevant index
is $[\PP\W:\PP{O}\Lambda]$.
\EndProof

Another common formulation of the
Generic Torelli question is open and recalled next. 
Let ${\frak M}^0_\Lambda$ be a connected component of ${\frak M}_\Lambda$, 
containing a marked Hilbert scheme of a K3 surface. 
Let $\Omega' \subset \Omega$ be the set of points, over which the fiber of
$P$ intersects the component ${\frak M}^0_\Lambda$ in a single point. 
It is not clear that $\Omega'$ is non-empty.
The (open) Generic Torelli question is:

\begin{question}
\label{question-generic-torelli}
Is the complement $\Omega\setminus \Omega'$ contained in the union of
some countable collection of proper analytic subsets of $\Omega$?
\end{question}

The Torelli Theorem, for $K3$ surfaces, provides an affirmative answer 
in case $n=1$ \cite{burns-rapoport}. 
The answer to Question \ref{question-generic-torelli} is independent of 
the image of the restriction homomorphism
$Mon(S^{[n]})\rightarrow OH^2(S^{[n]},\Integers)$, but the answer is 
related to its kernel $K^2$.

\begin{prop}
\label{prop-torelli-implies-K-2-is-trivial}
If the answer to Question \ref{question-generic-torelli} is affirmative,
then the restriction homomorphism
$Mon(S^{[n]})\rightarrow Mon^2(S^{[n]})$ is injective.
\end{prop}

Proposition \ref{prop-torelli-implies-K-2-is-trivial} will not be used in 
this paper, but it serves as a motivation for the results in section
\ref{sec-extension-classes}. 
\journal{
The Proposition is proven in the eprint version of this paper
using results of Beauville and Verbitsky \cite{markman-eprint-version}.
}
\preprint{
The following Proposition, due to Beauville and Verbitsky, 
is  key to the proof of Proposition 
\ref{prop-torelli-implies-K-2-is-trivial}.

\smallskip
\noindent
{\bf Proposition}
{\em
Let $X$ be an irreducible holomorphic symplectic manifold 
deformation equivalent to 
the Hilbert scheme $S^{[n]}$, $n\geq 1$, of a $K3$ surface $S$.
A holomorphic automorphism of $X$, 
which acts as the identity on $H^2(X,\Integers)$,  is the identity
automorphism.
}

\noindent
{\bf Proof:}
The case $X=S^{[n]}$ was proven by Beauville 
(\cite{beauville-automorphisms} Proposition 10).
The statement follows, for deformations of $S^{[n]}$,
by Corollary 6.9 in \cite{kaledin-verbitsky}.
\EndProof

{\bf Proof of Proposition \ref{prop-torelli-implies-K-2-is-trivial}:}
Let ${\frak M}^0_\Lambda$ be a component of the moduli space of marked 
irreducible holomorphic symplectic manifolds 
deformation equivalent to $S^{[n]}$. 
There exists a universal family $\pi:\X\rightarrow {\frak M}^0_\Lambda$, by 
the above Proposition. 
Let $(X,\phi)$ be a marked pair in ${\frak M}^0_\Lambda$ and
$Mon(X,\pi)$ the subgroup of $Mon(X)$ of monodromy operators 
of the local system $\oplus_iR^i_{\pi_*}\Integers$.
The subgroup $Mon(X,\pi)$ is precisely the kernel of 
$Mon(X)\rightarrow Mon^2(X)$. 
It suffices to show that the local systems $R^i_{\pi_*}\Integers$
are trivial. 
An affirmative answer to Question \ref{question-generic-torelli}
implies that any two points, in the same fiber of the
period map $P$, must be non-separated. It follows
that the local system $R^i_{\pi_*}\Integers$ over ${\frak M}^0_\Lambda$
must be the pull-back of some local system over the image of $P$ in
$\Omega$. The image is the whole of $\Omega$, by Huybrechts Surjectivity 
Theorem \cite{huybrects-basic-results}.
The local systems $R^i_{\pi_*}\Integers$ must thus be trivial, since
$\Omega$ is simply connected. 
\EndProof
}

The generalized Kummer varieties constitute another sequence of 
irreducible holomorphic symplectic manifolds 
\cite{beauville-varieties-with-zero-c-1}.
Namikawa found a counter example to the Torelli question for 
generalized Kummer fourfolds \cite{namikawa}.
The Hodge-isometry class of their 
weight $2$ Hodge structure does not determine their bimeromorphic class.
The kernel of the homomorphism $Mon(X)\rightarrow Mon^2(X)$,
for a generalized Kummer $X$ of dimension $\geq 4$, is non-trivial
and non-abelian (\cite{markman-monodromy-I}, section 4.2). 
Question \ref{question-generic-torelli} is open for the generalized Kummer 
varieties as well.

%
\subsection{The Mukai lattice as an invariant of a hyper-K\"{a}hler variety}
\label{sec-introduction-the-Mukai-lattice-as-a-quotient}
The new information in Theorem \ref{thm-introduction-Mon-2-is-W} 
is the inclusion $Mon^2\subset \W$. 
The proof of this inclusion identifies an 
obstruction for an isometry of $H^2(X,\Integers)$
to extend to a monodromy operator of the full cohomology ring. 
The obstruction arises from  
an invariant $\{Q^4(X,\Integers),\bar{c}_2(X)\}$
of the Hodge structure of these varieties, which is more refined 
than $H^2(X,\Integers)$, 
whenever $(\dim_\ComplexNumbers{X}-2)/2$ is not a prime power. 
The integral Hodge structure 
$Q^4(X,\Integers)$  
is the quotient of $H^4(X,\Integers)$ 
by $\Sym^2H^2(X,\Integers)$ and
$\bar{c}_2(X)$ is the image in $Q^4(X,\Integers)$ of the second Chern class
$c_2(TX)$. 

We consider also analogous quotients of $H^d(X,\Integers)$, 
$d\geq 4$, in order to prove Theorem 
\ref{thm-monodromy-operator-is-determined-by-its-weight-2-action}. 
Let $A_d \subset H^*(X,\Integers)$, $d\geq 0$,  be the graded subring of 
$H^*(X,\Integers)$, generated by classed in $H^i(X,\Integers)$, $i\leq d$. 
Set $(A_d)^k:=A_d\cap H^k(X,\Integers)$. 
The odd integral cohomology groups of $X$ vanish 
(\cite{markman-original-archive-version} or Theorem 
\ref{thm-integral-generators} below). 
We set $Q^2(X,\Integers) := H^2(X,\Integers)$. 
For an integer $i\geq 2$, we 
define
$Q^{2i}(X,\Integers)$ as the quotient in the short exact sequence
\begin{equation}
\label{eq-extension-of-Q-2i-by-A-2i-2}
\label{eq-C-d}
0 \rightarrow 
(A_{2i-2})^{2i}
\rightarrow 
H^{2i}(X,\Integers) \rightarrow Q^{2i}(X,\Integers)
 \rightarrow 0.
\end{equation}

Let $U$ be the hyperbolic plane, 
whose bilinear form is given by the matrix
{\scriptsize
$\left[
\begin{array}{cc}
0 & 1
\\ 
1 & 0
\end{array}
\right].
$
}
Denote by $(-E_8)$ the $E_8$ lattice, 
with a negative definite bilinear form. 
A non-zero element of a free abelian group is said to be 
{\em primitive}, if it is not a multiple of another 
element by an integer larger than $1$.

\begin{thm}
\label{thm-introduction-invariant-Q-4}
Let $X$ be an irreducible holomorphic symplectic manifold 
deformation equivalent to
the Hilbert scheme $S^{[n]}$, $n\geq 4$. 

\noindent
a) Let $i$ be an even integer in the range 
$2\leq i \leq \frac{1}{4}\dim_\ComplexNumbers X$.
\begin{enumerate}
\item
\label{thm-item-second-chern-class}
$Q^{2i}(X,\Integers)$ is torsion free.
Let $\bar{c}_{i}(X)\in Q^{2i}(X,\Integers)$ be 
the projection of the Chern class
$c_{i}(TX)$ of the tangent bundle. Then $\frac{1}{2}\bar{c}_{i}(X)$ 
is a  non-zero, integral, and primitive element.
\item
\label{thm-item-unimodular-pairing}
There exists a unique monodromy-invariant even unimodular 
symmetric bilinear form on $Q^{2i}(X,\Integers)$, satisfying 
\begin{equation}
\label{eq-chern-classes-encode-dimension}
\left(\frac{1}{2}\bar{c}_{i}(X),\frac{1}{2}\bar{c}_{i}(X)\right)
\ \ \ = \ \ \ \dim_\ComplexNumbers(X)-2.
\end{equation}
\item
\label{thm-intro-item-Q-4-is-the-mukai-lattice}
The lattice $Q^{2i}(X,\Integers)$ 
is isometric to the orthogonal direct sum 
$U^{\oplus 4}\oplus (-E_8)^{\oplus 2}$.
\item
\label{thm-item-pair-of-embeddings}
There exists a monodromy invariant pair of primitive lattice embeddings 
\[
e, -e \ : \ H^2(X,\Integers) \ \ \ \LongIsomRightArrow \ \ \
\bar{c}_{i}(X)^\perp \ \subset \ Q^{2i}(X,\Integers).
\]
The embedding $e$ induces an isometry, compatible with the Hodge 
structures, between $H^2(X,\Integers)$ and the corank $1$ sublattice
orthogonal to $\bar{c}_{i}(X)$. 
The character ${\rm span}_\Integers\{e\}$ of $Mon(X)$
is  non-trivial. 
\end{enumerate}

b) When $i$ is an odd integer in the above range, the above statements hold, 
with the following modifications. 
In part \ref{thm-item-second-chern-class} the class 
$\bar{c}_{i}(X)$ vanishes, however, there is a 
rank $1$ sublattice $Q^{2i}(X,\Integers)'$, 
which is a monodromy subrepresentation of $Q^{2i}(X,\Integers)$. 
Part \ref{thm-item-unimodular-pairing} holds, with
$\bar{c}_{i}(X)/2$ replaced by an integral generator of
$Q^{2i}(X,\Integers)'$.
Part \ref{thm-item-pair-of-embeddings} holds, with
$\bar{c}_{i}(X)$ replaced by $Q^{2i}(X,\Integers)'$,
and excluding the claim that ${\rm span}\{e\}$ in a non-trivial character. 
\end{thm}



Theorem \ref{thm-introduction-invariant-Q-4} is proven in section
\ref{sec-proof-of-thm-arithmetic-invariant}. 
We denote the pair $\{Q^4(X,\Integers),\bar{c}_2(X)\}$ by
\begin{equation}
\label{eq-short-notation}
\{Q^4,\bar{c}_2\}(X)
\end{equation}
for short. Two pairs $\{Q^4,\bar{c}_2\}(X)$ 
and $\{Q^4,\bar{c}_2\}(Y)$ are {\em Hodge-isometric},
if there is an isometry between $Q^4(X,\Integers)$ and $Q^4(Y,\Integers)$,
compatible with the Hodge structures, which sends $\bar{c}_2(X)$ to
$\bar{c}_2(Y)$. 
Hodge-isometry classes of such pairs are in one-to-one
correspondence with $\PP\W(\Lambda)$-orbits in the period domain $\Omega$,
by Lemma \ref{lemma-equivalent-data}. The invariant $\{Q^4,\bar{c}_2\}(X)$
seems to capture all the cohomological invariants of $X$,
by the following Theorem.

\begin{defi}
\label{def-Mon-X-Y}
{\rm
Let $X_1$ and $X_2$ be two deformation equivalent 
irreducible holomorphic symplectic manifolds. 
An isomorphism 
$f: H^*(X_1,\Integers)\rightarrow H^*(X_2,\Integers)$ is called a 
{\em parallel-transport operator}, if 
$X_1$ and $X_2$ are two fibers of a family 
$\pi:\X\rightarrow B$ as in 
Definition \ref{def-monodromy}, over two points $b_1$, $b_2$ in $B$,
and $f$ is associated to a path in $B$ from $b_1$ to $b_2$ via the
local system $R^*\pi_*\Integers$.
}
\end{defi}

A parallel-transport operator is in particular a graded ring isomorphism, 
which maps the Chern class $c_i(X_1)$ to $c_i(X_2)$, for all $i>0$.

\begin{thm}
\label{thm-invariant-Q-bar-c-captures-all}
Let $X_1$, $X_2$ be two irreducible holomorphic symplectic manifolds
deformation equivalent to $S^{[n]}$.
Assume that one of the following two conditions holds:

(i) $n\leq 3$ and there exists a Hodge-isometry 
$g:H^2(X_1,\Integers)\rightarrow H^2(X_2,\Integers)$.

(ii) $n\geq 4$ and there exists a Hodge-isometry 
$g:\{Q^4,\bar{c}_2\}(X_1)\rightarrow \{Q^4,\bar{c}_2\}(X_2)$.

\noindent
Then there exists a parallel-transport operator 
$f:H^*(X_1,\Integers)\rightarrow H^*(X_2,\Integers)$, inducing $g$, 
which is compatible with Hodge-structures.
\end{thm}

The Theorem is proven in section \ref{sec-symmetry-groups}.
Conjecture \ref{conj-Mon-isomorphic-to-Mon-2}
implies that such $f$ is uniquely determined by $g$. 
An explicit formula for $f$ is given in 
\cite{markman-monodromy-I} equation (1.22),
whenever each $X_i$ is a moduli space of sheaves on some $K3$ surface $S_i$,
$i=1,2$.
The Hodge-isometry classes of both $H^2(X,\Integers)$ and
$\{Q^4,\bar{c}_2\}(X)$ are 
constant throughout the bimeromorphic class of $X$.
This follows from Huybrechts' result, 
that bimeromorphic irreducible holomorphic symplectic manifolds correspond to
non-separated points in their 
moduli space \cite{huybrechts-kahler-cone}.

{\bf Reduction of the proof of Theorem \ref{thm-introduction-Mon-2-is-W}
to that of Theorem \ref{thm-introduction-invariant-Q-4}:
}
The inclusion $\W\subset Mon^2$ was proven in \cite{markman-monodromy-I}.
We prove the reverse inclusion.
Set $\widetilde{\Lambda}:=U^{\oplus 4}\oplus (-E_8)^{\oplus 2}$.
Choose any primitive embedding 
$\iota: H^2(X,\Integers)\hookrightarrow \widetilde{\Lambda}$.
The positive cone 
\begin{equation}
\label{eq-positive-cone}
\C_+ \ \ := \ \ \{\lambda\in  H^2(X,\RealNumbers) \ : \ 
(\lambda,\lambda)>0\}
\end{equation}
in $H^2(X,\RealNumbers)$ is homotopic to the two-sphere 
and comes with a distinguished generator of $H^2(\C_+,\Integers)$ 
(Definition \ref{def-distinguished-orientation}).
The subgroup $\W$ of $OH^2(X,\Integers)$ is characterized by the following two
properties, by Lemma 4.10 in \cite{markman-monodromy-I}
(see also Lemma \ref{lemma-on-residual-orthogonal-group} below). 

1) For every $g\in \W$, there exists 
$\tilde{g}\in O(\widetilde{\Lambda})$, such that 
$\tilde{g}\circ \iota=\iota\circ g$ ($g$ extends to 
an isometry of $\widetilde{\Lambda}$).

2) $g$ acts as the identity on $H^2(\C_+,\Integers)$.

\noindent
Property 2 is automatically satisfied by any element of $Mon^2$.
If $\dim_{\ComplexNumbers}X\geq 8$, then 
Theorem \ref{thm-introduction-invariant-Q-4} verifies 
Property 1 for any element of $Mon^2$, where $\iota$ is the embedding
$e$ in  Theorem 
\ref{thm-introduction-invariant-Q-4} part \ref{thm-item-pair-of-embeddings}. 
If $\dim_{\ComplexNumbers}X$ is $4$ or $6$, property 1 holds for any 
isometry, by Lemma 4.10 in
\cite{markman-monodromy-I} (see also Lemma 
\ref{lemma-on-residual-orthogonal-group} below).
The inclusion $Mon^2\subset \W$ follows.
\EndProof

%
\subsubsection*{The invariant $\{Q^4,\bar{c}_2\}(\M)$ 
for a moduli space $\M$ of sheaves}
The above lattice $Q^4(X,\Integers)$ is isometric to the Mukai lattice.
We recall its definition.
Let $S$ be a $K3$ surface and $K(S)$ its
Grothendieck $K$-ring generated by topological complex vector bundles
\cite{atiyah-book}. The 
Euler characteristic $\chi: K(S) \rightarrow \Integers$ is given by
\[
\chi(v) = 2\rank(v) + \int_S\left[c_1(v)^2/2-c_2(v)\right].
\]
$K(S)$ is a free abelian group
of rank $24$. Any class $x\in K(S)$ is the difference 
$[E]-[F]$ of classes of complex vector bundles, and
we denote by $x^\vee$ the class $[E^\vee]-[F^\vee]$ obtained from  the
dual vector bundles.
The bilinear form on $K(S)$, given by 
\begin{equation}
\label{eq-mukai-pairing-on-K-top}
(x,y) \ \ \ := \ \ \ -\chi(x^\vee\cup y),
\end{equation}
is called the {\em Mukai pairing}. The pairing is symmetric and unimodular,
and the resulting lattice is called the {\em Mukai lattice}. 
It is isometric to $U^{\oplus 4}\oplus (-E_8)^{\oplus 2}$
\cite{mukai-hodge}.
We define a polarized weight $2$ 
Hodge structure on $K(S)_\ComplexNumbers:=
K(S)\otimes_\Integers\ComplexNumbers$,
by setting $K^{1,1}(S)$ 
to be the subspace of $K(S)_\ComplexNumbers$
of classes $v$ with $c_1(v)$ of type $(1,1)$,
and $K^{2,0}(S)$ the subspace of $K^{1,1}(S)^\perp$
of classes $v$ with $c_1(v)$ of type $(2,0)$.

We describe next 
the invariant $\{Q^4,\bar{c}_2\}(X)$, when $X$ is
$S^{[n]}$ or a more general moduli space of sheaves on the 
projective K3 surface $S$.

\begin{defi}
\label{def-effective-class}
{\rm
A non-zero class $v\in K^{1,1}(S)$ will be called {\em effective},
if $(v,v)\geq -2$, 
$\rank(v)\geq 0$, and the following conditions hold.
If $\rank(v)=0$, then $c_1(v)$ is the
class of an effective (or trivial) divisor on $S$.
If both $\rank(v)$ and $c_1(v)$ vanish, then $\chi(v)>0$.
}
\end{defi}

Let $H$ be an ample line bundle on $S$
and $v$ a primitive and effective class in $K(S)$ of type $(1,1)$. 
If both $\rank(v)$ and $\chi(v)$ vanish, assume further that $c_1(v)$ 
generates the Neron-Severi group of $S$.
Denote by $\M_H(v)$ the moduli space 
of $H$-stable sheaves on $S$ with class $v$. 
Stability is in the sense of Gieseker, Maruyama, and Simpson
(Definition \ref{def-stability}). 
For a generic $H$, in a sense to be made precise in Definition
\ref{def-v-suitable}, $\M_H(v)$ is a non-empty, 
smooth, projective, connected, holomorphically symplectic
variety, which is deformation
equivalent to $S^{[n]}$, where $n=(v,v)/2+1$
(see Theorem \ref{thm-irreducibility} below for the relevant citations).
We assume, throughout the paper, that $\M_H(v)$ is such.

Set $\M:=\M_H(v)$.
There is a class $[\E]$ in $K(S\times \M)$, unique up to 
tensorization by the class of a topological line-bundle on
$\M$, and defined below in (\ref{eq-e-v}). 
The class $[\E]$ is associated to a (possibly twisted) 
universal sheaf $\E$. 
Denote by $f_i$, $i=1,2$, the projection from $S\times \M$
onto the $i$-th factor. Let 
$f_1^!:K(S)\rightarrow K(S\times \M)$ be the pullback 
homomorphism and
$f_{2_!}:K(S\times \M)\rightarrow K(\M)$
the Gysin homomorphism 
(\cite{bfm},\cite{karoubi} Proposition IV.5.24). 
The class $[\E]$ induces the homomorphism
\begin{equation}
\label{eq-u}
u \ : \ K(S) \ \ \ \longrightarrow \ \ \ K(\M),
\end{equation}
given by $u(x)=f_{2_!}(f_1^!(-x^\vee)\cup [\E])$, 
and a polynomial map
\begin{equation}
\label{eq-tilde-varphi-2i}
\tilde{\varphi}^{2i} \ : \ 
K(S) \ \ \  \longrightarrow  \ \ \ H^{2i}(\M,\Integers),
\end{equation}
given by $\tilde{\varphi}^{2i}(y)=c_i(u(y))$.
Let 
\begin{equation}
\label{eq-varphi-k}
\varphi^{2i} \ : \ 
K(S)  \ \ \ \longrightarrow \ \ \ Q^{2i}(\M,\Integers)
\end{equation}
be the composition of $\tilde{\varphi}^{2i}$ with the projection 
$H^{2i}(\M,\Integers)\rightarrow Q^{2i}(\M,\Integers)$.
Then $\varphi^{2i}$ is a linear homomorphism
(Proposition \ref{prop-varphi-k-is-independent-of-universal-sheaf}). 
The notation $a=\pm b$ means, $a=b$ or $a=-b$. 

\begin{thm}
\label{thm-introduction-arithmetic-invariant-of-moduli}
a) Let $i$ be an even integer in the range 
$2\leq i \leq \frac{1}{4}\dim_\ComplexNumbers \M$.
The two quadruples
\[
\begin{array}{c}
\left\{
Q^{2i}(\M,\Integers),H^2(\M,\Integers), \pm e, \bar{c}_{i}(\M)
\right\},
\\
\left\{K(S),v^\perp, \pm \iota, 2v\right\}
\end{array}
\]
are isomorphic,
where $K(S)$ is endowed with the Mukai pairing 
(\ref{eq-mukai-pairing-on-K-top}), and 
$\iota:v^\perp\hookrightarrow K(S)$ is the inclusion. 
The isomorphism is understood in the following sense. 
\begin{eqnarray}
\label{eq-varphi-4-maps-v-to-a-multiple-of-second-Chern-class}
\bar{c}_{i}(\M) & = & \varphi^{2i}(2v) \ \ \ \mbox{and}
\\
\label{eq-e-is-conjugated-to-iota}
\varphi^{2i}\circ \iota & = & \pm e\circ \varphi^2\circ \iota,
\end{eqnarray}
where both $\varphi^{2i}$, which is given in (\ref{eq-varphi-k}), 
and $\varphi^2\circ\iota$ are Hodge-isometries. 

b) When $i$ is an odd integer in the above range, 
the above statement  holds with 
$\bar{c}_{i}(\M)$ replaced by the rank $1$ subrepresentation 
$Q^{2i}(\M,\Integers)'$,
as in Theorem \ref{thm-introduction-invariant-Q-4} part b, and
$2v$ replaced by ${\rm span}_\Integers\{v\}$. 
\end{thm}

The theorem is proven in section \ref{sec-proof-of-thm-arithmetic-invariant}.
Theorem \ref{thm-introduction-arithmetic-invariant-of-moduli}
is a refinement of the well-known statement in 
Theorem \ref{thm-irreducibility}, that
$\varphi^2\circ \iota:v^\perp\rightarrow H^2(\M,\Integers)$ 
is a Hodge-isometry.

Three invariants are associated to an irreducible holomorphic 
symplectic manifold $X$, 
deformation equivalent to $S^{[n]}$, $n\geq 4$; the Hodge-isometry classes of
1) the transcendental sub-lattice $\Theta(X)$ of $H^2(X,\Integers)$, 
2) $H^2(X,\Integers)$, and 
3) $\{Q^4,\bar{c}_2\}(X)$. Each determines the previous one.
Theorem \ref{thm-introduction-arithmetic-invariant-of-moduli} 
implies that $\Theta(X)$ is Hodge-isometric to $\Theta(S')$,
for some projective $K3$ surface $S'$, if and only if 
$\{Q^4,\bar{c}_2\}(X)$ is Hodge-isometric to $\{Q^4,\bar{c}_2\}(\M_H(v))$,
for some primitive and effective class $v\in K(S')$
(Lemma \ref{lemma-realization-as-a-moduli-space}). 

Similarly, $H^2(X,\Integers)$ is Hodge isometric to 
$H^2(S^{[n]},\Integers)$, if and only if 
$\{Q^4,\bar{c}_2\}(X)$ is Hodge-isometric to $\{Q^4,\bar{c}_2\}(\M_H(v))$,
for some primitive and effective class $v\in K(S)$, with 
$c_1(v)=0$ (Proposition \ref{thm-non-birational-moduli-spaces}). 
For $n\geq 7$ and a suitably chosen 
$K3$ surface $S$, we show that 
the set of such moduli spaces 
contains  $2^{\rho(n-1)-1}$ distinct $\{Q^4,\bar{c}_2\}$-Hodge-isometry 
classes (and hence birational classes,  
Proposition \ref{thm-non-birational-moduli-spaces}).
For example, $S^{[7]}$ is not birational to the moduli space of 
rank $2$ sheaves with $c_1=0$ and $c_2=5$. 

\hide{
\subsection{Outline of the proofs}
\label{sec-outline-of-proofs}
{\bf Key ideas in the proofs of 
Theorems \ref{thm-introduction-invariant-Q-4} and 
\ref{thm-introduction-arithmetic-invariant-of-moduli}:}
We sketch the proofs of the case $i=1$ of both Theorems,
as the general case is similar. 
The homomorphism $\varphi^4$ is an isomorphism of free abelian groups, by
Proposition  
\ref{prop-varphi-k-is-independent-of-universal-sheaf}, which is
an easy consequence of 
our previous results,  
about the  
integral cohomology of the moduli space $\M$
(\cite{markman-original-archive-version} or 
Theorem \ref{thm-integral-generators} 
below). 
The equality $\varphi^4(v)=\frac{1}{2}\bar{c}_2(\M)$, 
in Theorem \ref{thm-introduction-arithmetic-invariant-of-moduli}, 
was proven in \cite{markman-monodromy-I} Lemma 4.9. 
The bilinear form on $Q^4(\M,\Integers)$
is defined as the pushforward of the Mukai pairing via $\varphi^4$.
Theorem \ref{thm-introduction-arithmetic-invariant-of-moduli} 
as well as parts \ref{thm-item-second-chern-class} and 
\ref{thm-intro-item-Q-4-is-the-mukai-lattice}
of Theorem \ref{thm-introduction-invariant-Q-4} follow.
The proof of Theorem \ref{thm-introduction-invariant-Q-4} part 
\ref{thm-item-unimodular-pairing} 
reduces to showing that the 
pushed forward bilinear form is $Mon(\M)$-invariant. 
The $Mon(\M)$-invariance implies, that the pushed forward 
bilinear form is well-defined for varieties $X$, deformation equivalent to 
$\M$. Set $e:=\varphi^4\circ (\varphi^2\circ\iota)^{-1}$.
The proof of Theorem \ref{thm-introduction-invariant-Q-4} part 
\ref{thm-item-pair-of-embeddings} reduces to showing that 
${\rm span}_\Integers\{e\}$ is $Mon(\M)$-invariant. 
The proofs of the $Mon(\M)$-invariance of the bilinear form,
and of ${\rm span}_\Integers\{e\}$, rely on further results in 
 \cite{markman-monodromy-I} and are postponed to the proof of
Theorem \ref{thm-introduction-invariant-Q-4}.
%
\EndProof

{\bf The strategy of the proof of 
Theorem \ref{thm-monodromy-operator-is-determined-by-its-weight-2-action}:}
Theorem \ref{thm-monodromy-operator-is-determined-by-its-weight-2-action}
is proven by induction. 
Let $i$ be an integer $\geq 2$. 
Consider $H^{2i}(X,\Integers)$ as the extension 
(\ref{eq-extension-of-Q-2i-by-A-2i-2}) of the {\em integral} 
$Mon(X)$-representation $Q^{2i}(X,\Integers)$ by $B^{2i}$.
Let $f$ be a monodromy operator, which acts as the identity on $B^{2i}$. 
Then $f$ belongs to the center of $Mon(X)$ and has order $2$,
by Proposition \ref{prop-mon-is-K2-times-W}, and is thus a semi-simple
element of $GL[H^*(X,\RationalNumbers)]$. 
In the induction step we need to show that $f$ acts as 
the identity also on $Q^{2i}(X,\Integers)$. 
Let $Mon'\subset Mon(X)$ and $\W'\subset \W$ be the derived subgroups.
Then $Mon'=\nu(\W')$, where 
$\nu:\W\rightarrow Mon(X)$ is given in (\ref{eq-nu}), by Proposition 
\ref{prop-mon-is-K2-times-W}. 
$Mon'$ is a normal subgroup of finite index in $Mon(X)$, by Proposition 
\ref{prop-mon-is-K2-times-W}, since $\W'$ has finite index in $\W$ 
(see Corollary 8.10 in
\cite{markman-monodromy-I} for the latter index). 

Set $n:=\dim_\RealNumbers(X)/4$.
We describe here the proof only in the range 
$4\leq 2i \leq n$, in which the $Mon'$-module 
$Q^{2i}(X,\Integers)$ is the extension of a trivial rank $1$ integral 
$Mon'$-module by a $Mon'$-sub-module $L$ of 
$Q^{2i}(X,\Integers)$ isomorphic to
$H^2(X,\Integers)$ (Theorem \ref{thm-introduction-invariant-Q-4}). 
$f$ must act via multiplication by $\pm 1$ on the 
$Mon'$-sub-representation $L\subset Q^{2i}(X,\Integers)$, 
since $L\otimes_\Integers\RationalNumbers$ is irreducible
of multiplicity $1$ in 
$Q^{2i}(X,\Integers)\otimes_\Integers\RationalNumbers$. 
Similarly, $f$ must act via multiplication by $\pm 1$ on the quotient
$Q^{2i}(X,\Integers)/L$. 
The order of the corresponding extension class 
$\epsilon_1\in \Ext^1_{Mon'}[Q^{2i}(X,\Integers)/L,L]$ is $2n-2$, 
which is larger than $2$ (Theorem 
\ref{thm-introduction-arithmetic-invariant-of-moduli} and Lemma 
\ref{lemma-orders-of-extension-classes-with-Mukai-lattice-in-the-middle}). 
Hence, $f$ must act via multiplication by $\pm 1$ on 
$Q^{2i}(X,\Integers)$. 

Let $\epsilon_2$ be the class in 
$\Ext^1_{Mon'}[Q^{2i}(X,\Integers),B^{2i}]$
of the $Mon'$-modules extension 
(\ref{eq-extension-of-Q-2i-by-A-2i-2}). 
If $f$ acts on $Q^{2i}(X,\Integers)$ as $-id$, then the
class $\epsilon_2$ has order at most $2$. 
We prove that 
the order of $\epsilon_2$ is $\geq 3$
(Proposition \ref{prop-extension-of-Q-bar-by-Z-hat-does-not-split}, which is
the main result of section \ref{sec-extension-classes}). 
Hence $f$ acts as the identity on $Q^{2i}(X,\Integers)$. 
Proposition \ref{prop-psi-is-G-equivariant} 
is key to the computation of the order of $\epsilon_2$;  
a $Mon'$-equivariant {\em rational} splitting of 
(\ref{eq-extension-of-Q-2i-by-A-2i-2})
is constructed for $X=S^{[n]}$, 
in terms of a normalization of the Chern character of the 
universal ideal sheaf
over $S\times S^{[n]}$.
\EndProof
}

\medskip
The paper is organized as follows.
In section \ref{sec-moduli-spaces-of-stable-sheaves} we 
review basic facts about the moduli spaces of 
stable sheaves on a $K3$ surface. 
Section \ref{sec-mukai-lattice-as-Q-i} is dedicated to the proof of Theorems 
\ref{thm-introduction-invariant-Q-4} and
\ref{thm-introduction-arithmetic-invariant-of-moduli}.
Section \ref{sec-degree-of-period-maps} 
starts with a careful comparison between
Hodge isometry classes of $H^2(X,\Integers)$ and of
$\{Q^{4},\bar{c}_2\}(X)$. A refined period map is defined and 
counter examples to a version of the Torelli question are exhibited.
Section \ref{sec-extension-classes} is dedicated to the
proof of Theorem 
\ref{thm-monodromy-operator-is-determined-by-its-weight-2-action}. 
Sections \ref{sec-degree-of-period-maps} and 
\ref{sec-extension-classes} are logically independent.

This is the final paper in a sequence of four, and it relies heavily
on the preceeding three paper 
\cite{markman-diagonal,markman-monodromy-I,markman-original-archive-version}.
In fact, the results of the current paper, with the exception of 
Theorem \ref{thm-monodromy-operator-is-determined-by-its-weight-2-action} and 
section \ref{sec-extension-classes}, 
were posted originally as the second part of the eprint 
arXiv:math.AG/0406016 v1
containing the 
motivation to and main application of the determination of the 
generators for the {\em integral} cohomology ring $H^*(S^{[n]},\Integers)$
in the first part.
The original version 
was split, following the recommendation of a referee, and the first part was
published as \cite{markman-original-archive-version}.
We are forced to include a somewhat lengthy summary of background material and
previous results in 
sections \ref{sec-moduli-spaces-of-stable-sheaves} and 
\ref{sec-resume-of-previous-results}.
All the geometric ingredients needed for the proofs 
are contained already in the three 
preceeding papers.  
The current paper assembles these geometric results, using
lattice theoretic techniques, to produce the concluding results of the 
project.

{\em Acknowledgments}: The work on this paper began while visiting 
the mathematics department of the University of Lille 1 during June 2003. 
I would like to thank Dimitri Markushevich and Armando Treibich for 
their hospitality. 

\subsection{Notation}
\label{sec-notation}

Given a smooth complex algebraic variety $X$, we denote by 
$K(X)$ the topological $K$-ring of vector bundles on $X$.
We denote by $x\cup y$ the product in $K(X)$. 
We let $K_{alg}X$ be its algebraic analogue. 
Given a morphism $f:X\rightarrow Y$,
we denote by $f^!:K(Y)\rightarrow K(X)$ the pullback. 
When $f$ is a proper morphism, we denote by
$f_!:K(X)\rightarrow K(Y)$ the topological Gysin map
(\cite{bfm}, \cite{karoubi} Proposition IV.5.24). 
We used above the assumption, that $X$ and $Y$ are smooth, 
which enables us to identify the topological $K$-cohomology groups with the 
$K$-homology groups $K_0^{top}$, for both spaces. 
Similarly, we identify the Grothendieck group $K_0^{alg}X$, of coherent 
sheaves, with $K_{alg}X$, replacing the class of a sheaf by that of
any of its locally free resolutions. 
The algebraic push-forward 
$f_!:K_{alg}X\rightarrow K_{alg}Y$ 
takes the class of a coherent sheaf $E$ on $X$, 
to the alternating sum $\sum_{i\geq 0} (-1)^i R^i_{f_*}E$ of 
the classes of the higher direct image sheaves on $Y$. 
When $Y$ is a point $\{pt\}$, 
we identify $K(\{pt\})$ with $\Integers$ and $f_!$ is the 
{\em Euler characteristic} 
$\chi:K(X)\rightarrow\Integers$.
Let
$
\alpha  :  K_{alg}X  \rightarrow  K(X)
$
be the natural homomorphism. 
The algebraic and topological Gysin maps are compatible via the
equality $f_!\circ\alpha=\alpha\circ f_!$ (see \cite{bfm}).


$S$ will denote a $K3$ surface and 
$\M$ the moduli space $\M_H(v)$ of sheaves on $S$ with class $v\in K(S)$,
which are stable with respect to a $v$-generic polarization $H$ 
(Definition \ref{def-v-suitable}).
The morphism $f_i$ is the projection from $S\times \M$ on the $i$-th factor,
$i=1,2$. 

$\widetilde{\Lambda}$ denotes the lattice 
$U^{\oplus 4}\oplus (-E_8)^{\oplus 2}$ in
Theorem \ref{thm-introduction-invariant-Q-4}, which is 
isometric to the Mukai lattice of a $K3$ surface $S$. 
The lattice $\Lambda$ will denote a fix copy of a lattice isomorphic
to $H^2(S^{[n]},\Integers)$, where the positive integer $n$ will be clear 
from the context.

%
\section{Moduli spaces of stable sheaves on a $K3$ surface}
\label{sec-moduli-spaces-of-stable-sheaves}

Let $S$ be a smooth and projective $K3$ surface and
$H$ an ample line bundle on $S$. We recall the definition of $H$-stability
for coherent sheaves, as it will be used in the proof of Proposition
\ref{thm-non-birational-moduli-spaces}.
The Hilbert polynomial of a coherent sheaf $F$ on $S$ is defined by
\[
P_F(n) \ \ := \ \ \chi(F\otimes H^n) := h^0(F\otimes H^n)-
h^1(F\otimes H^n)+h^2(F\otimes H^n).
\]
Let $r$ be the rank of $F$, $c_i:=c_i(F)$, $i=1,2$, its Chern classes, 
and $h:=c_1(H)$. 
Identify $H^4(S,\Integers)$ with $\Integers$ via the natural orientation of 
$S$.
The equality
\[
P_F(n) \ \ = \ \   (rh^2/2)n^2 + (h\cdot c_1)n + (c_1^2-2c_2)/2 + 2r
\]
follows from Hirzebruch-Riemann-Roch.
The degree $d$ of $P_F(n)$ is equal to 
the dimension of the support of $F$.
Let $l_0(F)/d!$ be 
the coefficient of $n^d$. 
Then $l_0(F)$ is a positive integer, provided $F\neq 0$. Explicitly, 
if $r>0$ then $l_0(F) = rh^2$. If $r=0$ and $d=1$ then
$l_0(F) =h\cdot c_1(F)$. If $r=0$ and $d=0$, then 
$l_0(F) =-c_2(F)$.
Given two polynomials $p$ and $q$ with real coefficients,
we say that $p \succ q$ (resp. $p \succeq q$) if
$p(n) > q(n)$ (resp. $p(n) \geq q(n)$) for all $n$ sufficiently large.

\begin{defi}
\label{def-stability}
A coherent sheaf $F$ on $S$ is called {\em $H$-semi-stable} (resp. 
{\em $H$-stable}) if it has support of pure dimension 
and any non-trivial subsheaf $F' \subset F$, $F' \neq (0)$, $F' \neq F$ 
satisfies 
\begin{equation}
\label{eq-inequality-of-normalized-Hilbert-polynomials}
\frac{P_{F'}}{l_0(F')} \preceq \frac{P_F}{l_0(F)} \ \ \ (\mbox{resp.} \prec).
\end{equation}
\end{defi}

Let 
$v\in K(S)$ be a class of rank 
$r\geq 0$ and first Chern class 
$c_1\in H^2(S,\Integers)$ of Hodge type $(1,1)$. 
The moduli space $\M_{H}(v)$, of isomorphism classes of 
$H$-stable sheaves of class $v$, 
is a quasi-projective scheme \cite{gieseker,simpson}. 

\begin{defi}
\label{def-v-suitable}
{\rm 
An ample line bundle $H$ is said to be $v$-{\em generic}, 
if every $H$-semi-stable sheaf with class $v$ is $H$-stable. 
}
\end{defi}

A $v$-generic polarization exists, if $v$ is assumed 
effective and primitive, and at least one of the following conditions
holds: $\rank(v)>0$, or $\chi(v)\neq 0$, or $c_1(v)$ generates the 
Neron-Severi group of $S$
(see \cite{yoshioka-chamber-str}, when $r>0$, and 
\cite{yoshioka-abelian-surface} Lemma 1.2, when $r=0$).
If $H$ is $v$-generic, then $\M_H(v)$ is projective \cite{gieseker,simpson}. 
{\em Caution:} The standard definition of the term $v$-{\em generic}
is more general, and does not assume that $v$ is primitive.

We assume, throughout the rest of the paper, that 
$v$ is an effective and primitive class in $K(S)$
of type $(1,1)$, and $H$ is $v$-generic.
Set $\M:=\M_H(v)$.

A {\em universal sheaf\/}  is a coherent sheaf $\E$ over $S\times \M$, 
flat over $\M$, 
whose restriction to $S\times \{m\}$,  $m\in\M$, is isomorphic 
to the sheaf $E_m$ on $S$ in the isomorphism class $m$.
The universal sheaf is canonical only up to tensorization by
the pull-back of a line-bundle on $\M$. 
A universal sheaf exists, if there exists a class $x$ in
$K_{alg}(S)$ satisfying $\chi(x\cup v)=1$ \cite{mukai-hodge}. 
Otherwise, there is a weaker notion of a 
{\em twisted} universal sheaf, denoted also by $\E$, 
where the twisting is encoded by a class $\alpha$ in \v{C}ech cohomology
$H^2(\M,\StructureSheaf{}^*)$, in the classical topology. 
For a triple $S$, $v$, $H$, as above, 
the class $\alpha$ is always topologically trivial;
it maps to $0$ in 
$H^2(\M,\F)$, where $\F$ is the sheaf of continuous complex valued 
invertible functions
\cite{markman-original-archive-version}. 
Consequently, $\E$ defines a class 
\begin{equation}
\label{eq-e-v}
[\E]  \ \ \in \ \ K(S\times\M),
\end{equation}
canonical up to tensorization by 
the pull-back of the class of a topological line-bundle on $\M$
\cite{markman-original-archive-version}. 

Let $v^\perp\subset K(S)$ be the sublattice orthogonal to $v$. 
Mukai introduced the natural homomorphism
\begin{eqnarray}
\label{eq-theta-v-from-v-perp}
\theta_v \ : \ v^\perp & \longrightarrow & H^2(\M,\Integers)
\\
\nonumber
\theta_v(x) & = & c_1\left[
-f_{2_!}\left(
f_1^!(x^\vee)\cup [\E]
\right)
\right].
\end{eqnarray}

\noindent
Note that $\theta_v$ is the restriction to $v^\perp$ of the homomorphism 
$\varphi^2$ given in
(\ref{eq-varphi-k}).

Beauville constructed an integral symmetric bilinear form on
the second cohomology of irreducible holomorphic symplectic manifolds
\cite{beauville-varieties-with-zero-c-1}. 
We will not need the intrinsic formula, but only the fact that it is 
monodromy invariant and its identification in 
the following theorem due to Mukai, Huybrechts, O'Grady, and Yoshioka:

\begin{thm} (\cite{ogrady-hodge-str}, 
\cite{yoshioka-abelian-surface} Theorem 8.1,
and \cite{yoshioka-note-on-fourier-mukai} Corollary 3.15)
\label{thm-irreducibility}
Let $v$ be a primitive and effective class and $H$ a $v$-generic ample
line-bundle. Set $\M:=\M_H(v)$.
\begin{enumerate}
\item 
$\M$ is a smooth, non-empty, irreducible symplectic, 
projective variety of dimension $\dim(v)=\langle v,v\rangle+2$. 
It is obtained by  deformations 
from the Hilbert scheme of $\frac{<v,v>}{2}+1$ points on $S$. 
\item
The homomorphism (\ref{eq-theta-v-from-v-perp}) 
is an isomorphism of weight 2 
Hodge structures with respect to Beauville's bilinear form on 
$H^2(\M,\Integers)$ when $\dim(v)\geq 4$.  When $\dim(v)=2$,
(\ref{eq-theta-v-from-v-perp}) factors through an isomorphism 
from $v^\perp/\Integers\cdot v$. 
\end{enumerate}
\end{thm} 

{\em Note:} Most of the cases of Theorem \ref{thm-irreducibility} are
proven in \cite{yoshioka-abelian-surface}
Theorem 8.1 under the additional assumption that $\rank(v)>0$
or $c_1(v)$ is ample. Corollary 3.15 in \cite{yoshioka-note-on-fourier-mukai}
proves the remaining cases ($\rank(v)=0$ and $c_1(v)$ not ample),
under the assumption that $\M$ is non-empty. The non-emptyness
is proven in \cite{yoshioka-note-on-fourier-mukai} Remark 3.4,
under the assumption that $c_1(v)$ is nef and in complete generality in an 
unpublished note communicated to the author by K. Yoshioka.

We recall next a result about generators for the cohomology ring
of $\M$. 
Given any cell complex $M$, the K\"{u}nneth Theorem provides an isomorphism
\[
K(S)\otimes K(M) \ \ \cong \ \ 
K(S\times M),
\]
given by the exterior product 
(\cite{atiyah-book}, Corollary 2.7.15). We used here the vanishing of 
$K^1(S)$. Use a basis  $\{x_1, \dots, x_{24}\}$ of $K(S)$
to write the K\"{u}nneth decomposition of the class (\ref{eq-e-v})
of the (possibly twisted) universal sheaf:
\begin{equation}
\label{eq-e-i}
[\E] \ \ \equiv \ \ 
\sum_{i=1}^{24} x_i\otimes u_i.
\end{equation}

Let $\{y_1,  \dots, y_{24}\}$ is a basis of $K(S)$ dual to
$\{x_1, \dots, x_{24}\}$ with respect to the Mukai pairing 
(\ref{eq-mukai-pairing-on-K-top}). Note that 
the K\"{u}nneth factor $u_i$ is simply $u(y_i)$, where
$u$ is the homomorphism given in (\ref{eq-u}). Indeed, 
the equality $f_{2_!}[f_1^!(x)\cup f_2^!(\xi)]=\chi(x)\cdot \xi$ holds, 
for any two classes $x\in K(S)$ and $\xi\in K(\M)$, 
by the definition of $\chi$ as a Gysin homomorphism, and
a standard property of Gysin maps (Property (M) in
section 4 of \cite{atiyah-hirzebruch-rr}). 
The K\"{u}nneth factors $u_i$ thus satisfy the equality 
$u_i=f_{2_!}[f_1^!(-y_i^\vee)\cup [\E]]=u(y_i)$.

\begin{thm}
\label{thm-integral-generators}
\cite{markman-original-archive-version}
\begin{enumerate}
\item
\label{thm-item-integral-generators}
The cohomology ring $H^*(\M,\Integers)$ is generated by the
Chern classes $c_j(u_i)$ of the classes $u_i\in K(\M)$, 
which are given in equation (\ref{eq-e-i}).
\item
\label{thm-item-odd-cohomology-vanishes}
The cohomology group $H^i(\M,\Integers)$ 
vanishes, for odd $i$, and is torsion free, when $i$ is even.
\end{enumerate}
\end{thm}

Generators for the cohomology ring
$H^*(X^{[n]},\RationalNumbers)$, with {\em rational} coefficients, were
found in \cite{lqw1} for any projective surface
$X$. A basis for the
{\em group} $H^*(X^{[n]},\Integers)/{\rm torsion}$
is given in
\cite{qin-wang} for many surfaces $X$, including $K3$'s.

We recall next the minimal degree of a non-trivial relation 
among these generators.
Let $\R(v)$ be the weighted polynomial ring generated by the vector spaces
$V_2:=v^\perp\otimes_\Integers \RationalNumbers$, and 
$V_{2i}:=K(S)\otimes_\Integers \RationalNumbers$, 
$2\leq i \leq n-1$, where $n:=\frac{1}{4}\dim_\RealNumbers\M$. 
Let $u:K(S) \rightarrow K(\M)$ be the homomorphism given in Equation
(\ref{eq-u}).
Define 
$h^{2i}:V_{2i}\rightarrow H^{2i}(\M,\RationalNumbers)$ by 
$
h^{2i}(x) =ch_i(u(x)),
$
and let 
\[
h:\R(v)\rightarrow H^*(\M,\RationalNumbers)
\] 
be the associated ring homomorphism. Denote by $I^d$ the degree $d$ summand of 
the relation ideal $\ker(h)$. 

\begin{new-lemma}
\label{lemma-10-in-markman-diagonal}
(Lemma 10 in \cite{markman-diagonal})
Assume $n\geq 2$. 
The homomorphism $h$ is surjective. It is injective in degree $\leq n$. 
If $n$ is odd, then $\dim(I^{n+1})=1$. 
If $n$ is even, then $\dim(I^{n+2})=24$. More generally, the dimension of 
$I^k$, $k$ even in the range $n< k \leq \frac{4n}{3}$, 
is the sum of the two Betti numbers
$
\dim(I^k) = 
b_{2(k-n-1)}(\M)+ b_{2(k-n-2)}(\M).
$
\end{new-lemma}

\begin{prop}
\label{prop-varphi-k-is-independent-of-universal-sheaf}
The homomorphism $\varphi^{2i}$, given in (\ref{eq-varphi-k}), 
is linear and surjective, for $i\geq 1$,   
and is independent of the choice of the 
class $[\E]$, for $i>1$. 
The homomorphism $\varphi^{2i}$ is an isomorphism, for
$4\leq 2i \leq \frac{1}{4}\dim_\RealNumbers\M$.
\end{prop}

\noindent
{\bf Proof:} 
The classes 
$c_i(f_{2_!}([\E]\cup f_1^!y))$ and 
$c_i([L]\cup f_{2_!}([\E]\cup f_1^!y))$ are equivalent  
modulo $A_2$, for any complex line bundle $L$ on $\M$.
The independence, of the choice of $[\E]$, follows 
in case $2i \geq 4$. 
The linearity of $\varphi^{2i}$ is verified via a direct calculation:
\[
c_i(x+y)-c_i(x)-c_i(y) \ = \ c_1(x)c_{i-1}(y)  + \ \cdots \ +  c_{i-1}c_1(y)
\ \in \ (A_{2i-2})^{2i},
\]
for $i\geq 2$, by the equality $c(x+y)=c(x)c(y)$ of Chern polynomials.
The surjectivity of $\varphi^{2i}$ follows from 
Theorem \ref{thm-integral-generators}.
$Q^{2i}(\M,\Integers)$ has rank $24$, for $2i$ in the specified range, 
by Lemma \ref{lemma-10-in-markman-diagonal}.
The injectivity in the specified range follows.
\EndProof

\hide{
The proof of Proposition \ref{prop-varphi-k-is-independent-of-universal-sheaf}
will depend on Lemmas \ref{lemma-10-in-markman-diagonal}
and \ref{lemma-chern-class-versus-chern-character}.

\begin{new-lemma}
\label{lemma-chern-class-versus-chern-character}
Let $X$ be a topological space and $y$ a class in $K(X)$. Then 
the degree $2i$ summands $ch_i(y)$ of the Chern character satisfy
$ch_2(y)  =  -c_2(y)+\frac{1}{2}c_1(y)^2$ and 

\begin{equation}
\label{eq-Girard-Formula}
ch_i(y) \ \ \ = \ \ \
\frac{(-1)^{i-1}}{(i-1)!} c_i(y) + \frac{(-1)^i}{(i-1)!} c_1(y)c_{i-1}(y) + 
\frac{1}{i!}g_i(c_1(y), \dots, c_{i-2}(y)),
\end{equation}
for $i\geq 3$, 
where $g_i$ is a universal polynomial with integral coefficients.
\end{new-lemma}

Lemma \ref{lemma-chern-class-versus-chern-character}
plays a minor role in the proof of Proposition
\ref{prop-varphi-k-is-independent-of-universal-sheaf}.
The main usage of the Lemma is in the proof of Theorem 
\ref{thm-monodromy-operator-is-determined-by-its-weight-2-action}
(see Lemma \ref{lemma-sigma-is-linear}). Particularly essential there is 
the fortunate equality, up to sign, of the two leading coefficients in
equation (\ref{eq-Girard-Formula}).

\noindent
{\bf Proof of Lemma \ref{lemma-chern-class-versus-chern-character}:} 
The product $(i!)ch_i(y)$ is, by definition, 
a polynomial in the Chern classes
$c_1(y), c_2(y), \dots, c_i(y)$. 
The coefficients are given by Girard's Formula 
(\cite{milnor} Ch. 16 Problem 16-A):

\noindent
${\displaystyle 
(-1)^i(i-1)!ch_i(y)  =
\sum_{d_1+2d_2+ \cdot + id_i=i}(-1)^{d_1+\cdots +d_i}
\frac{(d_1 + \cdots + d_i-1)!}{d_1!\cdots d_i!}
c_1(y)^{d_1}\cdots c_i(y)^{d_i}.
}$
%
\EndProof

}

\medskip
Let $F^{2i}K(S)$, $i\geq 0$, be the descending filtration, which is
the pullback of the weight filtration on $H^*(S,\Integers)$, 
via the Chern character. Then 
$F^0K(S)=K(S)$, $F^2K(S):=\{x \ : \ \rank(x)=0\}$, 
$F^4K(S):=\{x \ : \ \rank(x)=0 \ \mbox{and} \ c_1(x)=0\}$, and 
$F^{2i}K(S)=0$, for $i\geq 3$. 
Set $n:= \frac{1}{4}\dim_\RealNumbers\M$. 

\begin{thm}
\label{thm-vanishing-of-varphi-k}
{\rm
\cite{lehn-sorger-generators}
$H^*(\M,\RationalNumbers)$ is generated by 
$
c_i(f_{2_!}(f_1^!(x)\cup [\E])),
$
for $x\in F^{2j}K(S)$ and 
$0\leq i\leq \lfloor\frac{n}{2}\rfloor+j$.
}
\end{thm}

The Theorem implies, that the homomorphism 
(\ref{eq-varphi-k}) induces a surjective homomorphism 
\[
\varphi^{2\lfloor\frac{n}{2}\rfloor+2j}:
F^{2j}K(S)\otimes_\Integers\RationalNumbers\rightarrow
Q^{2\lfloor\frac{n}{2}\rfloor+2j}(\M,\RationalNumbers),
\] 
for any non-negative integer $j$.
Consequently, $\rank(Q^{2\lfloor\frac{n}{2}\rfloor+2})\leq 23$, 
$\rank(Q^{2\lfloor\frac{n}{2}\rfloor+4})\leq 1$, and
$\rank(Q^{2\lfloor\frac{n}{2}\rfloor+2i})=0$, for $i\geq 3$.

The statement of the Theorem is not sharp, at least for small $n$, as
$H^*(S^{[2]},\RationalNumbers)$ is generated by $H^2(S^{[2]})$
and $H^*(S^{[3]},\RationalNumbers)$ is generated by $H^2$ and $H^4$,
by Lemma \ref{lemma-10-in-markman-diagonal}. 

\hide{
The Conjecture would follow from its Hilbert scheme case
$\M=S^{[n]}$, by Theorem \ref{thm-irreducibility}.
Lehn and Sorger determined the ring structure of
$H^*(S^{[n]},\RationalNumbers)$ in \cite{lehn-sorger}.
Using these results, they easily reduced Conjecture
\ref{conj-vanishing-of-varphi-k} to the analogous statement about
$H^*((\ComplexNumbers^2)^{[n]},\RationalNumbers)$
\cite{lehn-sorger-generators}.
The latter statement is combinatorial in nature, as 
$H^*((\ComplexNumbers^2)^{[n]},\RationalNumbers)$ is isomorphic to the
center $Z(\RationalNumbers[\Sigma_n])$ of the group ring of the
symmetric group $\Sigma_n$ on $n$ letters.
The grading of $\RationalNumbers[\Sigma_n]$ is such, that the degree
of a permutation $\pi$ is $2s$, where
$s$ is the minimal length of a decomposition 
$\pi=\tau_1\cdot \cdots \cdot \tau_s$ as a product of transpositions. 
A detailed outline of the proof of the combinatorial statement is given in 
\cite{lehn-sorger-generators}.
}

%
\section{The Mukai lattice as a quotient Hodge structure}
\label{sec-mukai-lattice-as-Q-i}

Theorems \ref{thm-introduction-invariant-Q-4}
and \ref{thm-introduction-arithmetic-invariant-of-moduli}
are proven in section 
\ref{sec-proof-of-thm-arithmetic-invariant}, after 
we calculate an extension class in 
section \ref{sec-elementary-calculations-of-extension-classes},
and review known related results in section 
\ref{sec-resume-of-previous-results}.

%
\subsection{Elementary calculations of extension classes}
\label{sec-elementary-calculations-of-extension-classes}

Let $v$ be a class in $K(S)$ satisfying $(v,v)\geq 2$. 
Let $G$ be any subgroup of finite index in the group $OK(S)_v$ of isometries
stabilizing $v$.
Modules over the group ring $\Integers[G]$ will be referred to as $G$-modules.
Given $G$-modules $V$, $W$, 
we denote by $\Ext^i_G(V,W)$ the $i$-th extension group 
$\Ext^i_{\Integers[G]}(V,W)$ of the $\Integers[G]$-modules
(\cite{hilton-stammbach} Ch. IV section 7). 
Consider the short exact sequence
\begin{equation}
\label{eq-non-split-cyclic-extensions-withtrivial-G-action}
0\rightarrow d\Integers/de\Integers 
\rightarrow \Integers/de\Integers\rightarrow 
\Integers/d\Integers\rightarrow 0,
\end{equation}
where $d$, $e$ are integers, $d\neq 0$, 
and $G$ acts trivially on the three groups.

\begin{new-lemma}
\label{lemma-non-split-cyclic-extensions-withtrivial-G-action}
The order of the extension class of 
(\ref{eq-non-split-cyclic-extensions-withtrivial-G-action})
is $\gcd(d,e)$.
\end{new-lemma}

\noindent
{\bf Proof:} 
The order is equal to that of the extension class of $\Integers$-modules,
as the $G$-action is trivial. The latter order is well known
(\cite{hilton-stammbach}, Ch. III section 4).
\EndProof

\begin{new-lemma}
\label{lemma-characterization-of-pullback-pushforward}
Let $\begin{array}{ccccccc}
0\rightarrow & A_1 & \rightarrow & B_1 & \rightarrow & C_1 & \rightarrow 0
\\
 & \alpha \ \downarrow \ \hspace{1ex} & &
 \beta \ \downarrow \ \hspace{1ex} & &
 \gamma \ \downarrow \ \hspace{1ex} 
\\
0\rightarrow & A_2 & \rightarrow & B_2 & \rightarrow & C_2 & \rightarrow 0
\end{array}
$
be a commutative diagram of $G$-modules with short exact rows. 
Denote by $\epsilon_i$ the extension class of the $i$-th row.
\begin{enumerate}
\item
Assume that $A_1=A_2=A$ and $\alpha$ is the identity. Then 
$\epsilon_1=\gamma^*(\epsilon_2)$, where
$\gamma^*:\nolinebreak\Ext^1_G(C_2,A)\rightarrow \Ext^1_G(C_1,A)$
is the pullback homomorphism.
\item
Assume that $C_1=C_2=C$ and $\gamma$ is the identity. 
Then $\epsilon_2=\alpha_*(\epsilon_1)$, where 
$\alpha_*:\nolinebreak\Ext^1_G(C,A_1)\rightarrow \Ext^1_G(C,A_2)$ is 
the pushforward homomorphism.
\end{enumerate}
\end{new-lemma}

\noindent
{\bf Proof:}
Follows immediately from the definitions of the pushforward and pullback of 
extensions.
\EndProof

Identify the quotient $K(S)/\Integers v$ with 
$(v^\perp)^*$, using the fact that $K(S)$ is unimodular. 
Consider the two short exact sequences:
\begin{eqnarray}
\label{eq-extension-of-v-perp-dual-by-span-v}
& 0 \rightarrow \Integers v \rightarrow K(S) \rightarrow 
(v^\perp)^* \rightarrow 0,
\\
\label{eq-extension-ofspan-v-dual-by-v-perp}
& 0 \rightarrow v^\perp \rightarrow K(S) \rightarrow 
K(S)/v^\perp \rightarrow 0.
\end{eqnarray}

\begin{new-lemma}
\label{lemma-orders-of-extension-classes-with-Mukai-lattice-in-the-middle}
The order of the extension classes of the short exact 
sequences of $G$-modules 
(\ref{eq-extension-of-v-perp-dual-by-span-v}) and 
(\ref{eq-extension-ofspan-v-dual-by-v-perp}) is $(v,v)$.
\end{new-lemma}

\noindent
{\bf Proof:}
The order of the extension class of 
(\ref{eq-extension-ofspan-v-dual-by-v-perp}) divides $k$, if and only if
there exists $\varphi\in \Hom_G(K(S),v^\perp)$, which restricts to 
$v^\perp$ as multiplication by $k$. 
The order of the extension class of 
(\ref{eq-extension-of-v-perp-dual-by-span-v}) divides $k$, if and only if
there exists $\psi\in \Hom_G((v^\perp)^*,K(S))$, whose composition with
$K(S)\rightarrow (v^\perp)^*$ acts on $(v^\perp)^*$ as 
multiplication by $k$. Setting $\psi=\varphi^*$, we see that the two 
conditions are equivalent. Hence, the order of the two classes are equal.

The extension class of 
$
0\rightarrow \Integers v\rightarrow \Integers\frac{v}{(v,v)}\rightarrow
\Integers\frac{v}{(v,v)}/\Integers v\rightarrow 0
$
has order $(v,v)$, by the case $e=0$ and $d=(v,v)$ of Lemma
\ref{lemma-non-split-cyclic-extensions-withtrivial-G-action}. 
The extension (\ref{eq-extension-of-v-perp-dual-by-span-v})
is the pullback of the above extension
via the natural homomorphism 

\hspace{8ex}
$(v^\perp)^*\rightarrow (v^\perp)^*/v^\perp\cong 
K(S)/[\Integers v+v^\perp] \cong
\Integers\frac{v}{(v,v)}/\Integers v,
$\\
by Lemma \ref{lemma-characterization-of-pullback-pushforward}.
The pullback homomorphism 
$\Ext^1_G((v^\perp)^*/v^\perp,\Integers v)\rightarrow 
\Ext^1_G((v^\perp)^*,\Integers v)$ is injective, since its kernel is the 
image of $\Hom_G(v^\perp,\Integers v)$, which vanishes,
since $G$ has finite index in $OK(S)_v$. 
Hence, the order of the extension class of 
(\ref{eq-extension-of-v-perp-dual-by-span-v}) is $(v,v)$ as well.
\EndProof

%
\subsection{Summary of previous results on the monodromy group}
\label{sec-resume-of-previous-results}

%
\subsubsection*{Orientation characters}
Let $X$ be an irreducible holomorphic symplectic manifold, 
which is deformation equivalent to
the Hilbert scheme $S^{[n]}$, $n\geq 2$, of length $n$ subschemes
of a K3 surface $S$. 
The lattice $H^2(X,\RealNumbers)$ has signature $(3,20)$
(\cite{beauville-varieties-with-zero-c-1} or Theorem 
\ref{thm-irreducibility}).
A $3$-dimensional subspace of $H^2(X,\RealNumbers)$ is said to be 
{\em positive-definite}, if the 
Beauville-Bogomolov pairing restricts to it as a positive definite pairing. 
The unit $2$-sphere, in any positive-definite $3$-dimensional subspace, 
is a deformation retract of
the positive cone $\C_+\subset H^2(X,\RealNumbers)$, 
given in (\ref{eq-positive-cone}).
Hence, $H^2(\C_+,\Integers)$ is isomorphic to $\Integers$
and is a natural representation of the isometry group
$OH^2(X,\RealNumbers)$. We get an associated character, 
denoted by
\begin{equation}
\label{eq-orientation-character-of-H2}
\eta \ : \ OH^2(X,\RealNumbers) \ \ \ \longrightarrow \ \ \
\{\pm 1\}, 
\end{equation}
which we call the {\em orientation character}. 
\begin{defi}
\label{def-orientation-preserving}
{\rm
We denote the kernel of $\eta$ 
by $O^+H^2(X,\RealNumbers)$ and refer to its elements as
{\em orientation preserving isometries}. 
}
\end{defi}
Note, that $\eta(-id)=-1$ and the group
$\W$, given in (\ref{eq-W}), is a subgroup of $O^+H^2(X,\Integers)$.
Choose a K\"{a}hler class $\kappa$ and a holomorphic symplectic 
$2$-form $\sigma$.
Then $\{{\rm Re}(\sigma),{\rm Im}(\sigma),\kappa\}$ is a basis of 
a $3$-dimensional positive definite subspace of $H^2(X,\RealNumbers)$, so
it determines an orientation of $\C_+$. 
\begin{defi}
\label{def-distinguished-orientation}
{\rm
The above orientation of $\C_+$ is independent of the
choices of $\kappa$ and $\sigma$ \cite{huybrechts-mirror-symmetry}.
We will call it the {\em distinguished orientation} of $\C_+$.
}
\end{defi}

$Mon^2$ 
is also a subgroup of $O^+H^2(X,\Integers)$.
$O^+H^2(X,\Integers)$ is the symmetry group  
preserving the integral structure, 
the Beauville-Bogomolov bilinear form, and the orientation.

%
%

Let $OK(S)$ be the isometry group of the Mukai lattice and
$OK(S)_v$ the subgroup stabilizing $v$.
The positive cone $\widetilde{\C}_+$
of $K(S)\otimes_\Integers\RealNumbers$ is 
homotopic to the unit $3$-sphere, in any positive definite four-dimensional
subspace of $K(S)\otimes_\Integers\RealNumbers$. 
The orientation character 
\begin{equation}
\label{eq-tilde-eta}
\tilde{\eta} \ : \ OK(S) \ \ \ \longrightarrow  \ \ \ 
\{\pm1\}
\end{equation}
is determined by the action of $OK(S)$ on 
$H^3(\widetilde{\C}_+,\Integers)$. The two orientation characters,
$\tilde{\eta}$ of $OK(S)$ and 
$\eta$ of $O(v^\perp)$ given in (\ref{eq-orientation-character-of-H2}), 
are compatible, $\tilde{\eta}(g)=\eta(\restricted{g}{v^\perp})$,
under the restriction from the stabilizer 
$OK(S)_v$ to $O(v^\perp)$.
However, $\tilde{\eta}(-1)=-\eta(\restricted{-1}{v^\perp})$.

\hide{
We choose an orientation of the positive cone $\widetilde{\C}_+$ 
as follows. Let $U$ be the sublattice of $K(S)$ 
of classes $x$ with trivial $c_1(x)$. $U$ is isometric to 
the hyperbolic plane and its orthogonal complement $U^\perp$ 
is isometric to $H^2(S,\Integers)$. Thus
$K(S)$ decomposes as 
the orthogonal direct sum $U\oplus H^2(S,\Integers)$.  
The orientation of $\widetilde{\C}_+$ will be
determined by the choice of orientations of each of the two summands.
The positive cone of $U$ has two connected components, and an
orientation amounts to a choice of a component. 

\begin{defi}
\label{def-orientation-of-the-Mukai-lattice}
The orientation of the positive cone $\widetilde{\C}_+$ 
of the Mukai lattice $K(S)$ is the one determined by the
distinguished orientation of $H^2(S,\Integers)$
and the choice of the component of the 
positive cone of $U$
containing the class of ideal sheaves of subschemes of length $2$. 
\end{defi}

%
\subsubsection*{The monodromy representation of the Mukai groupoid}

We define first a groupoid $\G$, i.e., a category, all of whose morphisms are
isomorphisms. An object of $\G$ is the data needed
to construct a smooth and compact moduli space, i.e., a triple
$(S,v,H)$ consisting of a projective $K3$ surface $S$,
an effective and primitive class $v\in K^{1,1}S$, 
and a $v$-generic polarization $H$ on $S$.
A morphism $g\in \Hom_\G((S_1,v_1,H_1),(S_2,v_2,H_2))$
is an isometry $g:K(S_1)\rightarrow K(S_2)$, satisfying
$g(v_1)=v_2$. We do not assume $g$ to preserve the Hodge structure.
Nor do we assume any compatibility between $g$ and the polarizations $H_i$.

We define next the orientation character 
\[
cov \ : \ \G \ \ \ \rightarrow \ \ \ \{\pm 1\},
\]
a functor to the group of order $2$. 
All objects are sent to the group $\{\pm 1\}$.
A morphism $g$ is sent to $1$, if the corresponding isometry is
orientation preserving (Definition 
\ref{def-orientation-of-the-Mukai-lattice}), and to $-1$ otherwise.

Let $(S_i,v_i,H_i)$, $i=1,2$, be two objects of $\G$, and
$g\in \Hom_{\G}((S_1,v_1,H_1),(S_2,v_2,H_2))$ a morphism. 
Denote by 
\[
(g\otimes 1) \ : \ K(S_1)\otimes_\Integers K(\M_{H_1}(v_1))
\ \ \ \longrightarrow \ \ \ 
K(S_2)\otimes_\Integers K(\M_{H_1}(v_1))
\]
the homomorphism acting via the identity on the second factor.
We identify its domain with $K[S_1\times \M_{H_1}(v_1)]$ 
via the K\"{u}nneth Theorem. 
Choose universal classes $[\E_{v_i}]$ in 
$K(S_i\times \M_{H_i}(v_i))$, $i=1,2$, as in (\ref{eq-e-v}).
Let $D:K(S_2)\rightarrow K(S_2)$ be the involution, which
takes a class $x$ to its dual $x^\vee$. 
Set $m:=2+(v_1,v_1)$. 
Let $\pi_{ij}$ be the projection from 
$\M_{H_1}(v_1)\times S_2\times \M_{H_2}(v_2)$ to the product of the
$i$-th and $j$-th factors.
We define a class in the middle cohomology
$H^{2m}(\M_{H_1}(v_1)\times \M_{H_2}(v_2),\Integers)$ of the product:
\[
\monclass(g) \ \ := \ \ 
\left\{
\begin{array}{ccc}
c_m\left(
-\pi_{13_!}\left\{
\pi_{12}^!\left[
(g\otimes 1)[\E_{v_1}]
\right]^\vee\cup \pi_{23}^![\E_{v_2}]
\right\}
\right)
 & \mbox{if} & cov(g)=1,
\\
c_m\left(
-\pi_{13_!}\left\{
\pi_{12}^!\left[
(Dg\otimes 1)[\E_{v_1}]
\right]\cup \pi_{23}^![\E_{v_2}]
\right\}
\right)
 & \mbox{if} & cov(g)=-1.
\end{array}
\right.
\]
Denote by 
\begin{equation}
\label{eq-monrep-g}
\monrep(g) \ : \ H^*(\M_{H_1}(v_1),\Integers) \ \ \ 
\longrightarrow \ \ \ H^*(\M_{H_2}(v_2),\Integers)
\end{equation}
the homomorphism obtained from $\monclass(g)$
using the K\"{u}nneth and Poincare-Duality Theorems. 
Let ${\A}lg$ denote the category of associative algebras with a unit.
Denote by $ch:K(S_i)\rightarrow H^*(S_i,\Integers)$ the
Chern character isomorphism, and  the conjugate 
$ch\circ g \circ (ch^{-1})$ by 
$\bar{g}:H^*(S_1,\Integers)\rightarrow H^*(S_2,\Integers)$.

\begin{thm}
\label{thm-symmetries-of-moduli-spaces}
(\cite{markman-monodromy-I} Theorems 1.2 and 1.6)
Let $(S_i,v_i,H_i)$, $i=1,2$, and $g$ be as above. 
\begin{enumerate}
\item
The homomorphism $\monrep(g)$ is an algebra isomorphism and
a parallel-transport operator (Definition \ref{def-Mon-X-Y}). 
\item
The assignment 
\[
\monrep \ : \ \G \ \ \ \longrightarrow \ \ \ {\A}lg,
\]
sending an object $(S,v,H)$ of
$\G$ to the cohomology algebra $H^*(\M_{H}(v),\Integers)$, 
and a morphism $g$ as above to the isomorphism $\monrep(g)$,
is a functor. 
\item
\label{thm-item-mon-g-sends-a-universal-classes-to-such}
There exists a topological complex line bundle $\ell$ on
$\M_{H_2}(v_2)$ satisfying one of the following equations: 
\[
\begin{array}{ccccc}
(\bar{g}\otimes \monrep(g))(ch[\E_{v_1}]) & = & 
ch([\E_{v_2}]\cup f_2^*\ell),
 & \mbox{if} & cov(g)=1,
\\
(\bar{D}\bar{g}\otimes \monrep(g))(ch[\E_{v_1}]) & = & 
ch([\E_{v_2}]^\vee\cup f_2^*\ell),
 & \mbox{if} & cov(g)=-1.
\end{array}
\]
\end{enumerate}
\end{thm}
}
%
\subsubsection*{The monodromy representation of the stabilizer}
Let $S$ be a projective $K3$ surface, $v\in K^{1,1}(S)$
an effective and primitive class, 
and $H$ a $v$-generic polarization on $S$ 
(Definitions \ref{def-effective-class} and \ref{def-v-suitable}). 
Set $\M:=\M_H(v)$. 
Given an isometry $g\in OK(S)$, 
denote by $g\otimes 1$ the endomorphism of $K(S)\otimes_\Integers K(\M)$
acting via the identity on the second factor.
We identify $K(S)\otimes_\Integers K(\M)$ with $K(S\times \M)$ 
via the K\"{u}nneth Theorem. 
Choose a universal class $[\E]$ in 
$K(S\times \M)$ as in (\ref{eq-e-v}).
Let $D\in OK(S)$ be the involution, which
takes a class $x$ to its dual $x^\vee$. 
Set $m:=2+(v,v)$. 
Let $\pi_{ij}$ be the projection from 
$\M\times S\times \M$ to the product of the
$i$-th and $j$-th factors.
We define a class in the middle cohomology
$H^{2m}(\M\times \M,\Integers)$ of the product:
\[
\monclass(g) \ \ := \ \ 
\left\{
\begin{array}{ccc}
c_m\left(
-\pi_{13_!}\left\{
\pi_{12}^!\left[
(g\otimes 1)[\E]
\right]^\vee\cup \pi_{23}^![\E]
\right\}
\right)
 & \mbox{if} & \tilde{\eta}(g)=1,
\\
c_m\left(
-\pi_{13_!}\left\{
\pi_{12}^!\left[
(Dg\otimes 1)[\E]
\right]\cup \pi_{23}^![\E]
\right\}
\right)
 & \mbox{if} & \tilde{\eta}(g)=-1.
\end{array}
\right.
\]
Denote by 
\begin{equation}
\label{eq-monrep-g}
\monrep(g) \ : \ H^*(\M,\Integers) \ \ \ 
\longrightarrow \ \ \ H^*(\M,\Integers)
\end{equation}
the homomorphism obtained from $\monclass(g)$
using the K\"{u}nneth and Poincare-Duality Theorems. 

The group $Q^{2i}(\M,\Integers)$ 
is a representation of the monodromy group $Mon(\M)$ 
(Definition \ref{def-monodromy}). 
Let $Mon^{2i}$ be the 
image of $Mon(\M)$ in $GL(Q^{2i}(\M,\Integers))$. 
Denote by $ch:K(S)\rightarrow H^*(S,\Integers)$ the
Chern character isomorphism, and  the conjugate 
$ch\circ g \circ (ch^{-1})$, $g\in OK(S)$, by 
$\bar{g}\in GL[H^*(S,\Integers)]$.

\begin{thm} 
\label{thm-summary-of-monodromy-results}
\begin{enumerate}
\item
\label{thm-item-markman-monodromy-I-thm-1.5}
(\cite{markman-monodromy-I}, Theorem 1.6)
Let $v$ be an effective and primitive class in $K^{1,1}(S)$,
satisfying $(v,v)\geq 2$, and $H$ a $v$-generic polarization.
The homomorphism $\monrep(g)$ is a monodromy operator, 
for every $g\in OK(S)_v$, and the map 
\begin{equation}
\label{eq-mon-representation-of-stabilizer}
\monrep \ : \ OK(S)_v \ \ \ \longrightarrow \ \ \ 
Mon(\M)
\end{equation}
is a group homomorphism. The homomorphism $\monrep$
is injective, if $(v,v)\geq 4$, and its kernel
is generated by the involution $\rho_v$, given in (\ref{eq-W}), if $(v,v)=2$.
The image $\monrep[OK(S)_v]$ is a normal subgroup of finite
index in the monodromy group $Mon(\M)$. 
\item
\label{thm-item-markman-monodromy-I-lemma-4.4}
(\cite{markman-monodromy-I}, Lemma 4.7)
The kernel $K^2$ of $Mon\rightarrow Mon^2$  
is a subgroup of the center of $Mon(\M)$.
$K^2$ is finite of exponent $2$. 
\item
\label{thm-item-mon-g-sends-a-universal-classes-to-such}
(\cite{markman-monodromy-I}, Theorem 1.2, part 4 and equation (1.22))
There exists a topological complex line bundle $\ell_g$ on
$\M$, for each $g\in OK(S)_v$, 
satisfying one of the following equations: 
\[
\begin{array}{ccccc}
(\bar{g}\otimes \monrep(g))(ch[\E]) & = & 
ch([\E]\cup f_2^*\ell_g),
 & \mbox{if} & \tilde{\eta}(g)=1,
\\
(\bar{D}\bar{g}\otimes \monrep(g))(ch[\E]) & = & 
ch([\E]^\vee\cup f_2^*\ell_g),
 & \mbox{if} & \tilde{\eta}(g)=-1.
\end{array}
\]
\item
\label{thm-item-markman-monodromy-I-lemma-4.5}
(\cite{markman-monodromy-I}, Lemma 4.8)
Let $i$ be an integer in the range 
$2\leq i \leq \dim_\RealNumbers(\M)/8$. 
Denote by $\monrep^{2i}:OK(S)_v\rightarrow Mon^{2i}$ 
the composition of $\monrep$ with the projection $Mon\rightarrow Mon^{2i}$.
Then the following equations hold.
\begin{eqnarray}
\label{eq-pullback-of-Q-2i-to-stabilizer}
\monrep^{2i}(g) & = & 
\tilde{\eta}(g)^i\left[
\varphi^{2i}\circ g\circ (\varphi^{2i})^{-1}
\right], \ \ \ \mbox{for} \ i>1,
\\
\label{eq-pullback-of-Q-2-to-stabilizer}
\monrep^2(g) & = & 
\tilde{\eta}(g)\left[
\theta_v\circ g\circ (\theta_v)^{-1}
\right],
\end{eqnarray}
where $\varphi^{2i}$ is given in (\ref{eq-varphi-k}) and
$\theta_v$ in (\ref{eq-theta-v-from-v-perp}).
\item
\label{thm-item-N-2-is-W}
(\cite{markman-monodromy-I}, Lemma 4.10)
The image $\monrep^2[OK(S)_v]$ in $Mon^2$ is the subgroup 
$\W$ of $O^+[H^2(\M,\Integers)]$, given in (\ref{eq-W}).
\item
\label{thm-item-identification-of-residue-of-chern-class}
(\cite{markman-monodromy-I}, Lemma 4.9) The equality
$\bar{c}_{2i}(\M)=\varphi^{4i}(2v)$ holds in
$Q^{4i}(\M,\Integers)$, for $i$ an integer, $i\geq 1$.
\end{enumerate}
\end{thm}

\begin{rem}
\label{rem-varphi-in-terms-of-chern-character-vs-chern-class}
{\rm
The version of equality 
(\ref{eq-pullback-of-Q-2i-to-stabilizer}), proven in 
Lemma 4.8 of \cite{markman-monodromy-I}, used the isomorphism 
$\psi^{2i}$, which we now define.
Replace the Chern class $c_i$ in (\ref{eq-tilde-varphi-2i}),
by the $i$-th coefficient $ch_i$ of the Chern character of the same class,
to obtain the homomorphism 
$\tilde{\psi}^{2i}: K(S) \rightarrow H^{2i}(\M,\RationalNumbers)$.
Define 
$\psi^{2i}: K(S) \rightarrow Q^{2i}(\M,\RationalNumbers)$
as the composition of $\tilde{\psi}^{2i}$ with the projection. 
The equation 
$
\varphi^{2i} \ \ = \ \ (-1)^{(i-1)}(i-1)!\psi^{2i}
$
holds, since $ch_i$ is the sum of 
$(-1)^{(i-1)}\frac{c_i}{(i-1)!}$ and a polynomial in the Chern classes 
$c_j$, for $j<i$. 
Equation (\ref{eq-pullback-of-Q-2i-to-stabilizer}) thus follows.
}
\end{rem}

Let 
\begin{equation}
\label{eq-N}
N
\end{equation} 
be the image of $OK(S)_v$ in $Mon(\M)$ 
via the homomorphism $\monrep$ in 
(\ref{eq-mon-representation-of-stabilizer}),
and $N^j$ the image of $N$ in $Mon^j$. 
The subgroups $N$ of $Mon(\M)$ and
$N^j$ of $Mon^j$, being normal subgroups,
are well defined for any 
irreducible holomorphic symplectic manifold $X$ deformation equivalent to
$\M$. 
Assume $(v,v)\geq 6$.
Then $N\rightarrow N^4$ is injective, since equation 
(\ref{eq-pullback-of-Q-2i-to-stabilizer}) yields
$\monrep^4(g)=\varphi^4\circ g \circ (\varphi^4)^{-1}$.
$N^2$ is the subgroup $\W H^2(\M,\Integers)$ of
$O^+H^2(\M,\Integers)$, by Theorem
\ref{thm-summary-of-monodromy-results} part \ref{thm-item-N-2-is-W}. 
The following diagram summarizes the relationship between $N^2$
and $N^4$ spelled-out in Theorem \ref{thm-summary-of-monodromy-results}.
\begin{equation}
\label{eq-commutative-diagram-of-monodromy-groups}
\begin{array}{ccccc}
N^2 & \LongIsomLeftArrow & N & \LongIsomRightArrow & N^4
\\
Ad_{\theta_v} \ \downarrow \ \cong \hspace{1ex} & & & &
\hspace{1ex} \cong  \ \downarrow \ Ad_{\varphi^4} 
\\
\W(v^\perp)  & & \LongIsomLefttArrowOf{\mu} & & OK(S)_v.
\end{array}
\end{equation}
The homomorphism $\mu$, making the above diagram commutative,
is determined explicitly as follows. 
An isometry $g\in OK(S)_v$
is sent by $\mu$ to the restriction 
of $\tilde{\eta}(g)\cdot g$ to $v^\perp$.
$N\rightarrow N^2$ is an isomorphism, by the following Lemma.

\begin{new-lemma}
\label{lemma-mu-is-an-isomorphism}a)
The homomorphism $\mu$ is surjective, if $(v,v)\geq 2$. 
$\mu$ is an isomorphism if $(v,v)\geq 4$, and its kernel
is generated by the involution $\rho_v$, given in (\ref{eq-W}), 
if $(v,v)=\nolinebreak2$.
b) There exists a unique $Mon(\M)$-equivariant isomorphism
$\nu:\W{H}^2(\M,\Integers)\rightarrow N$, satisfying 
$\nu\circ\monrep^2=\monrep$.
\end{new-lemma}

\noindent
{\bf Proof:}
a) The surjectivity is precisely Theorem 
\ref{thm-summary-of-monodromy-results} 
part \ref{thm-item-N-2-is-W}.
Let $g$ be a non-trivial element 
in the kernel of $\mu$. Then $g$ satisfies $g(v)=v$ and it 
restricts as $-id$ to $v^\perp$. 
Such $g$ can exist only if the extension 
(\ref{eq-extension-ofspan-v-dual-by-v-perp}),
of $OK(S)_v$-modules, has class
$\epsilon$ of order at most $2$ in 
$\Ext^1_{OK(S)_v}(K(S)/v^\perp,v^\perp)$. 
The order of $\epsilon$ is $(v,v)$, 
by Lemma 
\ref{lemma-orders-of-extension-classes-with-Mukai-lattice-in-the-middle}. 

b) The homomorphism $\monrep^2:OK(S)_v\rightarrow \W{H}^2(\M,\Integers)=N^2$ 
is equal to $Ad_{\theta_v}^{-1}\circ\mu$
and is thus an isomorphism, if $(v,v)\geq 4$. 
Hence the projection $N\rightarrow \W{H}^2(\M,\Integers)$ is an isomorphism,
and we let $\nu$ be its inverse. If $(v,v)=2$,
the projection is still an isomorphism, since $H^*(\M,\RationalNumbers)$
is generated by $H^2(\M,\RationalNumbers)$.
\EndProof

Assume $(v,v)\geq 10$. Then $N\rightarrow N^6$ is an isomorphism, 
by equation (\ref{eq-pullback-of-Q-2i-to-stabilizer}).
We get the commutative diagram:
\begin{equation}
\label{eq-commutative-diagram-of-monodromy-groups-hlf-integer-case}
\begin{array}{ccccc}
N^2 & \LongIsomLeftArrow & N & \LongIsomRightArrow & N^6
\\
Ad_{\theta_v} \ \downarrow \ \cong \hspace{1ex} & & & &
\hspace{1ex} \cong  \ \downarrow \ Ad_{\varphi^6} 
\\
\W(v^\perp)  & & \LongIsomRightArrowOf{ext} & & ext[\W(v^\perp)],
\end{array}
\end{equation}
where the extension homomorphism 
$ext:\W(v^\perp)\rightarrow OK(S)$  
is determined by the equation $ext(\mu(h))=h\cdot \tilde{\eta}(h)$, 
for all $h\in OK(S)_v$. 
The homomorphism $ext$ sends $g\in \W(v^\perp)$ to the
unique isometry of $K(S)$, whose restriction to $v^\perp$ is $g$. 
The extension $ext(g)$ sends $v$ to $v$ or $-v$, 
depending on $\mu^{-1}(g)$ being orientation preserving or reversing.

Verbitsky constructed a representation of 
$Spin[H^2(X,\RealNumbers)]$ on the cohomology 
of an irreducible holomorphic symplectic manifold 
\cite{verbitsky,looijenga-lunts}. 
When the odd cohomology of $X$ vanishes, 
the representation factors through $SO[H^2(X,\RealNumbers)]$.
The representation $\monrep$, given 
in (\ref{eq-mon-representation-of-stabilizer}),
is related to Verbitsky's. 
The construction of the representation $\monrep$
in \cite{markman-monodromy-I} was independent of Verbitsky's, but
we used Verbitsky's result to prove 
that the image of $\monrep$ is a normal subgroup of 
the whole monodromy group $Mon(\M)$. 
We will further need the following proposition, part
\ref{prop-item-so-action-determines-Hodge-structure} 
of which is due to Verbitsky.

\begin{prop}
\label{prop-so-action-determines-Hodge-structure}
\begin{enumerate}
\item 
\label{prop-item-ad-I}
\label{prop-item-so-action-determines-Hodge-structure}
\cite{verbitsky,looijenga-lunts}
Let $X$ be an irreducible holomorphic symplectic manifold,
and $I$ the complex structure of $X$. 
Denote by $ad_I$ the semisimple endomorphism of $H^*(X,\ComplexNumbers)$,
with $H^{p,q}(X)$ an eigenspace with eigenvalue $\sqrt{-1}(p-q)$. 
Then $ad_I$ is an element of the Verbitsky-representation of the 
Lie-algebra $\LieAlg{so}[H^2(X,\ComplexNumbers)]$ on
$H^*(X,\ComplexNumbers)$. 
\item
\label{prop-item-N-is-Zariski-dense}
(\cite{markman-monodromy-I} Lemmas 4.13)
$N$ intersects the image of $SO[H^2(X,\RealNumbers)]$,
via Verbitsky's representation, in a subgroup of $N$ of finite index,
which has finite index in the image of $SO^+[H^2(X,\Integers)]$ as well.
\item
\label{prop-item-N-sub-reps-are-sub-Hodge-str}
Any $N$-subrepresentation, in the tensor product of 
copies of the cohomology groups $H^i(X,\ComplexNumbers)$ or their duals,
is a rational sub-Hodge-structure. 
\item
\label{prop-item-semi-simple}
The Zariski closure of $Mon(X)$ in
$GL[H^*(X,\ComplexNumbers)]$ is semi-simple. Its identity component is 
equal to the image of $SO[H^2(X,\ComplexNumbers)]$
via Verbitsky's representation and is
isomorphic to $SO[H^2(X,\ComplexNumbers)]$.
\end{enumerate}
\end{prop}

\noindent
{\bf Proof: }
\ref{prop-item-N-sub-reps-are-sub-Hodge-str}):
Let $T$ be the tensor product in the statement.
Any $N$-subrepresentation $V$ of $T$ is also a complex 
$SO^+(H^2(X,\RationalNumbers))$ representation, by
part \ref{prop-item-N-is-Zariski-dense} of the proposition
and the fact that any finite index subgroup of $SO^+[H^2(X,\Integers)]$ 
is Zariski dense in $SO[H^2(X,\ComplexNumbers)]$. All 
$SO^+(H^2(X,\RationalNumbers))$ representations are defined over 
$\RationalNumbers$. $V$ is invariant under 
$ad_I$, by part \ref{prop-item-ad-I} of the proposition. 
Hence, $V$ is a sub-Hodge-structure.

\ref{prop-item-semi-simple}) 
$Mon(X)$ is shown to be contained in the subgroup of
$GL[H^*(X,\RealNumbers)]$ generated by the image of $SO[H^2(X,\RealNumbers)]$,
via Verbitsky's representation, the operator $D_X$, acting by
$(-1)^i$ on $H^{2i}(X)$, and the subgroup $K^2$ in part 
\ref{thm-item-markman-monodromy-I-lemma-4.4} of Theorem 
\ref{thm-summary-of-monodromy-results} above, 
which is finite central of exponent $2$ 
(see the proof of Lemma 4.13 in \cite{markman-monodromy-I}). 
The characterization of the identity component of the
Zariski closure follows from 
part \ref{prop-item-N-is-Zariski-dense}.
Semi-simplicity follows from the fact that $K^2$ is finite and central.
\EndProof

We will use repeatedly in the proof of Theorem
\ref{thm-introduction-invariant-Q-4} the following easy lemma.
\begin{new-lemma}
\label{lemma-subrepresentations-with-respect-to-a-normalsubgroup}
Let $G$ be a group, $N<G$ a normal subgroup, $k$ a field, 
$V$ a $k$-vector space, $\rho:G\rightarrow GL(V)$ a
linear representation, and $U<V$  an irreducible 
$N$-subrepresentation. Assume, further, that $\Hom_N(U,V)$ is $1$-dimensional,
and $U$ is the unique 
irreducible $N$-subrepresentation of $V$ of dimension equal to
$\dim(U)$. Then $U$ is also a $G$-subrepresentation.
\end{new-lemma}

\subsection{Proof of Theorems \ref{thm-introduction-invariant-Q-4}
and \ref{thm-introduction-arithmetic-invariant-of-moduli}
and Proposition \ref{prop-mon-is-K2-times-W}}
\label{sec-proof-of-thm-arithmetic-invariant}
{\bf Proof of Theorem \ref{thm-introduction-invariant-Q-4}:}
Part \ref{thm-item-second-chern-class})
$Q^{2i}(X,\Integers)$ is torsion free, since 
$\varphi^{2i}$ is an isomorphism, when $X=\M_H(v)$, by Proposition
\ref{prop-varphi-k-is-independent-of-universal-sheaf}. 
If $i$ is even, the rest of part 
\ref{thm-item-second-chern-class}
follows from Theorem \ref{thm-summary-of-monodromy-results}
part \ref{thm-item-identification-of-residue-of-chern-class}.

The lattice $K(S)$ contains two irreducible integral 
subrepresentations of $OK(S)_v$, 
namely ${\rm span}_\Integers\{v\}$ and $v^\perp$. 
$Ad_{\varphi^{2i}}(N^{2i})=OK(S)_v$, if $i$ is even. If $i$ is odd, then 
$Ad_{\varphi^{2i}}(N^{2i})$ is the subgroup $ext[\W(v^\perp)]$ of
$OK(S)$ 
given in (\ref{eq-commutative-diagram-of-monodromy-groups-hlf-integer-case}).
In both cases, $-1$ and $Ad_{\varphi^{2i}}(N^{2i})$ generate the 
subgroup $OK(S)_{\pm v}$ of $OK(S)$ stabilizing ${\rm span}\{v\}$. 
It follows that $Q^{2i}(X,\Integers)$ consists of two irreducible 
subrepresentations of $N$;
a one dimensional representation $Q^{2i}(X,\Integers)'$ 
and a $23$-dimensional representation $Q^{2i}(X,\Integers)''$. 
Both $Q^{2i}(X,\Integers)'$ and $Q^{2i}(X,\Integers)''$ 
are also $Mon(X)$-subrepresentations, by
Lemma \ref{lemma-subrepresentations-with-respect-to-a-normalsubgroup}.

Part \ref{thm-item-unimodular-pairing}) {\bf Uniqueness}:
Let us assume the existence of a bilinear form on $Q^{2i}(X,\Integers)$,
satisfying the properties in part \ref{thm-item-unimodular-pairing},
and prove its uniqueness. We may assume $X=\M_H(v)$, for some
$S$, $v$, $H$, as in Theorem \ref{thm-irreducibility},
since the bilinear form is assumed monodromy-invariant. 
The form is invariant under the subgroup 
$\langle-1,N^{2i}\rangle$ generated by $-1$ and $N^{2i}$,
being bilinear and monodromy invariant.
The equality $Ad_{\varphi^{2i}}(\langle-1,N^{2i}\rangle)=OK(S)_{\pm v}$ 
implies, that it suffices to prove the following uniqueness statement:

\noindent
{\em The Mukai pairing on $K(S)$ is the unique integral 
even unimodular symmetric bilinear form on $K(S)$, 
which is invariant under $OK(S)_v$ and satisfies the equality
$(v,v)=\dim_\ComplexNumbers(\M_H(v))-2$. 
}

\noindent
Let $B'$ be any such form and denote Mukai's bilinear form by $B$. 
The  $OK(S)_v$-invariance implies, that $v$ and 
its Mukai-pairing-orthogonal sublattice $v^\perp$ 
are orthogonal also with respect to $B'$. 
Furthermore, the restriction $\restricted{B}{v^\perp}'$, of $B'$ to 
$v^\perp$, is equal to an integer multiple of $\restricted{B}{v^\perp}$. 
Unimodularity of $B'$ determines $\restricted{B}{v^\perp}'$ up to sign. 
Let $g\in GL(K(S))$ be the composition 
$K(S)\RightArrowOf{B'} K(S)^* \RightArrowOf{B^{-1}} K(S)$.
Then $g$ is $OK(S)_v$-equivariant, $g(v^\perp)=v^\perp$, 
the restriction $f$ of $g$ to $v^\perp$ is $\pm id$, 
$g(v)=v$, and so $g$ induces the identity automorphism
of $K(S)/v^\perp$. If $f=-id$, then the order of the extension class of
(\ref{eq-extension-ofspan-v-dual-by-v-perp}) is at most $2$. 
The order is $(v,v)$, by Lemma
\ref{lemma-orders-of-extension-classes-with-Mukai-lattice-in-the-middle}, 
and is thus $\geq 6$. Hence $f=id$ and so $g=id$ and 
$B'$ is equal to the Mukai pairing.

Part \ref{thm-item-unimodular-pairing}) {\bf Existence}:
It suffices to prove the existence of a $Mon(X)$-invariant such bilinear 
form, on one $X$ deformation equivalent to $S^{[n]}$, for each $n\geq 4$.
Let $X=\M_H(v)$ and  $B$ the pushforward of the Mukai pairing 
via $\varphi^{2i}:K(S)\rightarrow Q^{2i}(\M_H(v),\Integers)$. 
The form $B$ satisfies 
$B(\varphi^{2i}(v),\varphi^{2i}(v))=(v,v)=\dim_\ComplexNumbers(\M_H(v))-2$, by
the dimension formula 
given in Theorem \ref{thm-irreducibility}. 
When $i$ is even, the above equality translates to 
equation (\ref{eq-chern-classes-encode-dimension}), by
Theorem \ref{thm-summary-of-monodromy-results}
part \ref{thm-item-identification-of-residue-of-chern-class}.
We prove the $Mon(\M_H(v))$-invariance of $B$.
Let $g\in Mon^{2i}(\M_H(v))$ and $f\in N^{2i}$. Then
\[
B(gf(x),gf(y)) \ = \ B((gfg^{-1})g(x),(gfg^{-1})g(y)) \ = \
B(g(x),g(y)),
\]
where the second equality follows from 
the equality $Ad_{\varphi^{2i}}(\langle-1,N^{2i}\rangle)=OK(S)_{\pm v}$ 
and
the fact that $N^{2i}$ is a normal 
subgroup of $Mon^{2i}$ (Theorem \ref{thm-summary-of-monodromy-results} part 
\ref{thm-item-markman-monodromy-I-thm-1.5}). Set $(g_*B)(x,y):=B(g(x),g(y))$.
Then $g_*B$ is an $N^{2i}$-invariant unimodular symmetric bilinear form on
$Q^{2i}(\M_H(v),\Integers)$, which satisfied 
$g_*B(\varphi^{2i}(v),\varphi^{2i}(v))=(v,v)$.
The equality $B=g_*B$ now follows from 
the uniqueness of such a form, proven above. 

Part \ref{thm-intro-item-Q-4-is-the-mukai-lattice}) Follows from
the proof of existence above.

Part \ref{thm-item-pair-of-embeddings})
Assume first that $X=\M_H(v)$. 
Set 
$e:=\varphi^{2i}\circ (\theta_v)^{-1}:H^2(\M_H(v),\Integers)\rightarrow
Q^{2i}(\M_H(v),\Integers)$. Then $e$ spans in 
$\Hom[H^2(\M_H(v),\Integers),Q^{2i}(\M_H(v),\Integers)]$
a rank $1$ integral $N$-subrepresentation,
by equations (\ref{eq-pullback-of-Q-2i-to-stabilizer}) and
(\ref{eq-pullback-of-Q-2-to-stabilizer}). 
These equations also show, that the character ${\rm span}\{e\}$
of $N$ pulls back to the character $\tilde{\eta}^{\otimes i-1}$
of $OK(S)_v$, via the representation $\monrep$ in 
(\ref{eq-mon-representation-of-stabilizer}). 
The character $\tilde{\eta}^{\otimes i-1}$ is non-trivial if $i$ is even. 
Furthermore, $e$ is an isometry onto $\bar{c}_{i}(\M_H(v))^\perp$. 
The $N$-representation 
$\Hom[H^2(\M_H(v),\Integers),\bar{c}_{i}(\M_H(v))^\perp]$ 
decomposes over $\RationalNumbers$
into three irreducible $N$-representations, of three different ranks, 
corresponding to the decomposition of the tensor square of the standard 
irreducible representation of $SO(23,\RationalNumbers)$.
Thus, ${\rm span}\{e\}$ is also a $Mon(\M_H(v))$-subrepresentations, by
Lemma \ref{lemma-subrepresentations-with-respect-to-a-normalsubgroup}.
We get a $Mon(X)$-invariant pair $\{e,-e\}$ of such 
primitive integral embeddings,
for any irreducible holomorphic symplectic manifold deformation equivalent to 
$\M_H(v)$. 
${\rm Span}\{e\}$ is a Hodge substructure of
$\Hom[H^2(X,\Integers),Q^{2i}(X,\Integers)]$, by 
Proposition \ref{prop-so-action-determines-Hodge-structure} 
part \ref{prop-item-N-sub-reps-are-sub-Hodge-str}. 
%
This completes the proof of Theorem \ref{thm-introduction-invariant-Q-4}.
\EndProof

\medskip
{\bf Proof of Theorem 
\ref{thm-introduction-arithmetic-invariant-of-moduli}:}
Part a)
Equation (\ref{eq-varphi-4-maps-v-to-a-multiple-of-second-Chern-class})
follows from Theorem 
\ref{thm-summary-of-monodromy-results} part 
\ref{thm-item-identification-of-residue-of-chern-class}.
Equation (\ref{eq-e-is-conjugated-to-iota})
is the definition of $e$ above, but we need its proof to be independent of
the representation of $X$ as a moduli space of sheaves on a $K3$ surface,
which may not be unique.
Equation (\ref{eq-e-is-conjugated-to-iota})
follows from the equalities
$Ad_{\theta_v}(N^2)=\W(v^\perp)$ and 
$Ad_{\varphi^{2i}}(N^{2i})=OK(S)_v$ in 
diagram (\ref{eq-commutative-diagram-of-monodromy-groups})
by the following argument. 
The two equalities imply, that we can identify $N$ 
with both $\W(v^\perp)$ and $OK(S)_v$. 
Using this identification, each side of equation 
(\ref{eq-e-is-conjugated-to-iota}) is an $N$-equivariant embedding of 
$v^\perp$ in $Q^{2i}(\M_H(v),\Integers)$. The
multiplicity of the $N$-representation $v^\perp$ in $Q^{2i}(\M_H(v),\Integers)$
is one. Thus, the embedding is unique, up to sign, by Schur's Lemma. 
The equality (\ref{eq-e-is-conjugated-to-iota}) follows. 
The definition of the bilinear form on $Q^{2i}(\M_H(v),\Integers)$ implies, 
that $\varphi^{2i}$ is an isometry, and $\varphi^2\circ \iota$ is equal to
$\theta_v$, which is a Hodge-isometry by Theorem
\ref{thm-irreducibility}. It remains to prove that the isometry $\varphi^{2i}$ 
is compatible with Hodge structures. This is clear 
when a universal sheaf exists over
$S\times \M_H(v)$. In the absence of a universal sheaf, equations
(\ref{eq-varphi-4-maps-v-to-a-multiple-of-second-Chern-class}) and 
(\ref{eq-e-is-conjugated-to-iota})
imply that $\varphi^{2i}$ is a Hodge-isometry, as
both $e$ and $\theta_v$ are. 
The proof of part b is similar.
\EndProof


\medskip
{\bf Proof of Proposition \ref{prop-mon-is-K2-times-W}:}
$Mon^2=\W H^2(X,\Integers)$, by Theorem \ref{thm-introduction-Mon-2-is-W}. 
Let $N$ be the normal subgroup of $Mon(X)$, given in (\ref{eq-N}).
$N$ surjects onto $\W H^2(X,\Integers)$, by Theorem
\ref{thm-summary-of-monodromy-results} part \ref{thm-item-N-2-is-W}.
Thus $Mon(X)$ is generated by $K^2$ and $N$. 
$K^2$ is finite, central, and of exponent $2$, by 
Theorem \ref{thm-summary-of-monodromy-results} part 
\ref{thm-item-markman-monodromy-I-lemma-4.4}. 
$K^2$ and $N$ are both normal subgroups, which intersect trivially
and commute, so $Mon(X)$ is their direct product. 
$N$ is equal to the image of $\nu$, by Lemma
\ref{lemma-mu-is-an-isomorphism}.
\EndProof

%
\section{Period maps and monodromy groups}
\label{sec-degree-of-period-maps}

In section \ref{sec-symmetry-groups} we compare the isometry group of
$H^2(X,\Integers)$, with the subgroup stabilizing $\bar{c}_2(X)$ in the
isometry group of $Q^{4}(X,\Integers)$. 
Counter examples to the weight $2$ Generic Torelli question are 
provided in section 
\ref{sec-period-maps} (Theorem \ref{thm-non-bimeromorphic-classes}). 
The proof of Theorem \ref{thm-non-bimeromorphic-classes}
uses the construction of a {\em refined period map}, 
defined below in (\ref{eq-refined-period-map}).
The moduli space, of marked irreducible holomorphic symplectic manifolds
deformation equivalent to $S^{[n]}$, is disconnected in general,
by the proof of Corollary \ref{cor-introduction-counter-example-to-Torelli}.
The refined period map extracts an intrinsic invariant of a 
connected component.
In section \ref{sec-non-birational-moduli-spaces} we provide
explicit examples of non-birational pairs of projective 
moduli spaces of stable sheaves
on a $K3$ surface $S$, with Hodge-isometric weight $2$ cohomologies.

\subsection{Symmetry groups}
\label{sec-symmetry-groups}

Fix a pair of lattices $\Lambda\subset \widetilde{\Lambda}$ isometric to 
the pair $v^\perp\subset K(S)$, with respect to the Mukai pairing,
where $v\in K(S)$ is a primitive class and $(v,v)\geq 2$. 
Let $w$ be a generator of the rank $1$ sublattice $\Lambda^\perp$. 
It satisfies $(w,w)=(v,v)$. 
The symmetry group of the pair 
$(\widetilde{\Lambda},w)$ is the subgroup $O\widetilde{\Lambda}_w$ 
of $O\widetilde{\Lambda}$ stabilizing $w$. 
There is a natural homomorphism 
\begin{equation}
\label{eq-extension-homomorphism}
\mu \ : \ O\widetilde{\Lambda}_w \ \ \hookrightarrow \ \  O^+\Lambda
\end{equation}
sending $g\in O\widetilde{\Lambda}_w$ to the restriction 
of $\tilde{\eta}(g)\cdot g$ to $\Lambda$. 
Above, $\tilde{\eta}$ is the orientation character 
(\ref{eq-tilde-eta}) of $O\widetilde{\Lambda}$.
The homomorphism $\mu$ arises naturally in the context of monodromy 
representations (see diagram
(\ref{eq-commutative-diagram-of-monodromy-groups})).

The Hodge structure $Q^4(X,\Integers)$, 
together with the symmetric bilinear form and the 
class $\bar{c}_2(X)$, is a more refined invariant  
than the weight $2$ Hodge structure $H^2(X,\Integers)$. 
The refinement is by discrete and finite additional data. 
The period domains, for the two Hodge structures, are identical, but
the arithmetic symmetry group of the pair 
$\{Q^4,\bar{c}_2\}(X)$, forgetting its Hodge structure, 
is smaller in general than that of $H^2(X,\Integers)$. 
The arithmetic symmetry group $O\widetilde{\Lambda}_w$ 
of the former, is embedded in the symmetry group $O^+\Lambda$ of
the latter, via the homomorphism $\mu$ given in 
(\ref{eq-extension-homomorphism}).
Lemma \ref{lemma-index-of-smaller-symmetry-group} states, 
that $O\widetilde{\Lambda}_w$ is often a {\em proper subgroup} 
of $O^+\Lambda$. 
Let $(w,w)/2=p_1^{e_1}p_2^{e_2}\cdots p_r^{e_r}$ be the prime factorization
of $(w,w)/2$, where $p_1<p_2<\cdots <p_r$ are distinct primes, and $e_i$ 
are positive integers. Set $r=0$, if $(w,w)=2$. 

\begin{new-lemma}
\label{lemma-index-of-smaller-symmetry-group}
The image of $O\widetilde{\Lambda}_w$, via the homomorphism $\mu$ given in
(\ref{eq-extension-homomorphism}), is a normal subgroup of 
$O^+\Lambda$. The quotient 
$O^+\Lambda/\mu(O\widetilde{\Lambda}_w)$ is isomorphic to
$(\Integers/2\Integers)^d$, where 
\[
d \ = \ \left\{
\begin{array}{lcl}
0 & \mbox{if} &  (w,w)=2,
\\
r-1 & \mbox{if} & (w,w)\geq 4.
\end{array}
\right.
\]
In particular, $\mu$ surjects onto $O^+\Lambda$, if and only if 
$(w,w)=2$ or $(w,w)/2$ is a prime power.
\end{new-lemma}

The proof will depend on the following Lemma.
The bilinear form embeds the dual lattice $\Lambda^*$ in 
$\Lambda\otimes\RationalNumbers$. 
Since the lattice $\Lambda$ is even, we get a residual quadratic form 
\[
q \ : \ \Lambda^*/\Lambda \ \ \longrightarrow \ \ 
\RationalNumbers/2\Integers. 
\]
Denote by $O(\Lambda^*/\Lambda)$ the subgroup of $GL(\Lambda^*/\Lambda)$
leaving $q$ invariant. 
The following Lemma provides three additional characterizations of 
the image of $\mu$. 
Consider the natural homomorphism
\[
\pi \ : \ O\Lambda \ \ \ \longrightarrow \ \ \ 
O(\Lambda^*/\Lambda). 
\]

\begin{new-lemma}
\label{lemma-on-residual-orthogonal-group}
\begin{enumerate}
\item
\label{lemma-item-three-characterizations-of-W}
The following four subgroups of $O^+(\Lambda)$ are equal:
\begin{enumerate}
\item
\label{lemma-equal-subgroups-image-of-mu}
The image of $O\widetilde{\Lambda}_w$ via the homomorphism $\mu$ 
given in (\ref{eq-extension-homomorphism}).
\item
\label{lemma-equal-subgroups-W}
The subgroup
$\W(\Lambda)$ of $O^+\Lambda$ given in (\ref{eq-W}).
\item
\label{lemma-item-W-is-inverse-image-via-pi}
The intersection of $O^+\Lambda$ with the inverse image 
$\pi^{-1}\{1,-1\}$. 
\item
The subgroup of orientation preserving isometries, which can be extended to 
isometries of $\widetilde{\Lambda}$.
\end{enumerate}
\item
\label{lemma-item-mu-is-injective}
The homomorphism $\mu$ is injective, provided $(w,w)\geq 4$. 
\item
\label{lemma-item-pi-is-surjective}
The homomorphism $\pi$ maps $O^+\Lambda$ onto $O(\Lambda^*/\Lambda)$.
\item
\label{lemma-item-oder-of-residual-orthogonal-group}
The quotient group $O(\Lambda^*/\Lambda)$ is
isomorphic to $(\Integers/2\Integers)^r$. 
\end{enumerate}
\end{new-lemma}

\noindent
{\bf Proof:}
\ref{lemma-item-three-characterizations-of-W}) 
Part \ref{lemma-item-three-characterizations-of-W} is proven in
Lemma 4.10 in \cite{markman-monodromy-I}.

\noindent
\ref{lemma-item-mu-is-injective}) This part follows from  Lemma
\ref{lemma-mu-is-an-isomorphism}.

\noindent
\ref{lemma-item-pi-is-surjective}) Part \ref{lemma-item-pi-is-surjective}
follows from Theorem 1.14.2 in \cite{nikulin}.

\noindent
\ref{lemma-item-oder-of-residual-orthogonal-group})
This part of the Lemma is proven in \cite{oguiso}. 
\EndProof

\medskip
{\bf Proof of Lemma 
\ref{lemma-index-of-smaller-symmetry-group}:}
Parts \ref{lemma-item-mu-is-injective} and 
\ref{lemma-item-pi-is-surjective} of Lemma 
\ref{lemma-on-residual-orthogonal-group}
yield the commutative diagram with a short exact sequence at the
top horizontal row 
\begin{equation}
\label{eq-short-exact-at-the-top-horizontal-sequence}
\begin{array}{ccccclc}
0  \rightarrow & O\widetilde{\Lambda}_w & \LongRightArrowOf{\mu} &
O^+\Lambda & \longrightarrow & O(\Lambda^*/\Lambda)/\{1,-1\} & 
\rightarrow 0
\\
& & \cong \ \searrow \ \hspace{1em} & \cup & \hspace {1ex} \ \searrow \ \pi & 
\hspace{4ex} \uparrow 
\\
& & & \W(\Lambda) & & O(\Lambda^*/\Lambda).
\end{array}
\end{equation}
Lemma 
\ref{lemma-index-of-smaller-symmetry-group} follows from 
Lemma \ref{lemma-on-residual-orthogonal-group} part 
\ref{lemma-item-oder-of-residual-orthogonal-group} and the exactness
of the top horizontal row of the above diagram.
\EndProof

\medskip
{\bf Proof of Theorem 
\ref{thm-invariant-Q-bar-c-captures-all}:}
Assume that $n\geq 4$, $g:Q^4(X_1,\Integers)\rightarrow Q^4(X_2,\Integers)$
is an isometry, and $g(\bar{c}_2(X_1))=\bar{c}_2(X_2)$.
Then $g$ is induced by a parallel-transport operator 
$f:H^*(X_1,\Integers)\rightarrow H^*(X_2,\Integers)$,
by Theorem \ref{thm-introduction-Mon-2-is-W},
Theorem \ref{thm-introduction-invariant-Q-4}, and the equality of 
$\W[H^2(X_2,\Integers)]$ with the image of the stabilizer 
$O[Q^4(X_2,\Integers)]_{\bar{c}_2}$ in $O[H^2(X_2,\Integers)]$
(the two subgroups \ref{lemma-equal-subgroups-image-of-mu} and 
\ref{lemma-equal-subgroups-W} in 
Lemma \ref{lemma-on-residual-orthogonal-group}). 
Assume, further, that $g$ is a Hodge-isometry. 
Then the restriction $f_2$, of $f$ to $H^2(X_1)$,
is a Hodge-isometry, by Theorem \ref{thm-introduction-invariant-Q-4}. 
Let $ad_{I_{X_i}}\in \LieAlg{gl}[H^*(X_i,\ComplexNumbers)]$, $i=1,2$, be the 
graded endomorphism in Proposition 
\ref{prop-so-action-determines-Hodge-structure},
inducing the Hodge-decomposition.
Denote its restriction to $H^2(X_i,\ComplexNumbers)$
by $\overline{ad}_{I_{X_i}}$. Verbitsky's representation
$\rho_{X_i}:\LieAlg{so}[H^2(X_i,\ComplexNumbers)]\rightarrow 
\LieAlg{gl}[H^*(X_i,\ComplexNumbers)]$, given in
Proposition 
\ref{prop-so-action-determines-Hodge-structure}, 
takes $\overline{ad}_{I_{X_i}}$ back to $ad_{I_{X_i}}$. 
Being a parallel-transport operator, $f$ conjugates $\rho_{X_1}$ to
$\rho_{X_2}$. 
Hence
$f\circ ad_{I_{X_1}}\circ f^{-1}=
f\circ\rho_{X_1}[\overline{ad}_{I_{X_1}}]\circ f^{-1}
=\rho_{X_2}[f_2\circ \overline{ad}_{I_{X_1}} \circ f_2^{-1}] =
\rho_{X_2}(\overline{ad}_{I_{X_2}})=ad_{I_{X_2}}$.
Thus, $f$ is compatible with all the Hodge structures.
The proof of the case $n\leq 3$ is similar.
\EndProof

\medskip
The statement of Lemma 
\ref{lemma-the-index-of-W-versus-orbits-in-hyperbolic-plane} below will use 
Mukai's notation for the Mukai lattice, which we recall next.
Identify the group $K(S)$ with $H^*(S,\Integers)$, via the 
isomorphism  sending a class $E$ to its {\em Mukai vector} 
$ch(E)\sqrt{td_S}$. Using the grading of $H^*(S,\Integers)$, 
the Mukai vector is 
\begin{equation}
\label{eq-Mukai-vector}
(\rank(E),c_1(E),\chi(E)-\rank(E)),
\end{equation} 
where the rank is considered in
$H^0$ and $\chi(E)-\rank(E)$ in $H^4$ via multiplication by the 
orientation class of $S$. The homomorphism 
$ch(\bullet)\sqrt{td_S}:K(S)\rightarrow H^*(S,\Integers)$ 
is an isometry with respect 
to the Mukai pairing (\ref{eq-mukai-pairing-on-K-top}) on $K(S)$
and the pairing 
\[
\left((r',c',s'),(r'',c'',s'')\right) \ \ = \ \ 
\int_{S}c'\cup c'' -r'\cup s''-s'\cup r''
\]
on $H^*(S,\Integers)$
(by the Hirzebruch-Riemann-Roch Theorem). 
The Mukai vector in $H^*(S,\Integers)$, of the
ideal sheaf of a length $n$ subscheme, is $v:=(1,0,1-n)$.

Let $\widetilde{\Lambda}$ be the Mukai lattice $H^*(S,\Integers)$,
forgetting its Hodge structure, and 
set $\Lambda:=(1,0,1-n)^\perp$, where $n\geq 2$. 
Let $O(\Lambda,\widetilde{\Lambda})$ be the set of primitive 
isometric embeddings of $\Lambda$ into $\widetilde{\Lambda}$. 
$O(\Lambda,\widetilde{\Lambda})$  is endowed with a left
$O(\widetilde{\Lambda})$-action and a right $O(\Lambda)$-action.
The following Lemma 
\ref{lemma-the-index-of-W-versus-orbits-in-hyperbolic-plane}
provides a set of representatives for 
the $O(\widetilde{\Lambda})$-orbits in $O(\Lambda,\widetilde{\Lambda})$.
The lemma provides a second calculation for the index of $\W(\Lambda)$ in
$O^+(\Lambda)$. 

Let $\P_n$ be the subset of $\Integers\oplus\Integers$ given by
\begin{equation}
\label{eq-P-n}
\P_n \ \ := \ \ \left\{ (r,s) \ \ : \ \ 
-s\geq r>0, \ \ rs=1-n, \ \ \mbox{and} \ \ \gcd(r,s)=1
\right\}.
\end{equation}
The cardinality of $\P_n$ is clearly $2^{\rho(n-1)-1}$.
Let $U$ be the rank $2$ hyperbolic lattice given by the
pairing $((r_1,s_1),(r_2,s_2))=-r_1s_2-r_2s_1$. 
Then $O(U)$ is isomorphic to 
$\Integers/2\Integers\times \Integers/2\Integers$.
The set $\P_n$ consists of one representative from each
$O(U)$-orbit of  primitive elements of $U$ of square-length $2n-2$. 
For each $(r,s)$ in $\P_n$, let
$\iota_{r,s}:\Lambda\hookrightarrow \widetilde{\Lambda}$ be the isometric 
embedding, which restricts to $H^2(S,\Integers)$ as the identity, 
and sends $(1,0,n-1)$ to $(r,0,-s)$. The image of
$\iota_{r,s}$ is $(r,0,s)^\perp$.

\begin{new-lemma}
\label{lemma-the-index-of-W-versus-orbits-in-hyperbolic-plane}
\begin{enumerate}
\item
\label{lemma-item-P-n-represent-O-Lambda-tilde-orbits}
The map $(r,s)\mapsto \iota_{r,s}$ induces a one-to-one correspondence,
between the set $\P_n$ and the set of $O(\widetilde{\Lambda})$-orbits
in $O(\Lambda,\widetilde{\Lambda})$.
\item
\label{lemma-item-W-is-the-stabilizer-of-action-on-orbit-space}
$O(\widetilde{\Lambda})\times O(\Lambda)$ acts transitively on 
$O(\Lambda,\widetilde{\Lambda})$. The stabilizer 
in $O(\Lambda)$, of every point in the orbit space 
$O(\widetilde{\Lambda})\setminus O(\Lambda,\widetilde{\Lambda})$,
is generated by $\W(\Lambda)$ and $-1$.
\end{enumerate}
\end{new-lemma}

\journal{The Lemma is an easy consequence of Theorem 1.14.4 in
\cite{nikulin}. A detailed proof can be found in the preprint version
\cite{markman-eprint-version}.
}

\preprint{
\noindent
{\bf Proof:} 
\ref{lemma-item-P-n-represent-O-Lambda-tilde-orbits})
Suppose $\iota_{r_2,s_2}=g\circ \iota_{r_1,s_1}$, for some 
$g\in O(\widetilde{\Lambda})$. Then $g$ 
leaves $H^2(S,\Integers)$ invariant and restricts to $H^2(S,\Integers)$
as the identity. Hence, $g$ comes from an isometry of the hyperbolic plane
$U$, which takes $(r_1,-s_1)$ to $(r_2,-s_2)$. Consequently,
$(r_1,s_1)$ and $(r_2,s_2)$ also belong to the same $O(U)$ orbit.
Each being the unique representative in $\P_n$, we conclude 
the equality $(r_1,s_1)=(r_2,s_2)$. 

Let $\iota:\Lambda\hookrightarrow \widetilde{\Lambda}$ be an isometric
embedding. There is a unique $O(\widetilde{\Lambda})$-orbit of 
isometric embeddings of the unimodular lattice $H^2(S,\Integers)$ 
in $\widetilde{\Lambda}$ (\cite{nikulin} Theorem 1.14.4). 
Hence, there is an isometry 
$g\in O(\widetilde{\Lambda})$, such that $g\circ \iota$ restricts to
$H^2(S,\Integers)$ as the inclusion. Then $g(\iota(1,0,n-1))$ 
belongs to the orthogonal complement $H^2(S,\Integers)^\perp$ in
$\widetilde{\Lambda}$. Now $H^2(S,\Integers)^\perp$ is the hyperbolic plane
and $\widetilde{\Lambda}$ is an orthogonal direct sum 
$H^2(S,\Integers)\oplus H^2(S,\Integers)^\perp$.
Composing $g\circ \iota$ with an isometry of $H^2(S,\Integers)^\perp$, 
we get $\iota_{r,s}$, for some $(r,s)\in \P_n$. 

\ref{lemma-item-W-is-the-stabilizer-of-action-on-orbit-space})
Let $\iota:\Lambda\hookrightarrow \widetilde{\Lambda}$ be a primitive
isometric embedding. Then the orthogonal complement
$\iota(\Lambda)^\perp$, of its image, is generated by an element of 
square-length $2n-2$. There is a unique $O(\widetilde{\Lambda})$-orbit
of such primitive elements (\cite{nikulin} Theorem 1.14.4).
The transitivity of the $O(\Lambda)\times O(\widetilde{\Lambda})$ action 
follows. 
The description of the stabilizer follows from the
equality of $\W(\Lambda)$ and the image of the homomorphism
$\mu$ given in (\ref{eq-extension-homomorphism}) 
(see Lemma \ref{lemma-on-residual-orthogonal-group}).
\EndProof
}

\subsection{The refined period map}
\label{sec-period-maps}
We exhibit pairs of irreducible holomorphic symplectic
manifolds $X$ and $Y$, in the deformation class of $S^{[n]}$,
with isomorphic weight $2$ Hodge structures, 
but for which 
$\{Q^4,\bar{c}_2\}(X)$ and 
$\{Q^4,\bar{c}_2\}(Y)$ are {\em not} Hodge-isometric. 

\begin{condition}
\label{cond-stabilizer-of-the-period-is-contained-in-W}
The subgroup of $O^+H^2(X,\Integers)$, preserving the
Hodge structure, is contained in the subgroup $\W$
given in (\ref{eq-W}).
\end{condition}


Let $\Lambda$ be a fixed lattice isometric to $H^2(S^{[n]},\Integers)$, 
$n\geq 4$, and $\Omega$ the period domain given in equation
(\ref{eq-period-domain}). 
The group $O\Lambda$ naturally acts on $\Omega$. 
The index $2$ subgroup $O^+\Lambda$ acts faithfully. 
Any period $\ell\in\Omega$, in the complement of the (countable)
union of eigenspaces of isometries $g$ in $O\Lambda\setminus\{id,-id\}$,
satisfies Condition \ref{cond-stabilizer-of-the-period-is-contained-in-W},
since the stabilizer of $\ell$ in $O^+\Lambda$ is trivial. 
A criterion, for the Hilbert scheme $S^{[n]}$ to satisfy
condition \ref{cond-stabilizer-of-the-period-is-contained-in-W},
is provided in Lemma \ref{lemma-S-H-n-satisfying-condition}.

Let $r$ be the number of prime powers in the prime decomposition of $n-1$
as in Lemma \ref{lemma-index-of-smaller-symmetry-group}. 

\begin{thm}
\label{thm-non-bimeromorphic-classes}
Assume that $n\geq 4$. Fix a weight $2$ Hodge structure $h$
in $\Omega$ satisfying Condition 
\ref{cond-stabilizer-of-the-period-is-contained-in-W}. 
There is a set of $2^{r-1}$ pairwise non-bimeromorphic 
irreducible holomorphic symplectic manifolds, 
deformation equivalent to $S^{[n]}$, 
each with a weight $2$ Hodge structure Hodge-isometric to $h$.
Furthermore, the set may be chosen with invariants
$\{Q^4,\bar{c}_2\}(X)$, which are pairwise not Hodge-isometric.
\end{thm}

The proof of the Theorem will depend on 
Lemma \ref{lemma-non-bimeromorphic-classes} below.
Fix a lattice $\widetilde{\Lambda}$, isometric to the Mukai lattice 
$K(S)$ of a $K3$ surface, forgetting its Hodge structure.
Let $O(\Lambda,\widetilde{\Lambda})$ be the set of primitive 
isometric embeddings of $\Lambda$ into $\widetilde{\Lambda}$, introduced in 
Lemma \ref{lemma-the-index-of-W-versus-orbits-in-hyperbolic-plane}.
Let ${\frak M}_\Lambda$ be the moduli space of marked 
irreducible holomorphic symplectic manifolds, 
deformation equivalent to $S^{[n]}$, with $n\geq 4$.
Theorem \ref{thm-introduction-invariant-Q-4} enables us 
to refine the period map (\ref{eq-period-map}) 
and define 
\begin{equation}
\label{eq-refined-period-map}
\widetilde{P}  \ : \ {\frak M}_\Lambda \ \ \ \longrightarrow \ \ \ 
\Omega\times [O(\widetilde{\Lambda})\setminus O(\Lambda,\widetilde{\Lambda})].
\end{equation}
The first element in the pair $\widetilde{P}(X,\phi)$ 
is the period $P(X,\phi)$, given in (\ref{eq-period-map}).
The second element in the pair $\widetilde{P}(X,\phi)$ is 
defined as follows. Compose 
the inverse, of the marking $\phi$, with the pair of embeddings $\pm e_X$
given in Theorem \ref{thm-introduction-invariant-Q-4}, to determine
a pair of embeddings of $\Lambda$ into $Q^4(X,\Integers)$. 
The pair $\pm e_X\circ\phi^{-1}$ determines a well defined 
$O(\widetilde{\Lambda})$ orbit in $O(\Lambda,\widetilde{\Lambda})$.

Let ${\frak M}^0_\Lambda$ be a connected component of ${\frak M}_\Lambda$, 
$(X,\phi)$ a  marked pair in ${\frak M}^0_\Lambda$, and $f$ 
an element of $O^+\Lambda$. Choose a marked pair $(Y,\psi)$ 
in ${\frak M}^0_\Lambda$, whose period satisfies 
$
P(Y,\psi) = f(P(X,\phi)).
$
Such a pair $(Y,\psi)$ exists, by the Surjectivity Theorem
\cite{huybrects-basic-results}. 

\begin{new-lemma}
\label{lemma-non-bimeromorphic-classes}
Assume that $X$ satisfies Condition
\ref{cond-stabilizer-of-the-period-is-contained-in-W}.
If $f$ does not belong to $\W(\Lambda)$, 
then the pairs
$\{Q^4,\bar{c}_2\}(X)$ and 
$\{Q^4,\bar{c}_2\}(Y)$ are {\em not} 
Hodge-isometric. 
\end{new-lemma}

{\bf Proof of Theorem \ref{thm-non-bimeromorphic-classes}:}
The Theorem follows from Lemma \ref{lemma-non-bimeromorphic-classes}
via the computation of the index of $\W(\Lambda)$ in 
$O^+(\Lambda)$, carried out in 
Lemma \ref{lemma-on-residual-orthogonal-group}. 
\EndProof

The proof of Lemma \ref{lemma-non-bimeromorphic-classes}
will depend on Lemma \ref{lemma-equivalent-data}.
$O^+(\Lambda)$ orbits in the period domain $\Omega$ parametrize
Hodge-isometry classes of weight $2$ Hodge structures. 
We characterize in Lemma \ref{lemma-equivalent-data} the data  parametrized 
by $\W(\Lambda)$ orbits in $\Omega$. 
We say that two pairs $(\ell_j,\iota_j)$ in 
$\Omega\times O(\Lambda,\widetilde{\Lambda})$, $j=1,2$, 
are {\em Hodge-isometric}, if there are isometries $g\in O(\Lambda)$
and $f\in O(\widetilde{\Lambda})$, such that
$g(\ell_1)=\ell_2$ and $f\circ\iota_1=\iota_2\circ g$.

\begin{new-lemma}
\label{lemma-equivalent-data}
Data \ref{data-orbits-in-Omega-times-orbit-space-of-embeddings} and 
\ref{data-hodge-isometries-of-embeddings} below are equivalent. 
Upon a choice of an $O(\widetilde{\Lambda})$ orbit in 
$O(\Lambda,\widetilde{\Lambda})$, each of 
data \ref{data-W-orbits-in-Omega} and 
\ref{data-Hodge-isometry-classes-of-Q-v} 
are equivalent to data 
\ref{data-orbits-in-Omega-times-orbit-space-of-embeddings} as well as 
\ref{data-hodge-isometries-of-embeddings}.
\begin{enumerate}
\item
\label{data-orbits-in-Omega-times-orbit-space-of-embeddings}
$O^+(\Lambda)$ orbits in $\Omega\times 
[O(\widetilde{\Lambda})\setminus O(\Lambda,\widetilde{\Lambda})]$.
\item
\label{data-hodge-isometries-of-embeddings}
Hodge isometry classes of pairs $(\ell,\iota)$, consisting of 
a polarized Hodge structure 
$\ell\in \Omega$ on $\Lambda\otimes \ComplexNumbers$, and a 
primitive isometric 
embedding $\iota:\Lambda\hookrightarrow \widetilde{\Lambda}$.
\item
\label{data-W-orbits-in-Omega}
$\W(\Lambda)$ orbits in $\Omega$.
\item
\label{data-Hodge-isometry-classes-of-Q-v}
Hodge isometry classes of pairs $(Q,v)$, consisting of 
an integral polarized Hodge structure $Q$ on 
$\widetilde{\Lambda}\otimes_\Integers\ComplexNumbers$ of weight $2$,
and a primitive integral class $v$ in $\widetilde{\Lambda}$
of Hodge-type $(1,1)$ and square-length $2n-2$.
\end{enumerate}
\end{new-lemma}

\noindent
{\bf Proof:}
($\ref{data-orbits-in-Omega-times-orbit-space-of-embeddings}
\Leftrightarrow$ \ref{data-hodge-isometries-of-embeddings})
The equivalence of Data 
\ref{data-orbits-in-Omega-times-orbit-space-of-embeddings} and 
\ref{data-hodge-isometries-of-embeddings} is clear. 

(\ref{data-orbits-in-Omega-times-orbit-space-of-embeddings}
$\Leftrightarrow$ \ref{data-W-orbits-in-Omega})
$\W(\Lambda)$ is the stabilizer in $O^+(\Lambda)$ of any orbit 
$O(\widetilde{\Lambda})\cdot \iota$ in
$[O(\widetilde{\Lambda})\setminus O(\Lambda,\widetilde{\Lambda})]$, by Lemma 
\ref{lemma-the-index-of-W-versus-orbits-in-hyperbolic-plane}. 
The equivalence of Data 
\ref{data-orbits-in-Omega-times-orbit-space-of-embeddings}
and \ref{data-W-orbits-in-Omega} follows.

It remains to prove the equivalence of  Data
\ref{data-W-orbits-in-Omega} and 
\ref{data-Hodge-isometry-classes-of-Q-v} 
upon a choice of an $O(\widetilde{\Lambda})$ orbit in 
$O(\Lambda,\widetilde{\Lambda})$. 

(\ref{data-Hodge-isometry-classes-of-Q-v} $\Rightarrow$
\ref{data-W-orbits-in-Omega})
Given a pair $(Q,v)$
as in data \ref{data-Hodge-isometry-classes-of-Q-v},
the orthogonal complement $v^\perp$ is a sub-Hodge structure. 
Furthermore, $v^\perp$ is isometric to $\Lambda$, since the group
$O(\widetilde{\Lambda})$ acts transitively on the set of primitive
elements in $\widetilde{\Lambda}$ of square-length $2n-2$. 
Choose an isometry $\iota:\Lambda\rightarrow v^\perp$, in the
prescribed $O(\widetilde{\Lambda})$-orbit in $O(\Lambda,\widetilde{\Lambda})$.
This is possible, since $O(\Lambda)$ acts transitively on
$[O(\widetilde{\Lambda})\setminus O(\Lambda,\widetilde{\Lambda})]$,
by Lemma \ref{lemma-the-index-of-W-versus-orbits-in-hyperbolic-plane}.
The choice of $\iota$ determines a Hodge structure on 
$\Lambda\otimes \ComplexNumbers$. All the possible choices of $\iota$ 
vary in an orbit of the subgroup of $O(\Lambda)$, 
generated by $-1$ and $\W(\Lambda)$, by Lemma
\ref{lemma-the-index-of-W-versus-orbits-in-hyperbolic-plane}. 
Hence, we get a $\W$-orbit in $\Omega$, naturally associated to $(Q,v)$
and the orbit $O(\widetilde{\Lambda})\iota$. 

(\ref{data-W-orbits-in-Omega} $\Rightarrow$
\ref{data-Hodge-isometry-classes-of-Q-v})
Given a Hodge structure $\ell\in \Omega$, representing Data
\ref{data-W-orbits-in-Omega}, choose an isometric
primitive embedding $\iota:\Lambda\hookrightarrow \widetilde{\Lambda}$,
in the prescribed $O(\widetilde{\Lambda})$ orbit,
completing it to Data 
\ref{data-orbits-in-Omega-times-orbit-space-of-embeddings}.
Let $v$ be a generator of the orthogonal complement
$\iota(\Lambda)^\perp$. Let the Hodge structure $Q^{2,0}$ 
be the image $\iota(\ell)$. Then $v$ is of type $(1,1)$ and its
square-length is $2n-2$ by Lemma 
\ref{lemma-the-index-of-W-versus-orbits-in-hyperbolic-plane}. 
The pairs $(Q,v)$ and $(Q,-v)$ are Hodge isometric, via 
$-1\in O(\widetilde{\Lambda})$.
Let $\iota_1$ be another embedding 
in the prescribed $O(\widetilde{\Lambda})$ orbit, 
yielding a 
pair $(Q_1,v_1)$ chosen as above. 
Then $\iota_1=g\circ \iota$, for a unique element 
$g\in O(\widetilde{\Lambda})$, and either
$v_1=g(v)$ or $v_1=-g(v)$, by part
\ref{lemma-item-mu-is-injective} of Lemma 
\ref{lemma-on-residual-orthogonal-group}. 
Consequently, either $g$ or $-g$ is a Hodge isometry 
between $(Q,v)$ and $(Q_1,v_1)$.
If $\ell_1=f(\ell)$, $f\in \W(\Lambda)$, is another Hodge structure in 
$\Omega$, in the $\W(\Lambda)$-orbit of $\ell$, 
then $\iota_1:=\iota\circ f^{-1}$ is an embedding in the prescribed 
$O(\widetilde{\Lambda})$-orbit, by Lemma 
\ref{lemma-the-index-of-W-versus-orbits-in-hyperbolic-plane}.
Hence, the same pair $(Q,v)$ is associated to $\ell_1$.
We get a well defined Hodge isometry class of pairs in data
\ref{data-Hodge-isometry-classes-of-Q-v}, naturally associated to the
$\W(\Lambda)$ orbit of $\ell$ and the point in 
$O(\Lambda,\widetilde{\Lambda})$.
\EndProof

\medskip
{\bf Proof of Lemma \ref{lemma-non-bimeromorphic-classes}:} 
Since $(X,\phi)$ and $(Y,\psi)$ belong to the same
connected component ${\frak M}^0_\Lambda$, then the refined periods 
$\widetilde{P}(X,\phi)$ and $\widetilde{P}(Y,\psi)$
belong to the same connected component of 
$\Omega\times 
[O(\widetilde{\Lambda})\setminus O(\Lambda,\widetilde{\Lambda})]$.
Assume that the pairs $\{Q^4,\bar{c}_2\}(X)$ and 
$\{Q^4,\bar{c}_2\}(Y)$ are Hodge-isometric. 
Then the periods 
$P(X,\phi)$ and $P(Y,\psi)$ belong to the same $\W(\Lambda)$ orbit,
by Lemma \ref{lemma-equivalent-data}. 
This contradicts the fact, that 
the periods $P(X,\phi)$ and $f(P(X,\phi))$ do not belong to the same
$\W(\Lambda)$ orbit in $\Omega$, by the assumptions on $X$ and $f$. 
\EndProof

\begin{example}
\label{example-genus-2-K3}
{\rm
Let $S$ be a $K3$ surface with a line bundle $L$, whose degree 
$c_1(L)^2$ is $2$ or $4$.  We can choose $S$, for example, 
to be a double cover of $\PP^2$, branched along a sextic, with $L$ 
the pullback of $\StructureSheaf{\PP^2}(1)$. $S$ could also be
a smooth quartic in $\PP^3$. 
We describe in equation (\ref{eq-reflection-by-minus-4-vector}) 
an orientation 
preserving Hodge isometry $f$ of
$H^2(S^{[7]},\Integers)$, which does not belong to 
$\W H^2(S^{[7]},\Integers)$. 

The Mukai vector in $H^*(S,\Integers)$, of the
ideal sheaf of a length $7$ subscheme, is $v:=(1,0,-6)$ 
(see equation (\ref{eq-Mukai-vector})). 
The direct summand $H^2(S,\Integers)$ of $H^*(S,\Integers)$ is contained in
the sublattice $v^\perp$ orthogonal to $v$ in the Mukai lattice. 
We identify $H^2(S,\Integers)$ 
with its image in $H^2(S^{[7]},\Integers)$ 
under the isometry $H^2(S^{[7]},\Integers)\cong v^\perp$ 
(Theorem \ref{thm-irreducibility}). 
The class $\delta:=-(1,0,6)$ in $v^\perp$ corresponds to half the class in 
$H^2(S^{[7]},\Integers)$ of the big diagonal. We see that
$v^\perp$, and hence $H^2(S^{[7]},\Integers)$, 
admits an orthogonal decomposition 
$H^2(S,\Integers)\oplus {\rm span}\{\delta\}$, as sublattices of
$H^*(S,\Integers)$.


Let $w_0$ be the primitive element $2c_1(L)+\delta$ in 
$H^2(S^{[7]},\Integers)$. 
Then $(w_0,w_0)=-4$ if $\deg(L)=2$ and $(w_0,w_0)=4$ if 
$\deg(L)=4$. Furthermore, 
$(w_0,x)$ is even, for every $x\in H^2(S^{[7]},\Integers)$. 
Define $f$ to be the reflection of $H^2(S^{[7]},\Integers)$ 
given by
\begin{equation}
\label{eq-reflection-by-minus-4-vector}
f(x) \ \ \ := \ \ \ \frac{-4}{(w_0,w_0)}\cdot x + 
\frac{(x,w_0)}{2}\cdot w_0.
\end{equation}
Then $f$ is an orientation preserving Hodge isometry of
$H^2(S^{[7]},\Integers)$. The quotient group
$H^2(S^{[7]},\Integers)^*/H^2(S^{[7]},\Integers)$ is generated by
the coset of $\delta/12$. When $\deg(L)=2$, then 
the image of $f$, in the residual isometry group
$O[H^2(S^{[7]},\Integers)^*/H^2(S^{[7]},\Integers)]$, is multiplication by 
$-5$, because $f(\delta)=-(12c_1(L)+5\delta)$. 
When $\deg(L)=4$, then the image of $f$ is multiplication by $-7$. 
Since $H^2(S^{[7]},\Integers)^*/H^2(S^{[7]},\Integers)$
is isomorphic to $\Integers/12\Integers$, then $f$ does not belong to 
$\W H^2(S^{[7]},\Integers)$, and 
$f$ does not extend to an isometry of the Mukai lattice 
(Lemma \ref{lemma-on-residual-orthogonal-group} part
\ref{lemma-item-three-characterizations-of-W}). 

The action of $f$, on the period domain 
$\Omega$, fixes the period of $S^{[7]}$. The fixed locus is
the hyperplane section $\Omega\cap w_0^\perp$. 
The isometry $f$ acts as an involution of the base of the local 
Kuranishi family of deformations of $S^{[7]}$. 
For a generic point $b$ in the base, the fibers over $b$ and $f(b)$
are non-bimeromorphic irreducible holomorphic symplectic manifolds, 
with isomorphic weight $2$ Hodge structures 
(Lemma \ref{lemma-non-bimeromorphic-classes}). 
}
\end{example}

\subsection{Hodge-isometric yet non-birational moduli spaces}
\label{sec-non-birational-moduli-spaces}
The {\em transcendental lattice} $\Theta(X)$, 
of an irreducible holomorphic symplectic manifold $X$, is
the sublattice of $H^2(X,\Integers)$ orthogonal to
$H^{1,1}(X,\Integers)$, with respect to the Beauville-Bogomolov form.

\begin{new-lemma}
\label{lemma-realization-as-a-moduli-space}
Let $S$ be a projective $K3$ surface and $X$ an
irreducible holomorphic symplectic 
manifold deformation equivalent to $S^{[n]}$, $n\geq 4$. 
There exists a primitive and effective class $v\in K^{1,1}(S)$ and
a $v$-generic ample line bundle $H$ on $S$, such that 
$\{Q^4,\bar{c}_2\}(X)$ is Hodge-isometric to 
$\{Q^4,\bar{c}_2\}(\M_H(v))$, if and only if 
$\Theta(X)$ and
$\Theta(S)$ are Hodge-isometric.
\end{new-lemma}

\preprint{
\noindent
{\bf Proof:}
If $\{Q^4,\bar{c}_2\}(X)$ is Hodge isometric to 
$\{Q^4,\bar{c}_2\}(\M_H(v))$, then both are Hodge isometric to 
$\{K(S),2v\}$, by Theorem 
\ref{thm-introduction-arithmetic-invariant-of-moduli}.
$\Theta(S)$ is Hodge isometric to $K^{1,1}(S)^\perp$.
Similarly, $\Theta(X)$ is Hodge isometric to
$Q^{2,2}(X,\Integers)^\perp$, by Theorem
\ref{thm-introduction-invariant-Q-4}
part \ref{thm-item-pair-of-embeddings}. 
Hence $\Theta(X)$ is Hodge isometric to $\Theta(S)$. 

Conversely, 
let $g:\Theta(S)\rightarrow \Theta(X)$ be a Hodge-isometry.
Suffices to show, that there exists a primitive and effective class 
$v\in K^{1,1}(S)$, such that $\{Q^4,\bar{c}_2\}(X)$ is 
Hodge-isometric to $\{K(S),2v\}$, by Theorem 
\ref{thm-introduction-arithmetic-invariant-of-moduli}.
Denote by
$\alpha:\Theta(S)\hookrightarrow K(S)$ and
$\beta: \Theta(X)\hookrightarrow H^2(X,\Integers)$ the natural inclusions
and by $e_X:H^2(X,\Integers)\hookrightarrow Q^4(X,\Integers)$ 
one of the two primitive embeddings in Theorem
\ref{thm-introduction-invariant-Q-4} part \ref{thm-item-pair-of-embeddings}.
Both $\alpha$ and $e_X\circ \beta\circ g$
are primitive embeddings of $\Theta(S)$.
There exists an isometry
$
f'  :  K(S) \rightarrow Q^4(X,\Integers)
$
satisfying $f'\circ \alpha = e_X\circ \beta \circ g$, by Theorem 1.14.4
in \cite{nikulin}. Any such $f'$ is necessarily a Hodge-isometry. 
Set $w:=(f')^{-1}(\bar{c}_2(X))/2$. If $w$ is effective, 
set $v:=w$ and $f=f'$.
Otherwise, set $v=-w$ and $f=-f'$. Then 
$f$ is a Hodge-isometry 
between $\{K(S),2v\}$ and $\{Q^4,\bar{c}_2\}(X)$.
\EndProof
}

\journal{We will not use the Lemma, but it motivates the Proposition below.
The easy proof can be found in the preprint version of this paper
\cite{markman-eprint-version}.}
Let $S$ be a $K3$ surface, $H$ an ample line-bundle on $S$, 
and $n\geq 4$ an integer. Assume that $S^{[n]}$ satisfies 
Condition \ref{cond-stabilizer-of-the-period-is-contained-in-W}.
Let $\rho(n-1)$ be the number of primes dividing $n-1$, as in Lemma
\ref{lemma-index-of-smaller-symmetry-group}. 
Condition \ref{cond-stabilizer-of-the-period-is-contained-in-W} 
implies, that there are at least $2^{\rho(n-1)-1}$
pairwise non-birational irreducible holomorphic symplectic manifolds,
deformation equivalent to $S^{[n]}$, with the same weight $2$ 
Hodge structure as $S^{[n]}$
(Theorem \ref{thm-non-bimeromorphic-classes}). 
Proposition \ref{thm-non-birational-moduli-spaces}
identifies such a set consisting of $2^{\rho(n-1)-1}$ moduli spaces of 
sheaves on $S$ (as expected by 
Lemma \ref{lemma-realization-as-a-moduli-space}). 

Let $\P_n$ be the subset of $\Integers\oplus\Integers$ given 
in equation (\ref{eq-P-n}).
Let $\M_H(r,0,s)$ be the moduli space of $H$-stable sheaves with Mukai vector 
$(r,0,s)$. We use the notation, for vectors in the Mukai lattice,
introduced in equation (\ref{eq-Mukai-vector}). 

\begin{prop} 
\label{thm-non-birational-moduli-spaces}
\begin{enumerate}
\item
\label{cor-item-moduli-are-projective}
The moduli spaces $\M_H(r,0,s)$, $(r,s)\in \P_n$, are non-empty, smooth, 
and projective. Their second cohomology groups 
$H^2(\M_H(r,0,s),\Integers)$ are all Hodge-isometric to 
$H^2(S^{[n]},\Integers)$. 
\item
\label{cor-item-non-birational}
The moduli spaces $\M_H(r,0,s)$, for $(r,s)\in \P_n$, 
are pairwise non-birational, provided $H^2(S^{[n]},\Integers)$ 
satisfies Condition 
\ref{cond-stabilizer-of-the-period-is-contained-in-W}.
\item
\label{cor-item-moduli-spaces-represent-all-invariants}
Let $X$ be an irreducible holomorphic symplectic manifold, 
deformation equivalent to $S^{[n]}$,
such that $H^2(X,\Integers)$ is Hodge isometric to $H^2(S^{[n]},\Integers)$.
Then the pair $\{Q^4,\bar{c}_2\}(X)$ is Hodge-isometric to 
$\{Q^4,\bar{c}_2\}(\M_H(r,0,s))$, for some $(r,s)$ in $\P_n$.
\end{enumerate}
\end{prop}

Condition \ref{cond-stabilizer-of-the-period-is-contained-in-W}
is satisfied by the Hilbert scheme $S^{[n]}$ of 
any K3 surface $S$, with a cyclic Picard group
generated by a line bundle $H$ of degree 
$m^2(2n-2)$, for some integer $m$ 
(see Lemma \ref{lemma-S-H-n-satisfying-condition} for a stronger criterion). 
The generic K\"{a}hler $K3$ surface satisfies the condition, 
but at present the moduli spaces 
considered in Proposition \ref{thm-non-birational-moduli-spaces} 
are constructed only for projective K3 surfaces. 

\noindent
{\bf Proof of Proposition \ref{thm-non-birational-moduli-spaces}:}
\ref{cor-item-moduli-are-projective})
We prove first that $\M_H(r,0,s)$ is projective, 
regardless of the choice of $H$, by 
showing that every $H$-semistable sheaf 
with Mukai vector $(r,0,s)$ is $H$-stable.
Let $F$ be an $H$-semistable sheaf with Mukai vector $(r,0,s)$,
$F'$ a subsheaf, and assume that their normalized Hilbert polynomials, given 
in equation (\ref{eq-inequality-of-normalized-Hilbert-polynomials}), 
are equal.
Then $c_1(F')\cup c_1(H)$ must vanish. 
Hence, the coefficients of the Hilbert polynomial of $F'$ are determined by
the rank and Euler characteristic of $F'$. Now the rank $r$ 
and Euler characteristic $r+s$ of $F$ are relatively prime. 
Hence, the equality of the normalized Hilbert polynomials 
implies equality of the ranks and Euler characteristics. 
It follows, that $F'=F$.

$\M_H(r,0,s)$ is non-empty and 
$H^2(\M_H(r,0,s),\Integers)$ is Hodge isometric to 
the orthogonal complement $(r,0,s)^\perp$, in the Mukai lattice,
by Theorem \ref{thm-irreducibility}.
Hence, the lattice embedding $\iota_{r,s}$, of Lemma 
\ref{lemma-the-index-of-W-versus-orbits-in-hyperbolic-plane}, 
induces a Hodge isometry between $\Lambda:=H^2(S^{[n]},\Integers)$
and $H^2(\M_H(r,0,s),\Integers)$.

\ref{cor-item-non-birational}) 
Assume that $\M_H(r_1,0,s_1)$  and 
$\M_H(r_2,0,s_2)$ are birational and 
$(r_1,s_1)$ and $(r_2,s_2)$ are different elements of $\P_n$. 
We will show that $S^{[n]}$ does not satisfy Condition
\ref{cond-stabilizer-of-the-period-is-contained-in-W}.
The invariants
$\{Q^4,\bar{c}_2\}(\M_H(r_i,0,s_i))$,
$i=1,2$, are Hodge isometric, since the moduli
spaces are birational. The identification of these
invariants, in Theorem \ref{thm-introduction-arithmetic-invariant-of-moduli}, 
implies that there is a Hodge isometry $g$ 
of the Mukai lattice $\widetilde{\Lambda}:=H^*(S,\Integers)$, mapping 
$(r_1,0,s_1)$ to $(r_2,0,s_2)$. 
Let $f\in O(\Lambda)$ be the unique isometry satisfying
$
g\circ \iota_{r_1,s_1}\circ f^{-1}  =  \iota_{r_2,s_2}.
$
Then $f$ does not belong to the stabilizer of the orbit 
$O(\widetilde{\Lambda})\cdot g\circ \iota_{r_1,s_1}$ in 
$O(\Lambda,\widetilde{\Lambda})$, 
because the embeddings $g\circ \iota_{r_1,s_1}$ and $\iota_{r_2,s_2}$ 
belong to different $O(\widetilde{\Lambda})$ orbits, by part
\ref{lemma-item-P-n-represent-O-Lambda-tilde-orbits} of Lemma
\ref{lemma-the-index-of-W-versus-orbits-in-hyperbolic-plane}.
Hence, neither $f$, nor $-f$, belongs to $\W(\Lambda)$, by part 
\ref{lemma-item-W-is-the-stabilizer-of-action-on-orbit-space} of Lemma
\ref{lemma-the-index-of-W-versus-orbits-in-hyperbolic-plane}.

On the other hand, $f$ is clearly a Hodge isometry of $\Lambda$,
when $\Lambda$ is endowed with the Hodge structure of
$H^2(S^{[n]},\Integers)$. Hence, $S^{[n]}$ does not satisfy Condition
\ref{cond-stabilizer-of-the-period-is-contained-in-W}.

\journal{
\ref{cor-item-moduli-spaces-represent-all-invariants})
The routine proof of part  
\ref{cor-item-moduli-spaces-represent-all-invariants}  
can be found in the preprint version
\cite{markman-eprint-version}. 
}
\preprint{
\ref{cor-item-moduli-spaces-represent-all-invariants})
Choose an isometry $g:Q^4(X,\Integers)\rightarrow \widetilde{\Lambda}$. 
The orthogonal complement $\bar{c}_2(X)^\perp$ in
$Q^4(X,\Integers)$ is Hodge isometric to $H^2(X,\Integers)$,
by part \ref{thm-item-pair-of-embeddings} of Theorem
\ref{thm-introduction-invariant-Q-4}. By assumption, we can choose a Hodge 
isometry 
$f:H^2(S^{[n]},\Integers)\rightarrow \bar{c}_2(X)^\perp$. 
Recall that $\Lambda$ is identical to $H^2(S^{[n]},\Integers)$. 
Thus, $g\circ f$ is an isometric embedding of $\Lambda$ in 
$\widetilde{\Lambda}$. Hence, there exists an isometry 
$h\in O(\widetilde{\Lambda})$, satisfying the equality 
\[
h\circ g\circ f \ \ \ = \ \ \ \iota_{r,s},
\]
for some $(r,s)$ in $\P_n$, by part
\ref{lemma-item-P-n-represent-O-Lambda-tilde-orbits} of Lemma
\ref{lemma-the-index-of-W-versus-orbits-in-hyperbolic-plane}.
Now, $h\circ g: Q^4(X,\Integers)\rightarrow \widetilde{\Lambda}$
is necessarily a Hodge isometry, if $\iota_{r,s}$ is, and that is the case, 
when $\widetilde{\Lambda}$
is endowed with the Hodge structure of the Mukai lattice of $S$.
Furthermore, $h\circ g$ maps $\bar{c}_2(X)/2$ to 
$(r,0,s)$ or $-(r,0,s)$.
Theorem \ref{thm-introduction-arithmetic-invariant-of-moduli} implies, that 
$\{Q^4,\bar{c}_2\}(X)$ is Hodge-isometric
to $\{Q^4,\bar{c}_2\}(\M_H(r,0,s))$
via $\varphi^4\circ h\circ g$ or $-\varphi^4\circ h\circ g$
}
\EndProof

\journal{
\medskip
The routine proof of the following Lemma 
can be found in the preprint version
\cite{markman-eprint-version}.
}

\begin{new-lemma}
\label{lemma-S-H-n-satisfying-condition}
Assume that the Picard group of $S$ is
generated by a line bundle $H$ with $c_1(H)^2=m(2n-2)$, for some 
non-negative integer $m$. Assume further, either that $m$ is the square of
an integer, or that 
every unit in the ring 
$\Integers[\sqrt{m}]$ has trace $2$ or $-2$ (mod $4n-4$).
Then $S^{[n]}$ satisfies Condition
\ref{cond-stabilizer-of-the-period-is-contained-in-W}.
\end{new-lemma}

\preprint{
\noindent
{\bf Proof:} The notation $x=\pm y$ would mean $x=y$ or $x=-y$.
$H^2(S^{[n]},\Integers)$ is isomorphic to the sublattice 
$(1,0,1-n)^\perp$ of the Mukai lattice. The orthogonal complement
$(1,0,1-n)^\perp$ decomposes as $H^2(S,\Integers)\oplus \Integers\delta$,
where $\delta=(1,0,n-1)$. 
A Hodge isometry $g$ of $(1,0,1-n)^\perp$
induces a residual isometry $\pi(g)$ 
on $((1,0,1-n)^\perp)^*/(1,0,1-n)^\perp$. 
If $g$ maps $\delta$ to the $(1,1)$ class $a c_1(H)+b\delta$, 
then $\pi(g)$ acts by multiplication by $b$ (mod $2n-2$). 
Hence, $\pi(g)=\pm1$ if and only if $b\equiv\pm 1$ (mod $2n-2$). 
If $g$ is orientation preserving and $\pi(g)=\pm 1$, then $g$
belongs to  the subgroup $\W$,
by Lemma \ref{lemma-on-residual-orthogonal-group} part 
\ref{lemma-item-three-characterizations-of-W}. 

The image $g(\delta)$ has
square-length $2-2n$ and hence the linear combination
$a c_1(H)+b\delta$ satisfies the equation 
$2-2n=(2-2n)(b-\sqrt{m}a)(b+\sqrt{m}a)$. 
If $\sqrt{m}$ is a positive integer, then 
the only integral solution is $a=0$ and $b=1$ or $b=-1$.
If $m=0$, then $b^2=1$ as well. 
If $m$ is not a square, then $2b={\rm tr}(b+\sqrt{m}a)=\pm 2$ 
(mod $4n-4$), by assumption. Hence, $b\equiv \pm 1$ (mod $2n-2$). 
Thus, $g$ belongs to $\W$.
\EndProof
}

%
\section{Monodromy constraints via non-split extensions}
\label{sec-extension-classes}
In section
\ref{sec-reduction-of-monodromy-conjecture-to-generators-conjecture}
we prove Corollary 
\ref{cor-reduction-of-monodromy-conjecture-to-generators-conjecture}
about the injectivity of the homomorphism $Mon(X)\rightarrow Mon^2$. 
The rest of the section is dedicated to the proof of
Theorem 
\ref{thm-monodromy-operator-is-determined-by-its-weight-2-action}.
An outline of the proof is given in section 
\ref{subsec-Monodromy-constaints-via-non-split-extensions}.

%
\subsection{A monodromy operator is determined by its action on $H^2$}
\label{sec-reduction-of-monodromy-conjecture-to-generators-conjecture} 
{\bf Proof of Corollary
\ref{cor-reduction-of-monodromy-conjecture-to-generators-conjecture}:}
The statement is clear if $n=2$, 
as $H^*(S^{[2]},\RationalNumbers)$ is generated by 
$H^2(S^{[2]},\RationalNumbers)$, by Lemma \ref{lemma-10-in-markman-diagonal}. 
The statement follows from Theorem 
\ref{thm-monodromy-operator-is-determined-by-its-weight-2-action}
if $n=3$, because $H^*(S^{[3]},\RationalNumbers)$ is generated by 
$H^i(S^{[3]},\RationalNumbers)$, $i\leq 4$, 
by Lemma \ref{lemma-10-in-markman-diagonal}. 
Theorem \ref{thm-vanishing-of-varphi-k} states that 
$H^*(S^{[n]},\RationalNumbers)$ is generated by 
$H^d(S^{[n]},\RationalNumbers)$, 
for $d\leq 4+2\lfloor\frac{n}{2}\rfloor$.

If $n$ is even, then the kernel $K^2$ of
the restriction homomorphism 
$Mon(S^{[n]})\rightarrow Mon^2(S^{[n]})$ acts faithfully on 
$H^{n+4}(S^{[n]},\Integers)$, by Theorems 
\ref{thm-monodromy-operator-is-determined-by-its-weight-2-action}
and \ref{thm-vanishing-of-varphi-k}.
$K^2$ acts faithfully also on the representation 
$Q^{n+4}(S^{[n]},\RationalNumbers)$, 
since the $Mon(S^{[n]})$-action
is semi-simple (Proposition \ref{prop-so-action-determines-Hodge-structure}
part \ref{prop-item-semi-simple}). 
Theorem \ref{thm-vanishing-of-varphi-k} implies that 
$\rank(Q^{n+4}(S^{[n]}))\leq 1$. 
Hence, the order of $K^2$ is $\leq 2$. 
If $n$ is odd, the order of $K^2$ is $\leq 2$, by the same argument applied 
to $Q^{n+3}(S^{[n]},\RationalNumbers)$.

Let $v\in K(S)$ be the class of the ideal sheaf of a length $n$ subscheme
of $S$.
If $n\equiv 0$ modulo $4$, then $2\varphi^{n+4}(v)$ is equal to the image
$\bar{c}_{(n/2)+2}(S^{[n]})$ of $c_{(n/2)+2}(S^{[n]})$ in
$Q^{n+4}(S^{[n]},\Integers)$,
by part \ref{thm-item-identification-of-residue-of-chern-class} of Theorem
\ref{thm-summary-of-monodromy-results}. 
The homomorphism
$\varphi^{n+4}:K(S)\rightarrow Q^{n+4}(S^{[n]},\Integers)$
is $O^+K(S)_v$-equivariant. Hence, $\varphi^{n+4}(v^\perp)=(0)$
and $Q^{n+4}(S^{[n]},\RationalNumbers)$ is spanned by $\varphi^{n+4}(v)$.
We conclude that $Mon(S^{[n]})$ acts
trivially on $Q^{n+4}(S^{[n]},\RationalNumbers)$ and 
the restriction homomorphism 
$Mon(S^{[n]})\rightarrow Mon^2(S^{[n]})$ is injective.
If $n\equiv 1$ modulo $4$, we prove the injectivity by the same argument
applied to $Q^{n+3}(S^{[n]},\RationalNumbers)$.
\EndProof

{\bf Geometric implications of Conjecture 
\ref{conj-Mon-isomorphic-to-Mon-2}:}
The following is a fundamental theorem of Namikawa:

\begin{thm}
\cite{namikawa-deformations}
Let $\pi:X\rightarrow Y$ be a symplectic resolution of a 
projective symplectic variety\footnote{A projective symplectic variety
is a normal projective variety with rational Gorenstein singularities, 
which regular locus admits a non-degenerate holomorphic two-form.}
$Y$ of dimension $n$. Then the Kuranishi spaces
$Def(X)$ and $Def(Y)$ are both smooth of the same dimension.
There exists a natural map
$\pi_*:Def(X)\rightarrow Def(Y)$ and $\pi_*$ is a finite 
branched covering.
Moreover, $Y$ has a flat deformation to a smooth symplectic $n$-fold $Y_t$.
Any smoothing $Y_t$ of $Y$ is a symplectic $n$-fold obtained as a flat 
deformation of $X$. 
\end{thm}

Let $Def(Y)^0$ be the complement of the branch divisor and 
$\widetilde{Def}(X)^0\rightarrow  Def(Y)^0$ the Galois closure 
of $\pi_*$. Its Galois group is easily seen to act on
$H^*(X,\Integers)$ via monodromy operators, 
preserving the Hodge-structure 
(see \cite{markman-monodromy-I} section 2). 
These are examples of {\em local} monodromy operators. 
It is often not difficult to calculate their action 
on $H^2(X)$, but that on $H^*(X)$ is harder to calculate. 

Let $\tilde{g}\in Mon(X)$ be such a local monodromy operator. 
There is one more reason to calculate $\tilde{g}$.
The action of
$\tilde{g}$ on $H^*(X)$ is induced by a Lagrangian correspondence 
$Z_{\tilde{g}}\subset X\times X$ 
(see the proof of \cite{huybrects-basic-results}, Theorem 4.3). 
This realizes $\tilde{g}$ also as an operator in a Lagrangian 
{\em convolution algebra}. See \cite{chriss-ginzburg} 
for the definition of this algebra. 
A conjectural formula for $\tilde{g}$  is given by Nakajima 
in terms of the generators of a convolution algebra
in a particular example of moduli spaces of sheaves on a $K3$ surface
(see the last paragraph is section 3 of \cite{nakajima-representations}).

Let $g\in \W(X)=Mon^2(X)$ be  the restriction to $H^2(X)$
of a local monodromy operator $\tilde{g}\in Mon(X)$. 
Another lift of $g$ to an operator in $Mon(X)$ is provided by the 
homomorphism $\nu:\W(X)\rightarrow Mon(X)$, given in 
Lemma \ref{lemma-mu-is-an-isomorphism}.
Conjecture \ref{conj-Mon-isomorphic-to-Mon-2} implies the equality 
$\tilde{g}=\nu(g)$. This equality 
combines with 
Theorem 
\ref{thm-summary-of-monodromy-results} part 
\ref{thm-item-mon-g-sends-a-universal-classes-to-such} to 
explicitly describe how the local monodromy operator $\tilde{g}$ acts on a set
of generators for the cohomology ring of $X$, when $X$ is a moduli
space of sheaves. Let us illustrate this in an example.
\begin{example}
{\rm
Set $X=S^{[n]}$, $n\geq 2$, and let $g$ be the reflection 
$g(x)=x+\frac{(x,\delta)}{n-1}\delta$,
where $\delta$ is half the class of the big diagonal. 
Then $g$ indeed lifts to a local monodromy operator $\tilde{g}$, induced
by the  fiber product $Z_{\tilde{g}}:=[S^{[n]}\times_{S^{(n)}} S^{[n]}]$,
with respect to the Hilbert-Chow morphism $S^{[n]}\rightarrow S^{(n)}$
onto the symmetric product (\cite{markman-monodromy-I}, Lemma 2.7).
$S^{[n]}\times_{S^{(n)}} S^{[n]}$ is reducible and 
its irreducibe components are in bijection with ordered partitions of $n$.
Choose a basis $\{y_1, y_2, \dots, y_{24}\}$ of $K(S)$ satisfying
$y_j^\vee=y_j$, for $j=1, 2$, and
$y_j^\vee=-y_j$, for $3\leq j\leq 24$. 
Let $\E$ be the ideal sheaf of the universal subscheme in $S\times S^{[n]}$, 
and set $u_{y_j}:=u(y_j)$, where $u:K(S)\rightarrow K(S^{[n]})$ is
given in (\ref{eq-u}). 
Conjecture \ref{conj-Mon-isomorphic-to-Mon-2} and Theorem 
\ref{thm-summary-of-monodromy-results} part 
\ref{thm-item-mon-g-sends-a-universal-classes-to-such}
yield the equlity: 
$\tilde{g}(ch(u_{y_j}))=ch(u_{(y_j^\vee)}^\vee)$,
so \ \ \ 
\[
\tilde{g}(ch_i(u_{y_j})) = 
\left\{
\begin{array}{rcc}
(-1)^i ch_i(u_{y_j}) & \mbox{if} & j=1, 2,
\\
(-1)^{i+1}ch_i(u_{y_j}) & \mbox{if} & 3\leq j\leq 24.
\end{array}\right.
\]
The latter equality is the desired description of the action of the 
local monodromy operator 
$\tilde{g}$ on the generators $ch_i(u_{y_j})$ of the cohomology ring
(Theorem \ref{thm-integral-generators}).
}
\end{example}

%
\subsection{Extensions of integral representations}
\label{subsec-Monodromy-constaints-via-non-split-extensions}
Theorem 
\ref{thm-monodromy-operator-is-determined-by-its-weight-2-action}
is shown to easily follow from
Proposition \ref{prop-extension-of-Q-bar-by-Z-hat-does-not-split},
the statement and proof of which occupies the rest of this section.
The following notation will be used throughout the rest of section
\ref{sec-extension-classes}.
Let $v$ be an effective and primitive Mukai vector and $H$ a 
$v$-generic polarization. Set $\M:=\M_H(v)$ and 
$n:=\dim_\ComplexNumbers(\M)/2$. Then $(v,v)=2n-2$.

Assume that $n\geq 2$. 
Let $G$ be the group $O^+K(S)_v$ of orientation 
preserving isometries of the Mukai lattice stabilizing $v$.
We identify $G$ also with the subgroup of $\W H^2(\M,\Integers)$, which is
the image of $G$ via the homomorphism 
$\monrep^2$, introduced in Theorem \ref{thm-summary-of-monodromy-results}.
The homomorphism $\monrep^2$ 
is injective for $n\geq 3$ and its 
restriction to $G$ is injective for $n=2$ as well, by Lemma 
\ref{lemma-mu-is-an-isomorphism}. 

Let $(A_{2i-4})^{2i}$, $i\geq 2$, be the subgroup of $H^{2i}(\M,\Integers)$
used in the definition of $Q^{2i}(\M,\Integers)$ in 
equation (\ref{eq-C-d}). Set
\[
E^{2i} \ \ \  := \ \ \  H^{2i}(\M,\Integers)/ (A_{2i-4})^{2i}.
\]
Then $Q^{2i}(\M,\Integers)$ is a quotient of $E^{2i}$.
Let $\overline{Q}^{2i}$ be the quotient of $Q^{2i}(\M,\Integers)$
by its torsion subgroup $Q^{2i}(\M,\Integers)_{tor}$. 
Let $\widehat{Z}^{2i}$ be the kernel of the natural 
$G$-equivariant homomorphism below 
\begin{equation}
\label{eq-extension-of-Q-bar-by-Z-hat}
0\rightarrow \ \widehat{Z}^{2i} \ \longrightarrow \
E^{2i} \  \longrightarrow \ \overline{Q}^{2i} \ \rightarrow 0.
\end{equation} 

When a short exact sequence of $G$-modules is displayed in a 
numbered equation, 
such as (\ref{eq-extension-of-Q-bar-by-Z-hat}),
we denote by $\epsilon_{(\ref{eq-extension-of-Q-bar-by-Z-hat})}$
the corresponding $1$-extension class. 
See section \ref{sec-elementary-calculations-of-extension-classes} 
for the definition of the extension groups of $G$-modules.
The following Proposition is the main technical result of this section. 

\begin{prop}
\label{prop-extension-of-Q-bar-by-Z-hat-does-not-split}
Set $n:=\dim_{\ComplexNumbers}(\M)/2$ and assume that $n\geq 3$ and  
$\overline{Q}^{2i}$ does not vanish. 
\begin{enumerate}
\item
\label{prop-item-summary-order-of-extension-is-geq-3}
The order of  $\epsilon_{(\ref{eq-extension-of-Q-bar-by-Z-hat})}$
in $\Ext^1_G(\overline{Q}^{2i},\widehat{Z}^{2i})$
is $\geq 3$, provided 
$2\leq i \leq \frac{n+2}{2}$.
\item
\label{prop-item-order-of-extension-class-of-Q-bar-by-Z-hat}
When $i=2$, the order 
of $\epsilon_{(\ref{eq-extension-of-Q-bar-by-Z-hat})}$ is divisible by
$2n-2$, if $n$ is odd, and 
by $n-1$, if $n$ is even. 
The order of $\epsilon_{(\ref{eq-extension-of-Q-bar-by-Z-hat})}$
is divisible by $\frac{2n-2}{\gcd(i-1,2n-2)}$, if $i$ is an integer in the 
range $3\leq i \leq \frac{n+2}{2}$. 
\end{enumerate}
\end{prop}

Proposition 
\ref{prop-extension-of-Q-bar-by-Z-hat-does-not-split}
is proven in section \ref{sec-reduction-of-prop-extension-to-prop}. 

{\bf Proof of Theorem 
\ref{thm-monodromy-operator-is-determined-by-its-weight-2-action}:}
It suffices to prove the Theorem for a moduli space $\M$ as above. 
The case $n=2$ is clear, as $H^*(\M,\RationalNumbers)$ is generated by 
$H^2(\M,\RationalNumbers)$, by Lemma \ref{lemma-10-in-markman-diagonal}. 
Assume $n\geq 3$. 
Let $f\in Mon(\M)$ be a monodromy operator acting trivially on 
$H^2(\M,\Integers)$. We prove by induction on $1\leq i\leq \frac{n+2}{2}$,
that $f$ acts trivially on $H^{2i}(\M,\Integers)$.
The case $i=1$ is clear. 
Assume that $i\geq 2$ and 
$f$ acts trivially on $A_{2i-2}$.
It suffices to prove that $f$ acts trivially on 
$\overline{Q}^{2i}(\M,\Integers)$, since 
the cohomology ring $H^*(\M,\RationalNumbers)$ is a semi-simple 
representation of $Mon(\M)$, by part
\ref{prop-item-semi-simple} of Proposition
\ref{prop-so-action-determines-Hodge-structure}. 
$f$ commutes with 
$Mon(\M)$, by part \ref{thm-item-markman-monodromy-I-lemma-4.4}
of Theorem \ref{thm-summary-of-monodromy-results}. 
We need to treat two cases:

Case 1: Assume that $Q^{2i}(\M,\Integers)$ is torsion free. 
If $\varphi^{2i}(\Integers v)$ or $\varphi^{2i}(v^\perp)$ does not vanish, 
then it is an irreducible $G$-submodule, and 
$f$ must act on it via multiplication by a scalar equal to $1$ or $-1$.
If both do not vanish, then 
these two scalars are equal, by Lemma
\ref{lemma-orders-of-extension-classes-with-Mukai-lattice-in-the-middle}
and the $G$-equivariance of $\varphi^{2i}$
(otherwise, the order of 
$\epsilon_{(\ref{eq-extension-of-v-perp-dual-by-span-v})}$
would be $1$ or $2$).
If $f$ acts on $Q^{2i}(\M,\Integers)$ via multiplication by $-1$, then
$2\epsilon_{(\ref{eq-extension-of-Q-bar-by-Z-hat})}=0$, contradicting the fact
that the order of $\epsilon_{(\ref{eq-extension-of-Q-bar-by-Z-hat})}$
is larger than $2$, by
Proposition \ref{prop-extension-of-Q-bar-by-Z-hat-does-not-split}.

Case 2: Assume that both $Q^{2i}(\M,\Integers)_{tor}$ and 
$\overline{Q}^{2i}$ do not vanish. Then $\overline{Q}^{2i}$ is
an irreducible $Mon(\M)$-module, so $f$ acts on it via multiplication by $1$ 
or $-1$. The latter case is excluded, since the order of 
$\epsilon_{(\ref{eq-extension-of-Q-bar-by-Z-hat})}$ is larger than $2$, by
Proposition \ref{prop-extension-of-Q-bar-by-Z-hat-does-not-split}. 
\EndProof


\medskip
The rest of section \ref{sec-extension-classes} is dedicated to the proof
of Proposition \ref{prop-extension-of-Q-bar-by-Z-hat-does-not-split}.
The discussion is complicated by the need to consider special cases. 
The reader is advised to concentrate on the generic case
$3\leq i \leq \frac{n}{2}$ in the first reading.
In the generic case $Q^{2i}$ and $\widehat{Z}^{2i}$ are torsion free  
and many of the arguments simplify. 
The proof of the generic case is concluded in section
\ref{sec-proof-of-generic-case}. Example 
\ref{example-warm-up-for-non-generic-cases} 
illustrates the main idea of the proof of Proposition 
\ref{prop-extension-of-Q-bar-by-Z-hat-does-not-split} 
both in the case $i=2$ and in case $i=\frac{n+1}{2}$ and $n$ odd. 
The cases in which $i>\frac{n}{2}$ are more subtle, since 
the homomorphism $\varphi^{2i}$, given in (\ref{eq-varphi-k}), is
no longer an isomorphism. 
The case $i=2$ is special, since $ch_2$
does not fit the general pattern for $ch_i$, $i\geq 3$.
Write $(-1)^{i-1}i!ch_i$ as a polynomial in the Chern classes $c_i$.
For $i\geq 3$, the coefficients of the 
first two terms $(-1)^{i-1}i!ch_i=ic_i-ic_{i-1}(y)c_1(y)+\cdots$
are the negative of each other, while $-2ch_2=2c_2-c_1^2$
(this first comes up in Lemma \ref{lemma-sigma-is-linear}).

%
\subsection{Computing extension classes via rational splittings}
\label{sec-a-technique}
We formulate in this section the elementary 
Lemmas \ref{lemma-order-of-extension-class} and \ref{lemma-width} 
for calculating the order and width of an extension class of $G$-modules.
Let

\begin{equation}
\label{eq-diagram-of-three-extensions}
\begin{array}{ccccc}
Z & \LongRightArrowOf{\iota} & E & \LongRightArrowOf{j} & L
\\
z \ \downarrow \ \hspace{1ex} & & 
e \ \downarrow \ \hspace{1ex} & & 
= \ \downarrow \ \hspace{1ex}
\\
Z' & \LongRightArrowOf{\iota'} & E' & \LongRightArrowOf{j'} & L
\\
\zeta \ \downarrow \ \hspace{1ex} & & 
\eta \ \downarrow \ \hspace{1ex} & & 
\downarrow 
\\
Z'' & \LongRightArrowOf{\iota''} & E'' & \LongRightArrowOf{} & 0
\end{array}
\end{equation}
be a commutative diagram of $G$-modules, with short exact rows and columns. 
Note, that the second horizontal extension is necessarily 
the push-forward of the first 
via $z: Z \rightarrow Z'$, by Lemma 
\ref{lemma-characterization-of-pullback-pushforward}. 
Consider the long exact sequences of extension
groups associated to the rows and columns in
diagram (\ref{eq-diagram-of-three-extensions}) via the right-derived
functors of the left-exact covariant functor $\Hom_G(L,\bullet)$. 
The first row $(\ref{eq-diagram-of-three-extensions})_{top}$ 
of the diagram yields
\begin{equation}
\label{eq-long-exact-sequence-of-extensions-for-first-row}
0\rightarrow \Hom_G(L,Z) \RightArrowOf{\iota_*} \Hom_G(L,E) 
\RightArrowOf{j_*}  \Hom_G(L,L) 
\RightArrowOf{\delta}  \Ext^1_G(L,Z) \rightarrow \cdots 
\end{equation}
The left column yields
\begin{equation}
\label{eq-long-exact-sequence-of-extensions-for-left-column}
0\rightarrow \Hom_G(L,Z) \RightArrowOf{z_*} \Hom_G(L,Z') 
\RightArrowOf{\zeta_*}  \Hom_G(L,Z'') 
\RightArrowOf{f}  \Ext^1_G(L,Z) \rightarrow \cdots 
\end{equation}

We would like to describe the connecting homomorphism $\delta$ 
in terms of $f$, 
assuming that the middle horizontal sequence 
of diagram (\ref{eq-diagram-of-three-extensions})
admits a $G$-equivariant 
splitting. 
The second and third rows of the diagram yield the commutative diagram 
with short exact rows
\begin{equation}
\label{eq-snake-lemma-diagram}
\begin{array}{ccccccc}
0  \rightarrow &
\Hom_G(L,Z') & \RightArrowOf{\iota'_*} & \Hom_G(L,E') &
\RightArrowOf{j'_*} & \Hom_G(L,L) &
\rightarrow 0
\\
& \zeta_* \ \downarrow \ \hspace{2ex} & & 
\eta_* \ \downarrow \ \hspace{2ex} & &
\downarrow 
\\
0  \rightarrow &
\Hom_G(L,Z'') & \RightArrowOf{\iota''_*} & \Hom_G(L,E'') &
\rightarrow & 0
\end{array}
\end{equation}
The kernel of $\zeta_*$ is $\Hom_G(L,Z)$ and its co-kernel 
is a subgroup of $\Ext^1_G(L,Z)$, since the left column of diagram
(\ref{eq-diagram-of-three-extensions}) is short exact. 
The analogous statement for $\eta_*$ holds as well.

\begin{new-lemma}
\label{lemma-computation-of-extension-class}
The connecting homomorphism $\delta$, 
in equation (\ref{eq-long-exact-sequence-of-extensions-for-first-row}), 
is equal to the homomorphism
from $\Hom_G(L,L)$ to $\Ext^1_G(L,Z)$, obtained by the Snake Lemma
(\cite{hilton-stammbach}, Ch III Lemma 5.1),
from the kernel of the right vertical (zero) homomorphism in
diagram (\ref{eq-snake-lemma-diagram}) to the co-kernel of 
the left vertical homomorphism $\zeta_*$. 
\end{new-lemma}

\journal{A straightforward proof can be found in the
preprint version of this paper
\cite{markman-eprint-version}.
}

\preprint{
\noindent
{\bf Proof:}
One easily proves a slightly more general statement. Consider 
a commutative diagram
$$
\begin{array}{ccccc}
Z & \LongRightArrowOf{\iota} & E & \LongRightArrowOf{j} & L
\\
\tilde{z} \ \downarrow \ \hspace{1ex} & & 
\tilde{e} \ \downarrow \ \hspace{1ex} & & 
\tilde{\ell} \ \downarrow \ \hspace{1ex}
\\
Z_1 & \LongRightArrowOf{\iota_1} & E_1 & \LongRightArrowOf{j_1} & L_1
\\
\tilde{\zeta} \ \downarrow \ \hspace{1ex} & & 
\tilde{\eta} \ \downarrow \ \hspace{1ex} & & 
\downarrow 
\\
Z_2 & \LongRightArrowOf{\iota_2} & E_2 & \LongRightArrowOf{j_2} & L_2
\end{array}
\eqno({*})
$$
of $G$-modules with short exact rows and columns.
Assume that the middle (and hence also the bottom) row admits a 
$G$-equivariant splitting.
Clearly, diagram (\ref{eq-diagram-of-three-extensions}) is a special case.
Another special case is used in the definition of the connecting
homomorphism $\delta:\Hom_G(L,L)\rightarrow \Ext^1_G(L,Z)$, 
where one proves the existence of such a diagram ($*$)
with $Z_1$, $E_1$, and $L_1$ injective $G$-modules,
so that the columns are injective presentations
(\cite{hilton-stammbach}, Ch. III Lemma 5.4 and exercise 5.5). 
The analogue of diagram 
(\ref{eq-snake-lemma-diagram}) is 
$$
\begin{array}{ccccccc}
0  \rightarrow &
\Hom_G(L,Z_1) & \rightarrow
& \Hom_G(L,E_1) &
\rightarrow
& \Hom_G(L,L_1) &
\rightarrow 0
\\
& \tilde{\zeta}_* \ \downarrow \ \hspace{2ex} & & 
\downarrow
& &
\downarrow 
\\
0  \rightarrow &
\Hom_G(L,Z_2) & 
\rightarrow
& \Hom_G(L,E_2) &
\rightarrow & 
\Hom_G(L,L_2) & \rightarrow  0.
\end{array}
\eqno(**)
$$
The connecting homomorphism $\delta$ is determined, by definition, using the
Snake Lemma applied to ($**$).
The more general form of Lemma \ref{lemma-computation-of-extension-class}
states, that the Snake Lemma yields the same connecting homomorphism $\delta$,
even if $Z_1$, $E_1$, and $L_1$ are not injective. 
For simplicity of notation, we proceed to prove 
Lemma \ref{lemma-computation-of-extension-class} in its original formulation,
assuming $Z_1$, $E_1$, and $L_1$ above are injective.
There exists a $G$-equivariant homomorphism $\beta_1:E'\rightarrow E_1$,
such that $\beta_1\circ e=\tilde{e}$, since $E_1$ is injective. 
Let $\alpha_1:Z'\rightarrow Z_1$ be the restriction of $\beta_1$ to $Z'$.
Set $\gamma_1:=\tilde{\ell}:L\rightarrow L_1$. 
Let $\alpha_2:Z''\rightarrow Z_2$, and $\beta_2:E''\rightarrow E_2$ be 
induced by $\alpha_1$ and $\beta_1$. 
Set $\alpha_0:=id:Z\rightarrow Z$. We get a homomorphism $\alpha_\bullet$
from the left column
of (\ref{eq-diagram-of-three-extensions}) to the 
left column of ($*$)
The homomorphism 
$\alpha_{2,*}$ induces the inclusion homomorphism from 
$\coker(\zeta_*)$ to $\coker(\tilde{\zeta}_*)$, 
under the identification of both cokernels with subgroups of
$\Ext^1_G(L,Z)$, as the 
induced homomorphism is $\alpha_{0,*}$. 
Diagram (\ref{eq-snake-lemma-diagram})
maps to ($**$)
via $\alpha_i$, $\beta_i$, $\gamma_1$, and all squares are commutative. 
One now easily checks, that the two connecting homomorphism, 
obtained via the Snake Lemma from (\ref{eq-snake-lemma-diagram}) and ($**$),
are the same.
\EndProof
}

Choose a $G$-equivariant splitting $\psi\in \Hom_G(L,E')$
of the second row of diagram (\ref{eq-diagram-of-three-extensions}), so that 
$j'\circ \psi=id$. Then the extension class 
$\epsilon \in \Ext^1_G(L,Z)$ of the first row of
diagram (\ref{eq-diagram-of-three-extensions}) is given via the equation
\begin{equation}
\label{eq-extension-class-in-terms-of-f}
\epsilon \ := \ \delta(id) \ = \ f[(\iota'')^{-1}\circ \eta\circ \psi],
\end{equation}
by Lemma \ref{lemma-computation-of-extension-class}.

The following Lemma will be used to calculate orders of 
extension classes. 
Let $\sigma:L\rightarrow E$ be a splitting of 
$(\ref{eq-diagram-of-three-extensions})_{top}$, so that
$j\circ \sigma=id$. We do not assume $\sigma$ to be $G$-equivariant.
\begin{new-lemma}
\label{lemma-order-of-extension-class}
Let $k$ be an integer. The class $k\epsilon\in \Ext^1_G(L,Z)$ is
trivial, if and only if
$k(\psi- e\circ \sigma)$ belongs to 
$\Hom(L,Z)+\Hom_G(L,Z').$
\end{new-lemma}

\noindent
{\bf Proof:}
The general case reduces to the case $k=1$ after replacing the top and
middle horizontal extensions in diagram 
(\ref{eq-diagram-of-three-extensions}) by their pullback via the
multiplication map $k:L\rightarrow L$. 
The splittings $\psi$ and $\sigma$ induce splittings $\psi_k$ and
$\sigma_k$ of the pulled-back 
extensions, and a short calculation yields
the equality $\psi_k-e\circ\sigma_k=k(\psi-e\circ\sigma)$.

It suffices prove the case $k=1$.
The equality $j'\circ (\psi-e\circ\sigma)=0$ implies that 
$\psi-e\circ\sigma$ belongs to $\Hom(L,Z')$. Now 
$\zeta\circ[\psi-e\circ\sigma]$ is $G$-equivariant,
since 
$\iota''\circ\zeta\circ[\psi-e\circ\sigma]=
\eta\circ[\psi-e\circ\sigma]=
\eta\circ\psi$. 
Using the latter equality we get
\[
\epsilon 
\ \ \stackrel{(\ref{eq-extension-class-in-terms-of-f})}{=} \ \ 
f((\iota'')^{-1}\circ \eta\circ \psi) = 
f(\zeta\circ [\psi-e\circ\sigma]). 
\]
Now $\ker(f)=\zeta_*[\Hom_G(L,Z')]$ and 
$\zeta\circ [\psi-e\circ\sigma]$ belongs to $\zeta_*[\Hom_G(L,Z')]$,
if and only if
$\psi-e\circ\sigma$ belongs to 
$\Hom(L,Z)+\Hom_G(L,Z')$.
\EndProof

\subsubsection*{The width of an extension class}
This subsection is not needed for the ``generic'' case 
$3\leq i\leq \frac{n}{2}$
of Proposition \ref{prop-extension-of-Q-bar-by-Z-hat-does-not-split}.
Let $A$ be an abelian group. If $A$ is finitely generated, denote by 
\begin{equation}
\label{eq-sharp}
\sharp(A)
\end{equation} 
the cardinality of a set of generators for $A$,
with a minimal number of generators. 
Set $\sharp(A)=\infty$, if $A$ is not finitely generated.

We define next the width of an extension.
Let $\epsilon$ be the class of an extension of $G$-modules 
\begin{equation}
\label{eq-extension-of-L-by-Z}
0\rightarrow Z \rightarrow E \rightarrow L\rightarrow 0.
\end{equation}
Fix a $G$-module $D$ and a $G$-equivariant embedding $z:Z\hookrightarrow D$,
such that $z_*(\epsilon)=0$, and identify $Z$ with $z(Z)$. 
Given a $G$-submodule $A\subset D$, we say that 
$A$ {\em splits} $\epsilon$, if $A$ contains $Z$ and 
the inclusion $Z\rightarrow A$ pushes forward $\epsilon$ to
$0$ in $\Ext^1_G(L,A)$. 
Define the {\em width $\sharp_{D}(\epsilon)$ of $\epsilon$ in $D$} by 
\begin{equation}
\label{eq-width-epsilon-in-D}
\sharp_{D}(\epsilon) \ \ \ := \ \ \ 
\min\{\sharp(A/Z) \ : \ A\subset D \  \mbox{and} \ A \ \mbox{splits} \ 
\epsilon\},
\end{equation}
where $\sharp(A/Z)$ ignores the $G$-module structure of the abelian group 
$A/Z$. Clearly, $\sharp_D(\epsilon)=0$, if and only if $\epsilon=0$.
Pullbacks of extensions can not increase their width.
Note\footnote{
Indeed, $\epsilon$ is the 
image of some homomorphism $\ell\in \Hom_G(L,D/Z)$ via the connecting
homomorphism $\Hom_G(L,D/Z)\rightarrow \Ext^1_G(L,Z)$, and we may
choose $A$ to be the pre-image of ${\rm Im}(\ell)\subset D/Z$ in $D$.},
that $\sharp_{D}(\epsilon)\leq \sharp(L)$. 
For any two injective $\Integers[G]$-modules $I$ and $J$, the equality
$\sharp_{I}(\epsilon)=\sharp_{J}(\epsilon)$ is easily verified.
Consequently, for $I$ injective, $\sharp_{I}(\epsilon)$ is equal to
\begin{equation}
\label{eq-width-in-an-injective-G-module}
\min\{\sharp(A/a(Z)) \ \mid \ a:Z\hookrightarrow A \ \mbox{is a} \ 
G\!-\!\mbox{equivariant embedding and} \ a_*(\epsilon)=0\}.
\end{equation}

When $Z$ is torsion free and the order of $\epsilon$ 
is finite, we set $D=Z_\RationalNumbers:=Z\otimes_\Integers\RationalNumbers$,
and define\footnote{
In our application, 
$G$ is Zariski dense in a semisimple group $G_\RationalNumbers$,
the $G$-modules $Z_\RationalNumbers$ and $L_\RationalNumbers$ will be
restrictions of $G_\RationalNumbers$-modules, and both $Z$ and $L$ 
will be torsion-free. 
Under these assumptions, one can show that 
$\sharp_{Z_\RationalNumbers}(\epsilon)$ is equal to 
(\ref{eq-width-in-an-injective-G-module}). The invariant 
(\ref{eq-width-in-an-injective-G-module})
is a more natural definition of the width $\sharp(\epsilon)$
(avoiding a choice of $D$). We used the equality 
(\ref{eq-width-epsilon}) as a definition, in order to save us its proof,
necessary if the more natural definition of 
$\sharp(\epsilon)$ is used.
} 
the {\em width} of $\epsilon$ by
\begin{equation}
\label{eq-width-epsilon}
\sharp(\epsilon):=\sharp_{Z_\RationalNumbers}(\epsilon).
\end{equation}

Let $\epsilon$, $\sigma$, and $\psi$ be as in 
Lemma \ref{lemma-order-of-extension-class}. 
Assume given a $G$-module $D$ and a $G$-equivariant embedding 
$Z'\hookrightarrow D$.
Given a subgroup $\widehat{Z}$ of $D$ containing $Z$,
set $\widehat{Z}'$ to be the subgroup $Z'+\widehat{Z}$ of $D$.
Let $\nu:Z\rightarrow \widehat{Z}$, $\nu':Z'\rightarrow \widehat{Z}'$ and
$\hat{z}:\widehat{Z}\rightarrow \widehat{Z}'$
be the inclusions. 
%
\begin{new-lemma}
\label{lemma-width}
Assume, that there exist integers $r\geq 1$ and $w\geq 0$ satisfying:
\begin{enumerate}
\item
\label{lemma-cond-pushforward-by-hat-z-is-surjective}
The homomorphism $\hat{z}_*:\Hom_G(L,\widehat{Z})\rightarrow 
\Hom_G(L,\widehat{Z}')$ is surjective, for every 
$G$-submodule $\widehat{Z}$ of $D$, such that 
$\sharp(\widehat{Z}/Z)\leq w$. 
\item
The following inequality holds for $1\leq k <r$.
\begin{equation}
\label{eq-inequality-between-numbers-of-generators}
\sharp\left(
[Z+k\cdot {\rm Im}(\psi-e\circ \sigma)]/Z
\right) \ \ \ > \ \ \ w.
\end{equation}
\end{enumerate}
Then the width $\sharp_D(k\epsilon)$, defined in 
(\ref{eq-width-epsilon-in-D}), satisfies 
$\sharp_D(k\epsilon)>w$, for $1\leq k <r$.
If, furthermore, $r(\psi-e\circ \sigma)$ belongs to $\Hom(L,Z)$,
then the order of the class
$\nu_*\epsilon\in \Ext^1_G(L,\widehat{Z})$ is $r$, for
every $G$-submodule $\widehat{Z}$ of $D$ containing $Z$, such that 
$\sharp(\widehat{Z}/Z)\leq w$.
\end{new-lemma}

\noindent
{\bf Proof:} Let $\widehat{Z}\subset D$ be a $G$-submodule of $D$ 
containing $Z$, such that $\sharp(\widehat{Z}/Z)\leq w$.
Set $\gamma:=\nu'\circ(\psi-e\circ \sigma):L\rightarrow \widehat{Z}'$.
Let 
\[
\begin{array}{ccccccc}
0  \rightarrow & \widehat{Z} & \rightarrow & \widehat{E} & \rightarrow & 
L & \rightarrow 0
\\
 & \hat{z} \downarrow \hspace{1ex} & & \hat{e} \downarrow \hspace{1ex} &&
= \downarrow  \hspace{1ex}
\\
0  \rightarrow & \widehat{Z}' & \rightarrow & \widehat{E}' & \rightarrow & 
L & \rightarrow 0
\end{array}
\]
be the pushforward of the top two horizontal extensions in diagram
(\ref{eq-diagram-of-three-extensions}) via $\nu$ and $\nu'$. 
The class $k\nu_*\epsilon$ is trivial, if and only if $k\gamma$
belongs to 
$\hat{z}_*\Hom(L,\widehat{Z})+\Hom_G(L,\widehat{Z}')$,
by Lemma \ref{lemma-order-of-extension-class}
(applied with the compositions 
$L\RightArrowOf{\psi}E'\rightarrow \widehat{E}'$ and 
$L\RightArrowOf{\sigma} E \rightarrow \widehat{E}$). 
The latter sum is equal to $\hat{z}_*\Hom(L,\widehat{Z})$, since
we assumed the equality
$\Hom_G(L,\widehat{Z})=\Hom_G(L,\widehat{Z}')$.
If $k\gamma$ belongs to $\hat{z}_*\Hom(L,\widehat{Z})$,
then ${\rm Im}(k\gamma)$ is contained in $\widehat{Z}$ and 
the left hand side of (\ref{eq-inequality-between-numbers-of-generators})
is $\leq \sharp(\widehat{Z}/Z)\leq w$, violating 
the inequality (\ref{eq-inequality-between-numbers-of-generators}), unless
$k\equiv 0$ (modulo $r$). 
\EndProof

\begin{example}
\label{example-warm-up-for-non-generic-cases}
{\rm
Let us sketch the proof of Proposition 
\ref{prop-extension-of-Q-bar-by-Z-hat-does-not-split} in the first
non-trivial case $n=3$. 
We need only consider the case $i=2$. 
The extension (\ref{eq-extension-of-Q-bar-by-Z-hat}) becomes
\begin{equation}
\label{eq-n-equal-3-extension-Q-bar-by-Z-hat}
0\rightarrow \ \widehat{Z}^{4} \ \longrightarrow \
H^4(S^{[3]},\Integers) \  \longrightarrow \ \overline{Q}^{4} \ \rightarrow 0.
\end{equation}
We need to prove that the order of 
$\epsilon_{(\ref{eq-n-equal-3-extension-Q-bar-by-Z-hat})}$ is divisible by $4$.
Consider the related extension
\begin{equation}
\label{eq-n-equal-3-extension-Q-4-by-Z-4}
0 \rightarrow \Sym^2H^2(S^{[3]},\Integers) \rightarrow H^4(S^{[3]},\Integers)
\rightarrow Q^4(S^{[3]},\Integers)\rightarrow 0.
\end{equation}
Let $\pi:Q^4(S^{[3]},\Integers)\rightarrow \overline{Q}^{4}$ be the quotient
homomorphism
and $\nu:\Sym^2H^2(S^{[3]},\Integers)\hookrightarrow \widehat{Z}^{4}$
the inclusion. 
Then $\nu_*(\epsilon_{(\ref{eq-n-equal-3-extension-Q-4-by-Z-4})})=
\pi^*(\epsilon_{(\ref{eq-n-equal-3-extension-Q-bar-by-Z-hat})})$,
by Lemma \ref{lemma-characterization-of-pullback-pushforward}.
Now $Q^4(S^{[3]},\Integers)$ has rank $23$, by Lemma
\ref{lemma-10-in-markman-diagonal}, and its torsion 
$Q^4(S^{[3]},\Integers)_{tor}$ is thus cyclic, 
since $\varphi^4:K(S)\rightarrow Q^4(S^{[3]},\Integers)$ is surjective,
by Proposition \ref{prop-varphi-k-is-independent-of-universal-sheaf}. 
$\widehat{Z}^{4}$
is the inverse image of $Q^4(S^{[3]},\Integers)_{tor}$
in $H^4(S^{[3]},\Integers)$. Set $Z:=\Sym^2H^2(S^{[3]},\Integers)$.
Then $\widehat{Z}^{4}/Z$ is isomorphic to $Q^4(S^{[3]},\Integers)_{tor}$
and is thus cyclic. 
Hence, it suffices to prove that the width of $k\cdot 
\epsilon_{(\ref{eq-n-equal-3-extension-Q-4-by-Z-4})}$ satisfies 
$\sharp\left(k\cdot 
\epsilon_{(\ref{eq-n-equal-3-extension-Q-4-by-Z-4})}\right)>1$,
for $1\leq k <4$.

Let $v$ be the class of the ideal sheaf of a length $3$ subscheme, and
\begin{equation}
\label{eq-n-equal-3-extension-v-perp-by-Z}
0\rightarrow Z\rightarrow E \rightarrow v^\perp\rightarrow 0
\end{equation}
the pullback of (\ref{eq-n-equal-3-extension-Q-4-by-Z-4}) via 
$v^\perp\rightarrow K(S)\RightArrowOf{\varphi^4}Q^4(S^{[3]},\Integers)$. 
Note that the latter homomorphism is $G$-equivariant and of positive rank, 
by Lemma \ref{lemma-10-in-markman-diagonal}, and thus injective, 
since $v^\perp$ is irreducible. The homomorphism 
$E\rightarrow H^4(S^{[3]},\Integers)$ is thus injective. 
We identify $E$ with its image.
It suffices to prove that 
$\sharp\left(
k\cdot\epsilon_{(\ref{eq-n-equal-3-extension-v-perp-by-Z})}\right)>1$, 
for $1\leq k <4$. 

Choose the universal sheaf over $S\times S^{[3]}$ to be the ideal sheaf of 
the universal subscheme. Let $u:K(S)\rightarrow K(S^{[3]})$ 
be the homomorphism given in (\ref{eq-u}), and set $u_x:=u(x)$.
Let $\tau:v^\perp\rightarrow v^\perp$ be multiplication by $2$.
Set $\sigma:v^\perp\rightarrow H^4(S^{[3]},\Integers)$
by $\sigma(x):=-2ch_2(u_x)=2c_2(u_x)-c_1(u_x)^2$. 
Then $\sigma$ is an integral splitting of the pullback 
$\tau^*$(\ref{eq-n-equal-3-extension-v-perp-by-Z})
of the extension (\ref{eq-n-equal-3-extension-v-perp-by-Z})
via $\tau$.
Let $w\in v^\perp$ be the class with Mukai vector $(1,0,2)$, 
using the notation (\ref{eq-Mukai-vector}).
Let $a:\Sym^2H^2(S^{[3]},\Integers)\rightarrow 
\Sym^2H^2(S^{[3]},\RationalNumbers)$ be the inclusion.
Define $\psi:v^\perp\rightarrow H^4(S^{[3]},\RationalNumbers)$
by $\psi(x):=2c_2(u_x)-c_1(u_x)^2+\frac{1}{2}c_1(u_w)c_1(u_x)$.
We will prove in section 
\ref{sec-extension-classes-via-universal-Chern-character} that
$\psi$ is a $G$-equivariant splitting of 
$a_*\tau^*$(\ref{eq-n-equal-3-extension-v-perp-by-Z}) 
(Proposition \ref{prop-psi-is-G-equivariant}).

Set $D:=\Sym^2H^2(S^{[3]},\RationalNumbers)$.
The difference $(\psi-\sigma):v^\perp\rightarrow D$ is given by
$(\psi-\sigma)(x)=\frac{1}{2}c_1(u_w)c_1(u_x)$ and $\Hom_G(v^\perp,D)=0$.
We conclude that the order of 
$\tau^*\epsilon_{(\ref{eq-n-equal-3-extension-v-perp-by-Z})}$
is $2$, by Lemma \ref{lemma-order-of-extension-class}. 
But $\tau^*\epsilon_{(\ref{eq-n-equal-3-extension-v-perp-by-Z})}=
2\epsilon_{(\ref{eq-n-equal-3-extension-v-perp-by-Z})}$.
Hence, the order of $\epsilon_{(\ref{eq-n-equal-3-extension-v-perp-by-Z})}$
is $4$. Thus 
$\sharp\left(\epsilon_{(\ref{eq-n-equal-3-extension-v-perp-by-Z})}\right)=
\sharp\left(3\epsilon_{(\ref{eq-n-equal-3-extension-v-perp-by-Z})}\right)\geq 
\sharp\left(2\epsilon_{(\ref{eq-n-equal-3-extension-v-perp-by-Z})}\right)$.
Condition
\ref{lemma-cond-pushforward-by-hat-z-is-surjective} of Lemma 
\ref{lemma-width} is satisfied trivially,
and we have the equality
\[
\sharp\left(
\left[
Z+{\rm Im}(\psi-\sigma)
\right]/Z
\right) \ \ \ = \ \ \ \mbox{rank}(v^\perp) \  = \  23.
\]
Hence, 
$\sharp\left(2\epsilon_{(\ref{eq-n-equal-3-extension-v-perp-by-Z})}\right)=
\sharp\left(\tau^*\epsilon_{(\ref{eq-n-equal-3-extension-v-perp-by-Z})}\right)
\geq 23$, by Lemma \ref{lemma-width}.
\EndProof

%
\hide{
We need a rational splitting of 
(\ref{eq-n-equal-3-extension-v-perp-by-Z}).
Let $Z':=\Sym^2(v^\perp)^*$ be the quotient of 
$(v^\perp)^*\otimes\nolinebreak(v^\perp)^*$
by the saturated anti-symmetric submodule.
Let $\bar{\lambda}$ be a generator of $(v^\perp)^*/v^\perp$. 
Let $t:v^\perp\rightarrow Z'/Z$ be the composition of 
multiplication 
$v^\perp\rightarrow v^\perp\otimes[(v^\perp)^*/v^\perp]$ by $\bar{\lambda}$,  
followed by the natural homomorphism 
$
v^\perp\otimes[(v^\perp)^*/v^\perp] \rightarrow 
Z'/Z.
$
Note that $t$ is $G$-equivariant, since $(v^\perp)^*/v^\perp$ is a trivial 
$G$-module. 
Let $q:Z'\rightarrow  Z'/Z$ is the quotient homomorphism. 
It turns out\footnote{
We omit the proof, since we do not need this fact in the actual proof of 
Proposition \ref{prop-extension-of-Q-bar-by-Z-hat-does-not-split}. 
We state this fact here only in order to present the example 
before the cohomological calculations of section 
\ref{sec-extension-classes-via-universal-Chern-character}.
}
that the $G$-module $E$ above is isomorphic to the kernel in the short exact 
sequence
\[
0\rightarrow E \LongRightArrowOf{e} Z'\oplus v^\perp 
\LongRightArrowOf{(q,t)} Z'/Z\rightarrow 0,
\]
for a suitable choice of $\bar{\lambda}$.
The extension (\ref{eq-n-equal-3-extension-v-perp-by-Z}) fits as the top row
in the commutative diagram
\[
\begin{array}{cccccccc}
0\rightarrow & Z & \rightarrow & E & \rightarrow & v^\perp & \rightarrow 0
\\
& \downarrow & & e \ \downarrow \ \hspace{1ex} & & 
= \ \downarrow \ \hspace{1ex}
\\
0 \rightarrow & Z' & \rightarrow & Z'\oplus v^\perp & \rightarrow & v^\perp
& \rightarrow 0.
\end{array}
\]
Choose a lift of $\bar{\lambda}$ to $\lambda\in (v^\perp)^*$. 
Let $\tilde{\sigma}:v^\perp\rightarrow Z'\oplus v^\perp$ be given by 
$\tilde{\sigma}(x)=(-\lambda x,x)$. 
Then $(q,t)(\tilde{\sigma}(x))=-q(\lambda x)+t(x)=0$. 
Thus, the image of $\tilde{\sigma}$ is contained in $E$ 
and $\tilde{\sigma}=e\circ \sigma$, where $\sigma:v^\perp\rightarrow E$ is a
non-equivariant 
splitting of (\ref{eq-n-equal-3-extension-v-perp-by-Z}). 
Let $\psi:v^\perp\rightarrow Z'\oplus v^\perp$ be the $G$-equivariant
splitting $\psi(x)=(0,x)$. Then $(\psi-e\circ\sigma)(x)=-\lambda x$. 
We get the equality
\[
\sharp\left(
\left[
Z+k \cdot {\rm Im}(\psi-e\circ\sigma)
\right]/Z
\right) \ \ \ = \ \ \ \mbox{rank}(v^\perp) \  = \  23,
\]
for $1\leq k <4$. Set $D:=\Sym^2H^2(S^{[3]},\RationalNumbers)$.
Then $\Hom_G(v^\perp,D)=0$, so condition
\ref{lemma-cond-pushforward-by-hat-z-is-surjective} of Lemma 
\ref{lemma-width} is satisfied trivially.
Hence, the width of $k\epsilon_{(\ref{eq-n-equal-3-extension-v-perp-by-Z})}$
is $\geq 23$, for $1\leq k <4$, by Lemma \ref{lemma-width}.
\EndProof
}
}
\end{example}

%
\subsection{A rational splitting via the universal Chern character}
\label{sec-extension-classes-via-universal-Chern-character}
Set $Q^{2i}:=Q^{2i}(\M,\Integers)$, given in (\ref{eq-C-d}), and 
\[
Z^{2i} \ \ \ := \ \ \ (A_{2i-2})^{2i}/(A_{2i-4})^{2i}.
\]
The quotient homomorphism $(A_{2i-2})^{2i}\rightarrow Z^{2i}$ pushes-forward
(\ref{eq-C-d}) to the short exact sequence of $G$-modules
\begin{equation}
\label{eq-extension-of-Q-by-Z-2i}
0 \rightarrow \ Z^{2i} \ \LongRightArrowOf{} \
E^{2i} \  \LongRightArrowOf{q} \ Q^{2i} \ \rightarrow 0,
\end{equation}
where $E^{2i}:=H^{2i}(\M,\Integers)/(A_{2i-4})^{2i}$ as in 
(\ref{eq-extension-of-Q-bar-by-Z-hat}).
We have the surjective homomorphisms
$\varphi^{2i}: K(S) \rightarrow Q^{2i}$, given in 
(\ref{eq-varphi-k}), 
\begin{equation}
\label{eq-mult}
m \ : \ 
Q^{2i-2}(\M,\Integers) \otimes_\Integers H^2(\M,\Integers) \ \ 
\longrightarrow \ \ Z^{2i},
\end{equation}
induced by cup-product, and the surjective composition
\begin{equation}
\label{eq-surjective-composition-of-m-with-varphi-2i-2}
K(S)\otimes_\Integers H^2(S,\Integers) \ \ 
\LongRightArrowOf{\varphi^{2i-2}\otimes id} \ \ 
Q^{2i-2}(\M,\Integers)\otimes_\Integers H^2(S,\Integers) \ \ 
\LongRightArrowOf{m}  \ \ Z^{2i}.
\end{equation}

Let $L$ be the Mukai lattice $K(S)$ or one of its 
sublattices $\Integers v$ or $v^\perp$. 
Denote by $(\ref{eq-extension-of-L-by-Z-2i})_L$ the pullback 
\begin{equation}
\label{eq-extension-of-L-by-Z-2i}
0 \rightarrow Z^{2i} \RightArrowOf{\iota} E^{2i}_L \RightArrowOf{j} L 
\rightarrow 0
\end{equation}
of (\ref{eq-extension-of-Q-by-Z-2i}) 
via the restriction of $\varphi^{2i}$ 
to $L$, so that $E^{2i}_L\subset E^{2i}\oplus L$ is the kernel of 
$
q-\restricted{\varphi^{2i}}{L}  : 
E^{2i}\oplus L \rightarrow Q^{2i},
$
$\iota$ is the composition 
$Z^{2i}\hookrightarrow E^{2i}\hookrightarrow E^{2i}\oplus L$, 
sending $x$ to $(x,0)$, and $j$ is the composition of the inclusion
$E^{2i}_L\hookrightarrow E^{2i}\oplus L$ with the projection to $L$. 

We apply Lemma \ref{lemma-order-of-extension-class} to calculate the order of
$\epsilon_{(\ref{eq-extension-of-L-by-Z-2i})}$.
We construct first a non-equivariant splitting 
$\sigma:L\rightarrow E^{2i}_L$ of 
(\ref{eq-extension-of-L-by-Z-2i}), for $i\geq 3$.
It suffices to construct a homomorphism
\begin{equation}
\label{eq-bar-sigma}
\bar{\sigma} \ : \  K(S) \ \ \ \longrightarrow \ \ \ 
E^{2i},
\end{equation}
such that $q\circ \bar{\sigma}=\varphi^{2i}$. 
Let $u:K(S)\rightarrow K(\M)$ be the homomorphism given in (\ref{eq-u})
and set $u_x:=u(x)$.
Let $\bar{\sigma}(x)$ be the coset 
\[
\bar{\sigma}(x) \ \ := \ \ 
\left\{
\begin{array}{ccc}
c_i(u_x) - c_1(u_x)c_{i-1}(u_x)
\ \  + \ \  (A_{2i-4})^{2i},
& \mbox{if} & i\geq 3,
\\
2c_2(u_x)-c_1(u_x)^2,
& \mbox{if} & i = 2.
\end{array}
\right.
\]
in $E^{2i}$.
Define the non-equivariant homomorphism 
$\sigma:L\rightarrow E^{2i}\oplus L$ by
\begin{equation}
\label{eq-sigma-non-equivariant-splitting}
\sigma(x) \ \ \ := \ \ \ 
\left\{
\begin{array}{ccc}
(\bar{\sigma}(x),x), 
& \mbox{if} & i\geq 3,
\\
 (\bar{\sigma}(x),2x), 
& \mbox{if} & i=2.
\end{array}
\right.
\end{equation}
The following lemma implies that $\sigma(x)$ has values in 
$E^{2i}_L$, for $i\geq 2$, and that it is a splitting for
$i\geq 3$.

\begin{new-lemma}
\label{lemma-sigma-is-linear}
The map $\bar{\sigma}$ is a linear homomorphism 
and \ \ 
$
q\circ \bar{\sigma} \ \ = \ \ 
\left\{
\begin{array}{ccc}
\varphi^{2i}, & \mbox{if} & i\geq 3,
\\
2\varphi^4, & \mbox{if} & i=2.
\end{array}
\right.
$
\end{new-lemma}

\noindent
{\bf Proof:}
Only the linearity of $\bar{\sigma}$ is non-obvious. 
If $i=2$, then $\bar{\sigma}(x)=-2ch_2(x)$ is clearly linear.
Assume that $i\geq 3$. 
The formula for $\bar{\sigma}$ uses the first two terms of 
$(-1)^{i-1}(i-1)!ch_i$, as a polynomial in the Chern classes. 
Some of the remaining coefficients are rational, so we compute
explicitly using 
the product formula $c(x+y)=c(x)c(y)$:\\
$
[c_i(x+y)-c_{i-1}(x+y)c_1(x+y)]-
[c_i(x)-c_{i-1}(x)c_1(x)]-[c_i(y)-c_{i-1}(y)c_1(y)]=\\
\epsilon[c_2(x)c_{i-2}(y)+\cdots +c_{i-2}(x)c_2(y)]
-c_1(x+y)[c_1(x)c_{i-2}(y)+\cdots+c_{i-2}(x)c_1(y)]\in (A_{2i-4})^{2i},
$\\
where $\epsilon=0$, for $i=3$, and $\epsilon=1$, for $i\geq 4$.
\EndProof

\smallskip
Next we push-forward the extension (\ref{eq-extension-of-L-by-Z-2i})
to an extension (\ref{eq-extension-of-L-by-widetilde-Z-2i}) and
show that it admits an equivariant splitting 
(Proposition \ref{prop-psi-is-G-equivariant}). 
Define $\widetilde{Z}^{2i}$ as follows:

Case $i=2$: Set $\widetilde{Z}^4:=\Sym^2[H^2(\M,\Integers)^*]$; 
the quotient of $[H^2(\M,\Integers)^*]\otimes_\Integers[H^2(\M,\Integers)^*]$ 
by the saturated anti-symmetric submodule.

Case $i\geq 3$:
Regard $H^2(\M,\Integers)$ as a subgroup of $H^2(\M,\Integers)^*$, via
the integral Beauville-Bogomolov pairing.
Let 
\[
b \ : \ Q^{2i-2}(\M,\Integers)\otimes_\Integers H^2(\M,\Integers)
\ \ \ \rightarrow \ \ \ 
Q^{2i-2}(\M,\Integers)\otimes_\Integers 
\{H^2(\M,\Integers)+(i-1)H^2(\M,\Integers)^*\}
\]
be the natural homomorphism. 
Set 
\[
\widetilde{Z}^{2i} \ \ := \ \ \left[
Q^{2i-2}(\M,\Integers)\otimes_\Integers 
\{H^2(\M,\Integers)+(i-1)H^2(\M,\Integers)^*\}
\right]/b(\ker(m)),
\]
where $m$ is given in equation (\ref{eq-mult}).

Let 
\begin{equation}
\label{eq-homomorphism-z}
z:Z^{2i}\rightarrow \widetilde{Z}^{2i}
\end{equation}
 be the induced
homomorphism,
\begin{equation}
\label{eq-extension-of-Q-by-widetilde-Z-2i}
0 \rightarrow \widetilde{Z}^{2i} \rightarrow \widetilde{E}^{2i} 
\RightArrowOf{\tilde{q}} Q^{2i} \rightarrow 0
\end{equation}
the push-forward of (\ref{eq-extension-of-Q-by-Z-2i}), 
and 
\begin{equation}
\label{eq-extension-of-L-by-widetilde-Z-2i}
0 \rightarrow \widetilde{Z}^{2i} \RightArrowOf{\iota} \widetilde{E}^{2i}_L 
\RightArrowOf{j} L \rightarrow 0
\end{equation}
the push-forward of (\ref{eq-extension-of-L-by-Z-2i}) via $z$.

The formulae (\ref{eq-psi-bar}) and (\ref{eq-psi-a-G-equivariant-splitting}),
for the $G$-equivariant splitting of 
(\ref{eq-extension-of-L-by-widetilde-Z-2i}), will involve the class $w$, 
which we now define.
Let $w\in v^\perp$ be the class in $K(S)$, orthogonal to $v$,
which is determined by the equation $c_1(u_w)=c_1(u_v)$. 
The existence and uniqueness of $w$ follows from Theorem 
\ref{thm-irreducibility} via the equality 
$c_1(u_x)=\theta_v(x)$, for all $x\in v^\perp$,
where $\theta_v$ is Mukai's isomorphism (\ref{eq-theta-v-from-v-perp}). 
The class $w$ depends on the choice of the universal class 
$[\E]$, since $c_1(u_v)$ does. 
Set $K(S,\RationalNumbers):=K(S)\otimes_\Integers\RationalNumbers$.

\begin{new-lemma}
\label{lemma-v-perp-and-w-over-length-square-of-v-generate}
The pairing of the class $\frac{w}{(v,v)}$ with elements of $v^\perp$
has integral values. 
Furthermore, $v^\perp$ and $\frac{w}{(v,v)}$ span the image of the embedding 
$(v^\perp)^*\hookrightarrow K(S,\RationalNumbers)$, 
induced by the Mukai pairing. If $v=(1,0,1-n)$ and $\E$ is the
universal ideal sheaf over $S\times S^{[n]}$, then $w=(1,0,n-1)$.
\end{new-lemma}

\journal{The easy proof can be found in the preprint version of this paper
\cite{markman-eprint-version}.
}

\preprint{
\noindent
{\bf Proof:} 
We show first that 
the coset $w+(v,v)v^\perp$ is independent of the choice of $[\E]$.
If we replace $[\E]$ by $[\E]\cup [F]$, where $F$ is a complex line-bundle 
on $\M$, then $c_1(u_v)$ gets replaced by $c_1(u_v)+(v,v)c_1(F)$,
while $c_1(u_w)$ is determined by $w$, and is independent of $F$.
For the equality $c_1(u_w)=c_1(u_v)$ to hold, 
we need to replace $w$ by $w+(v,v)f$, where $f\in v^\perp$ is
the class, such that $c_1(u_f)=c_1(F)$. 

The isometry group $O^+K(S)$ acts transitively on the 
set of primitive classes $v$, with fixed $(v,v)$. 
Let $v'\in K^{1,1}(S)$ be another primitive and effective class, 
satisfying $(v,v)=(v',v')$, $H'$
a $v'$-generic polarization, and set $\M':=\M_{H'}(v')$.
Choose $g\in O^+K(S)$, satisfying $g(v)=v'$, and
universal classes $[\E]\in K(S\times \M)$ and $[\E']\in K(S\times \M')$. 
Then there is a parallel-transport operator 
$p_g:K(\M)\rightarrow K(\M')$
and a topological complex line bundle $\ell_g$ on $\M_{H'}(v')$,
such that $(g\otimes p_g)([\E])=[\E']\cup f_2^*\ell_g$, by 
\cite{markman-monodromy-I}, Theorem 1.2 part 4 and Theorem 1.6. Hence, 
it suffices to prove the lemma for one choice of $v$,
satisfying $(v,v)=2n-2$, for each integer $n\geq 2$, by 
the independence of $w+(v,v)v^\perp$ of the choice of $[\E]$.

Let $v=(1,0,1-n)$ be the Mukai vector of the ideal sheaf of a length $n$
subscheme of $S$, so that $\M=S^{[n]}$, and 
$\E$ the universal ideal sheaf over $S\times S^{[n]}$. 
Then $c_1(\E)=0$. 
Let $y\in K(S)$ be the class of the sky-scraper sheaf over a point 
$s\in S$.
Mukai's notation for $y$ is $(0,0,1)$. 
Then $y=y^\vee$ and $f_1^!(y)=\iota_!(1)$, where 
$\iota:S^{[n]}\hookrightarrow S\times S^{[n]}$ is the 
embedding via $x\mapsto (s,x)$. We get the equalities
\[
u_y:=f_{2_!}(f^!(-y^\vee)\cup [\E])=-f_{2_!}(\iota_!(1)\cup [\E])=
-f_{2_!}(\iota_!(\iota^![\E]))=-\iota^![\E].
\]
Thus $c_1(u_{(0,0,1)})=-\iota^!c_1([\E])=0.$
It follows, that $w$ in this case is the primitive class 
$(1,0,n-1)=v+(v,v)(0,0,1)$ 
orthogonal to $v$. The rest of the proof is a straightforward calculation.
\EndProof
}

\smallskip
The following proposition would provide the desired splitting of
(\ref{eq-extension-of-L-by-widetilde-Z-2i}). 
The proposition is the {\em key} to the proof of the main technical result, 
Proposition
\ref{prop-extension-of-Q-bar-by-Z-hat-does-not-split}.
The group $H^2(\M,\Integers)^*$ is spanned by 
$H^2(\M,\Integers)$ and $\frac{c_1(u_w)}{(v,v)}$, by Lemma 
\ref{lemma-v-perp-and-w-over-length-square-of-v-generate}.
Define the homomorphism
\[
\bar{\psi} \ : \ K(S) \ \ \ \longrightarrow \ \ \ 
\widetilde{E}^{2i},
\]
sending $x\in K(S)$ to the coset $\bar{\psi}(x)$ of 
\begin{equation}
\label{eq-psi-bar}
\left\{
\begin{array}{ccc}
c_i(u_x) - c_1(u_x)c_{i-1}(u_x) +
\frac{(i-1)}{(v,v)}c_1(u_w)c_{i-1}(u_x),
& \mbox{if} & i\geq 3
\\
2c_2(u_x) - c_1(u_x)^2 +
\frac{2}{(v,v)}c_1(u_w)c_{1}(u_x) -
\frac{(v,x)}{(v,v)^2}c_1(u_w)^2,
& \mbox{if} & i=2
\end{array}
\right.
\end{equation}
in $\widetilde{E}^{2i}$. The coset is constructed as follows:
the first two terms, in each of the two cases above, 
determine a coset in $E^{2i}$, and the rest of the terms 
determine a coset in $\widetilde{Z}^{2i}$, by Lemma
\ref{lemma-v-perp-and-w-over-length-square-of-v-generate}, 
and $\widetilde{E}^{2i}$
is the quotient of $E^{2i}\oplus \widetilde{Z}^{2i}$ by the image of
$Z^{2i}\hookrightarrow E^{2i}\oplus \widetilde{Z}^{2i}$,
sending $x$ to $(x,-z(x))$.
The linearity of $\bar{\psi}$ follows from Lemma
\ref{lemma-sigma-is-linear}. 

\begin{prop}
\label{prop-psi-is-G-equivariant}
The homomorphism $\bar{\psi}$ 
is $G$-equivariant and satisfies 
\[
\tilde{q}\circ \bar{\psi} \ \ = \ \ 
\left\{
\begin{array}{ccc}
\varphi^{2i}, & \mbox{if} & i\geq 3,
\\
2\varphi^{4}, & \mbox{if} & i=2,
\end{array}
\right.
\] 
where $\tilde{q}$ is given in (\ref{eq-extension-of-Q-by-widetilde-Z-2i}). 
\end{prop}

The proof will require the following lemma. 

\begin{new-lemma}
\label{lemma-c-i-of-alpha-tensor-a-line-bundle}
(\cite{markman-original-archive-version}, Lemma 9)
Let $X$ be a topological space, $i>0$ an integer, 
$\alpha$ a class of rank $r$ in $K(X)$, and
$\ell\in K(X)$ the class of a line-bundle.
Then 
\[
c_i(\alpha\otimes \ell) \ \ \ \equiv \ \ \ 
\left\{
\begin{array}{ccc}
c_i(\alpha) + (r+1-i)c_{i-1}(\alpha)c_1(\ell), 
& \mbox{if} & i\neq 2,
\\
c_2(\alpha) + (r-1)c_1(\alpha)c_1(\ell)+\frac{r(r-1)}{2}c_1(\ell)^2,
& \mbox{if} & i=2,
\end{array}
\right.
\]
modulo the subalgebra $A_{2i-4}(X)$, generated by classes of degree 
$\leq 2i-4$.
\end{new-lemma}


{\bf Proof of Proposition
\ref{prop-psi-is-G-equivariant}:} 
The $G$-equivariance is the only non-obvious
statement. 
The natural homomorphism 
$K(\M)\rightarrow 
K(\M)\otimes_\Integers\RationalNumbers$ is injective,
since $H^*(\M,\Integers)$ is torsion free 
(Theorem \ref{thm-integral-generators}).
Denote by
\[
\monrep_g \ : \ K(\M) \ \ \ \longrightarrow \ \ \ K(\M),
\]
$g\in G$, the restriction to $K(\M)$ of the conjugate
of 
\[
\monrep(g) \ : \ H^*(\M,\RationalNumbers) \ \ \ \longrightarrow \ \ \ 
H^*(\M,\RationalNumbers)
\]
via the Chern character isomorphism
$ch:K(\M)\otimes_\Integers\RationalNumbers \rightarrow 
H^*(\M,\RationalNumbers)$. Theorem 
\ref{thm-summary-of-monodromy-results} part 
\ref{thm-item-mon-g-sends-a-universal-classes-to-such}
implies the equality
\[
(g\otimes \monrep_g)([\E_v]) \ \ \ = \ \ \ [\E_v]\cup f_2^!\ell_g,
\]
for $g\in G$, and for the class $\ell_g\in K(\M)$ 
of some topological complex line-bundle.
We get the equality 
\begin{equation}
\label{eq-mon-g-of-e-g-1-x}
\monrep_g(u_{g^{-1}(x)}) \ \ \ = \ \ \ u_x\cup \ell_g,
\end{equation}
for $x\in K(S)$. Set 
$\bar{\ell}_g:=c_1(\ell_g)\in H^2(\M,\Integers)$.
Then
\begin{equation}
\label{eq-mon-g-of-c-1-e-x}
\monrep_g(c_1(u_{g^{-1}(x)})) \ \ \ = \ \ \ c_1(u_x)+(v,x)\bar{\ell}_g,
\end{equation}
since $\rank(u_x)=(v,x)$. In particular, 
$\monrep_g(c_1(u_v))=c_1(u_v)+(v,v)\bar{\ell}_g$, and thus also
\begin{equation}
\label{eq-mon-g-of-c-1-e-w}
\monrep_g(c_1(u_w)) \ \ \ = \ \ \ c_1(u_w)+(v,v)\bar{\ell}_g.
\end{equation}

The polynomial maps 
$c_i:K(\M)\rightarrow H^*(\M,\Integers)$,
defined by the Chern classes, are $Mon(\M)$-equivariant. 
For $i\neq 2$ we get the equivalence
\begin{equation}
\label{eq-mon-g-of-c-i-1-g-inverse-x}
\monrep_g(c_{i-1}(u_{g^{-1}(x)})) \ = \  
c_{i-1}(\monrep_g(u_{g^{-1}(x)})) \ 
\stackrel{(\ref{eq-mon-g-of-e-g-1-x})}{=} \ 
c_{i-1}(u_x\cup\ell_g) \ \equiv \
c_{i-1}(u_x),
\end{equation}
where the last equivalence is modulo $(A_{2i-4})^{2i-2}$. 
Similarly, we get
\begin{eqnarray}
\label{eq-mon-g-c-i-g-inverse-x}
\monrep_g(c_{i}(u_{g^{-1}(x)}) &
\stackrel{(\ref{eq-mon-g-of-e-g-1-x})}{=} &
c_{i}(u_x\cup\ell_g) \ \equiv \ 
\\
& \equiv &
\left\{
\begin{array}{ccc}
c_{i}(u_x)+[(v,x)+1-i]c_{i-1}(u_x)\bar{\ell}_g,
& \mbox{if} & i\neq 2, 
\\
\nonumber
c_2(u_x) + [(v,x)-1]c_1(u_x)\bar{\ell}_g +
\frac{(v,x)[(v,x)-1]}{2}\bar{\ell}_g^2,
& \mbox{if} & i = 2, 
\end{array}
\right.
\end{eqnarray}
where the last equivalence is taken modulo $(A_{2i-4})^{2i}$, 
and follows from 
Lemma \ref{lemma-c-i-of-alpha-tensor-a-line-bundle}.

The equality $\monrep_g(\bar{\psi}(g^{-1}(x))=\bar{\psi}(x)$
now follows by a lengthy but 
straightforward substitution of each of the four Chern 
classes in (\ref{eq-psi-bar}) via the the corresponding 
equation among (\ref{eq-mon-g-of-c-1-e-x}), 
(\ref{eq-mon-g-of-c-1-e-w}), 
(\ref{eq-mon-g-of-c-i-1-g-inverse-x}), and 
(\ref{eq-mon-g-c-i-g-inverse-x}). 
\EndProof

\medskip
Define 
\begin{equation}
\label{eq-psi-a-G-equivariant-splitting}
\psi \ : \  L \ \ \ \hookrightarrow \ \ \ \widetilde{E}^{2i}_L
\end{equation}
by $\psi(x) = (\bar{\psi}(x),x) \in \widetilde{E}^{2i}\oplus L$,
if $i\geq 3$, and $\psi(x) = (\bar{\psi}(x),2x)$, if $i=2$.
Then $\psi$ has values in $\widetilde{E}^{2i}_L$, for $i\geq 2$, 
by Proposition \ref{prop-psi-is-G-equivariant}, and $\psi$ is 
a $G$-equivariant splitting of (\ref{eq-extension-of-L-by-widetilde-Z-2i}), 
for $i\geq 3$.

Let $\widetilde{\sigma}:L\rightarrow \widetilde{E}^{2i}_L$ 
be the composition of
$\sigma$, given in (\ref{eq-sigma-non-equivariant-splitting}), 
with the homomorphism $E^{2i}_L\rightarrow \widetilde{E}^{2i}_L$. 
The difference
\[
[\psi-\widetilde{\sigma}] \ : \ L \ \ \ \longrightarrow \ \ \ 
\widetilde{Z}^{2i}
\]
is given by sending 
$x$ to the coset $\left[\frac{i-1}{(v,v)}c_{i-1}(u_x)c_1(u_w)\right]$,
for $i\geq 3$. 
We get the equality
\begin{equation}
\label{eq-psi-minus-e-sigma-in-terms-of-varphi-2i-2}
[\psi-\widetilde{\sigma}](x) \ \ \ = \ \ \ 
\left\{
\begin{array}{ccc}
\left[\frac{i-1}{(v,v)}\restricted{\varphi^{2i-2}}{L}(x)c_1(u_w)\right],
& \mbox{if} & i\geq 3,
\\
\frac{2}{(v,v)}\restricted{\varphi^2}{L}(x)c_1(u_w)-
\frac{(v,x)}{(v,v)^2}c_1(u_w)^2,
& \mbox{if} & i=2.
\end{array}
\right.
\end{equation}

%
\subsection{Proof of the generic case of Proposition 
\ref{prop-extension-of-Q-bar-by-Z-hat-does-not-split}}
\label{sec-proof-of-generic-case}
We prove Proposition 
\ref{prop-extension-of-Q-bar-by-Z-hat-does-not-split} in case
$3\leq i\leq \frac{n}{2}$, since the proof is simple and 
illustrates the general idea. 

\begin{new-lemma}
\label{lemma-if-A-mod-v-perp-is-small}
Let $A\subset v^\perp\otimes \RationalNumbers$ be a $G$-invariant submodule 
containing $v^\perp$, such that $A/v^\perp$ is finite and can be generated 
by less that $23=\rank(v^\perp)$ elements. Then the natural homomorphism
$\Hom_G(v^\perp,v^\perp)\rightarrow \Hom_G(v^\perp,A)$ is an isomorphism.
\end{new-lemma}

\noindent 
{\bf Proof:}
$\Hom_G(v^\perp,v^\perp\otimes \RationalNumbers)=\RationalNumbers\cdot id$.
If $v^\perp\subset \lambda v^\perp$, for some $\lambda\in \RationalNumbers$, 
then $\lambda=\frac{1}{k}$, for some non-zero integer $k$. 
If $\Abs{k}\neq 1$, then $(\lambda v^\perp)/v^\perp$ is not generated by 
less than $23$ 
elements. Thus the image of $\Hom_G(v^\perp,A)$ in
$\Hom_G(v^\perp,v^\perp\otimes \RationalNumbers)$ must be equal to the image 
of $\Hom_G(v^\perp,v^\perp)$.
\EndProof

Part \ref{prop-item-summary-order-of-extension-is-geq-3} of
Proposition \ref{prop-extension-of-Q-bar-by-Z-hat-does-not-split}
follows from 
part \ref{prop-item-order-of-extension-class-of-Q-bar-by-Z-hat}.
In the generic case $3\leq i\leq \frac{n}{2}$ the group
$Q^{2i}(\M,\Integers)$ is torsion free and the extensions 
(\ref{eq-extension-of-Q-bar-by-Z-hat}) and (\ref{eq-extension-of-Q-by-Z-2i})
coincide. Extension 
$(\ref{eq-extension-of-L-by-Z-2i})_{v^\perp}$, 
with $L=v^\perp$, is the pullback of 
(\ref{eq-extension-of-Q-by-Z-2i}). 
Hence, it suffices to prove that the order of 
$\epsilon_{(\ref{eq-extension-of-L-by-Z-2i})_{v^\perp}}$
is divisible by $\frac{2n-2}{\gcd(i-1,2n-2)}$.

The order of $\epsilon_{(\ref{eq-extension-of-L-by-Z-2i})_{v^\perp}}$
is equal to the order of the coset of 
$[\psi-\widetilde{\sigma}]$ in 
\[
\Hom(v^\perp,\widetilde{Z}^{2i})/\left[
\Hom(v^\perp,Z^{2i})+\Hom_G(v^\perp,\widetilde{Z}^{2i})
\right],
\] 
by Lemma \ref{lemma-order-of-extension-class}.
The homomorphism $K(S)\otimes v^\perp\rightarrow Z^{2i}$, given in
(\ref{eq-surjective-composition-of-m-with-varphi-2i-2}), 
is an isomorphism in our case, by Lemma \ref{lemma-10-in-markman-diagonal}.
Clearly, $\Hom_G(v^\perp,K(S)\otimes v^\perp)=
\Hom_G(v^\perp,v\otimes v^\perp)$.
We conclude the equality
$\Hom_G(v^\perp,\widetilde{Z}^{2i})=\Hom_G(v^\perp,Z^{2i})$, 
by Lemma \ref{lemma-if-A-mod-v-perp-is-small}.
Hence, the order of the coset of $[\psi-\widetilde{\sigma}]$
is $\frac{(v,v)}{\gcd(i-1,(v,v))}$, by Equation 
(\ref{eq-psi-minus-e-sigma-in-terms-of-varphi-2i-2}).
This completes the proof, since $(v,v)=2n-2$.
\EndProof

%
\subsection{The order of three extension classes}
\label{sec-proof-of-prop-three-non-trivial-extensions}
Proposition \ref{prop-extension-of-Q-bar-by-Z-hat-does-not-split}
is a consequence of the following more detailed proposition. 
Consider the push-forward
\begin{equation}
\label{eq-a-wide-extension}
0 \rightarrow Z^{2i}/tor \rightarrow E^{2i}_{v^\perp}/tor \rightarrow v^\perp 
\rightarrow 0
\end{equation}
of $(\ref{eq-extension-of-L-by-Z-2i})_{v^\perp}$, with $L=v^\perp$, 
via $Z^{2i}\rightarrow Z^{2i}/tor$. 
The width 
$\sharp\left(\epsilon_{(\ref{eq-a-wide-extension})}\right)$ 
of the extension was defined
in equation (\ref{eq-width-epsilon}).
The extension $(\ref{eq-extension-of-L-by-Z-2i})_{\Integers v}$, 
with $L=\Integers v$, splits $G$-equivariantly 
(the proof
is omitted, as we do not use this fact). In contrast, we have:

\begin{prop}
\label{prop-three-non-trivial-extensions}
Set $n:=\dim_{\ComplexNumbers}(\M)/2$ and assume that $n\geq 2$. 

\noindent
\newcounter{prop-item-generic-i}
\setcounter{prop-item-generic-i}{1}
\Alph{prop-item-generic-i}) 
The following statements hold for any 
integer $i$ in the range $3\leq i \leq \frac{n+2}{2}$.
\begin{enumerate}
\item
$\frac{2n-2}{\gcd(i-1,2n-2)}$ divides the order of the 
$G$-modules extension class of 
(\ref{eq-extension-of-Q-by-Z-2i}).
\item
The order of $\epsilon_{(\ref{eq-extension-of-L-by-Z-2i})_{v^\perp}}$
is divisible by $\frac{2n-2}{\gcd(i-1,2n-2)}$.
\item
\label{prop-item-width}
The width $\sharp\left(k\epsilon_{(\ref{eq-a-wide-extension})}\right)$
is $> 20$, for any integer $k$ between $1$ and 
$\frac{2n-2}{\gcd(i-1,2n-2)}$.
\end{enumerate}
\newcounter{prop-item-i-equal-2}
\setcounter{prop-item-i-equal-2}{2}
\Alph{prop-item-i-equal-2}) When $i=2$, then $n-1$ divides the order of 
$2\epsilon_{(\ref{eq-extension-of-Q-by-Z-2i})}$.
If $i=2$, then 
the order of $2\epsilon_{(\ref{eq-extension-of-L-by-Z-2i})_{v^\perp}}$
is equal to $n-1$. 
\hide{
If $(n,i,L)=(3,2,v^\perp)$, then the width 
$\sharp\left(2\epsilon_{(\ref{eq-extension-of-L-by-Z-2i})}\right)$
is $23$. 
}
\end{prop}

The following Lemmas will be needed in the proof of Proposition
\ref{prop-three-non-trivial-extensions}.

Let $\widehat{B}$ be the saturation of 
$B:=\Integers v\otimes v^\perp+[v^\perp\otimes v^\perp]^G$ in
$K(S)\otimes v^\perp$. 
Let $\eta:v^\perp\oplus \Integers v\rightarrow B$ be the isomorphism sending 
$(x,0)$ to $v\otimes x$ and
$(0,v)$ to the image of $-(v,v)\cdot id_{v^\perp}$ via 
$(v,v)\cdot [(v^\perp)^*\otimes v^\perp]^G\cong [v^\perp\otimes v^\perp]^G$.
\journal{The routine proof of the following Lemma 
can be found in \cite{markman-eprint-version}.}

\begin{new-lemma}
\label{lemma-mukai-lattice-in-K-top-S-otimes-v-perp}
The quotient $\widehat{B}/B$ is cyclic of order $(v,v)$. 
The isomorphism $\eta$ extends to a $G$-equivariant isomorphism $\tilde{\eta}$
from $K(S)$ onto $\widehat{B}$. 
\end{new-lemma}

\preprint{
\noindent 
{\bf Proof:} $\widehat{B}/B$ is isomorphic to the torsion subgroup of the 
quotient
\[
\frac{K(S)\otimes v^\perp}{
[\Integers v\otimes v^\perp]+[v^\perp\otimes v^\perp]^G} \cong
\frac{(v^\perp)^*\otimes v^\perp}{[v^\perp\otimes v^\perp]^G}
\supset \frac{[(v^\perp)^*\otimes v^\perp]^G}{[v^\perp\otimes v^\perp]^G}
\cong \frac{\Integers}{(v,v)\Integers}.
\]
Thus, both $\widehat{B}/B$ 
and $K(S)/[v^\perp\oplus \Integers v]$ 
are cyclic of order $(v,v)$ and 
the extension $\tilde{\eta}$ is an isomorphism, if it exists.

We may assume, that $v=(1,0,1-n)$. Set $w:=(1,0,n-1)$ and 
$H:=H^2(S,\Integers)$. Then 
$id_{v^\perp}=-\frac{w\otimes w}{(v,v)}+id_H$. 
Note, that $w-v=(v,v)(0,0,1)$ and $(0,0,1)$ maps to a generator of
$K(S)/[v^\perp\oplus \Integers v]$. 
Now $\eta(w-v)=v\otimes w+(v,v)id_{v^\perp}=(v-w)\otimes w+(v,v)id_H=
(v,v)\left[(0,0,-1)\otimes w+id_H\right]$.
Hence, $\tilde{\eta}(0,0,1)$ 
is an integral element of $\widehat{B}$.
\EndProof
}

\begin{new-lemma}
\label{lemma-equivariant-Hom-v-perp-Z-equal-Hom-v-perp-tilde-Z}
Let $A\subset [K(S)\otimes v^\perp]$ 
be a saturated $G$-submodule of rank $\leq 24$ 
and $\widehat{A}$ its saturation in $K(S)\otimes (v^\perp)^*$.
Set 
$Z:=[K(S)\otimes v^\perp]/A$  and 
$Z':=[K(S)\otimes\nolinebreak{(v^\perp)^*]/\widehat{A}}$. 
Let $Y\subset Z\otimes\RationalNumbers$ be a subgroup containing $Z$
such that $Y/Z$ is finite and $\sharp(Y/Z)\leq 21$ ($\sharp$ is defined in
(\ref{eq-sharp})).
Then $\Hom_G(v^\perp,Z)\rightarrow \Hom_G(v^\perp,Y+Z')$ 
is an isomorphism.
\end{new-lemma}

\noindent 
{\bf Proof:} If $A$ contains $\Integers v\otimes v^\perp$, then 
both $\Hom_G(v^\perp,Z)$ and $\Hom_G(v^\perp,Y+Z')$ vanish.
Otherwise, $A=[v^\perp\otimes v^\perp]^G$. 
Let $\widehat{I}$ be the saturation of the image $I$ of 
$\Integers v\otimes v^\perp$ in $Z$.
Then $\widehat{I}/I\cong \Integers/(v,v)\Integers$, by Lemma 
\ref{lemma-mukai-lattice-in-K-top-S-otimes-v-perp}.
The image $I'$ of 
$\Integers v\otimes (v^\perp)^*$ in $Z'$ is saturated, contains $I$, and
$I'/I$ is cyclic of order $(v,v)$. 
Thus, $\widehat{I}$ maps isomorphically onto $I'$. 
Let $\widehat{I'}$  be the saturation of $I'$ in $Y+Z'$. 
Then\\
$
\sharp(\widehat{I'}/I) \leq 
\sharp(\widehat{I'}/I')+\sharp(I'/I) \leq
\sharp([Y+Z']/Z')+1 = 
\sharp(Y/[Y\cap Z'])+1 \leq \sharp(Y/Z)+1\leq 22<23.
$\\
Thus, the composition 
$\Hom_G(v^\perp,I)\rightarrow 
\Hom_G(v^\perp,\widehat{I})\rightarrow \Hom_G(v^\perp,\widehat{I'})$
is an isomorphism, by Lemma \ref{lemma-if-A-mod-v-perp-is-small}.
The statement follows, since
$\Hom_G(v^\perp,\widehat{I})\rightarrow\Hom_G(v^\perp,Z)$ and
$\Hom_G(v^\perp,\widehat{I'})\rightarrow\Hom_G(v^\perp,Y+Z')$ 
are both isomorphisms. 
\EndProof

\journal{The routine proof of the following Lemma 
can be found in \cite{markman-eprint-version}.}

\begin{new-lemma}
\label{lemma-torsion-subgroup-is-cyclic}
\begin{enumerate}
\item
\label{lemma-item-torsion-subgroup-is-cyclic}
Let $A\subset [K(S)\otimes v^\perp]$ be a $G$-invariant 
saturated subgroup of rank $\leq 24$. Then the torsion subgroup
$T_A$ of $[K(S)\otimes v^\perp]/[A+(v^\perp\otimes w)]$
is cyclic, where $w\in v^\perp$ is the class in Lemma 
\ref{lemma-v-perp-and-w-over-length-square-of-v-generate}. 
\item
\label{lemma-item-sym-2-plus-w-times-v-perp-is-saturated}
Let $\Wedge{2}(v^\perp)$ be the saturated anti-symmetric subgroup
of $v^\perp\otimes v^\perp$. Then the subgroup 
$\Wedge{2}(v^\perp)+[v^\perp\otimes w]$ is saturated in 
$v^\perp\otimes v^\perp$.
\end{enumerate}
\end{new-lemma}

\preprint{
\noindent
{\bf Proof:} \ref{lemma-item-torsion-subgroup-is-cyclic})
$A$ is necessarily a subgroup of the saturation of 
$B:=[\Integers v\otimes v^\perp]+[v^\perp\otimes v^\perp]^G$. 
The saturation $\widehat{B}$ of $B$ in $K(S)\otimes v^\perp$ intersects
$v^\perp\otimes w$ trivially. Thus the kernel of
$\frac{K(S)\otimes v^\perp}
{\widehat{B}+(v^\perp\otimes w)}\rightarrow 
\frac{K(S)\otimes v^\perp}{A+(v^\perp\otimes w)}$ is isomorphic to 
$\widehat{B}/A$, which is torsion free. So $T_A$ injects into 
$T_{\widehat{B}}$. It suffices to treat the case $A=\widehat{B}$. 
In that case, we have
\[
\frac{K(S)\otimes v^\perp}{\widehat{B}+[v^\perp\otimes w]} \cong
\frac{(v^\perp)^*\otimes v^\perp}{
[(v^\perp)^*\otimes v^\perp]^G+(v^\perp\otimes w)}
\supset \frac{[H\otimes H]^G+\Integers\frac{w\otimes w}{(v,v)}}{
[(v^\perp)^*\otimes v^\perp]^G+\Integers(w\otimes w)} \cong
\Integers/(v,v)\Integers,
\]
where the left isomorphism is due to Lemma 
\ref{lemma-mukai-lattice-in-K-top-S-otimes-v-perp}, and 
the right isomorphism is due to the equality
$id_{v^\perp}=id_H-\frac{w\otimes w}{(v,v)}$ in $(v^\perp)^*\otimes v^\perp$.
The quotient by this torsion subgroup is clearly torsion free, and 
$T_{\widehat{B}}$ is thus cyclic of order $(v,v)$. 

\ref{lemma-item-sym-2-plus-w-times-v-perp-is-saturated})
Set $v:=(1,0,1-n)$, $w:=(1,0,n-1)$, and $H:=H^2(S,\Integers)$.
Then 
\[
(v^\perp\otimes v^\perp)/[\Wedge{2}v^\perp+(v^\perp\otimes w)]\cong
\Sym^2(v^\perp)/[v^\perp\otimes w]\cong \Sym^2H
\] 
is torsion free.
\EndProof

Note, that if $A$ is $[v\otimes v^\perp]+\Wedge{2}v^\perp$, then
$T_A$ is a direct sum of $22$ copies of $\Integers/(v,v)\Integers$.
Hence, some condition on the rank of $A$ is necessary in Lemma
\ref{lemma-torsion-subgroup-is-cyclic}.
}

{\bf Proof of Proposition
\ref{prop-three-non-trivial-extensions}:}\\
\newcounter{step-order-of-extension-of-L-by-Z-2i}
\setcounter{step-order-of-extension-of-L-by-Z-2i}{1}
{\bf Step \Roman{step-order-of-extension-of-L-by-Z-2i}:} 
We consider first the order of 
$\epsilon_{(\ref{eq-extension-of-L-by-Z-2i})_{v^\perp}}$ and
the width of $k\epsilon_{(\ref{eq-a-wide-extension})}$.

{\bf Case $i=2$:} 
$Z^4=\Sym^2H^2(\M,\Integers)$ and 
$\widetilde{Z}^4=
\Sym^2[H^2(\M,\Integers)^*]$, which is equal to 
$\Sym^2[H^2(\M,\Integers)+\Integers\frac{1}{(v,v)}c_1(u_w)]$, by Lemma
\ref{lemma-v-perp-and-w-over-length-square-of-v-generate}. 
Consider the difference $\psi-\widetilde{\sigma}$ given in equation
(\ref{eq-psi-minus-e-sigma-in-terms-of-varphi-2i-2}), with $L=v^\perp$.
The order of $2\epsilon_{(\ref{eq-extension-of-L-by-Z-2i})_{v^\perp}}$
is equal to the order of the coset of $\psi-\widetilde{\sigma}$ in
$
\Hom(v^\perp,\widetilde{Z}^4)/[\Hom(v^\perp,Z^4)+
\Hom_G(v^\perp,\widetilde{Z}^4)],
$
by Lemma \ref{lemma-order-of-extension-class}.
$\Hom_G(v^\perp,\widetilde{Z}^4)=0$, and 
$[\psi-\widetilde{\sigma}](w)=\frac{2}{(v,v)}c_1(u_w)^2$. 
We conclude, that the order of the coset of $[\psi-\widetilde{\sigma}]$ is 
$\frac{(v,v)}{2}$, which is also the order of 
$2\epsilon_{(\ref{eq-extension-of-L-by-Z-2i})_{v^\perp}}$.

\hide{
We calculate the width 
$\sharp\left(2\epsilon_{(\ref{eq-extension-of-L-by-Z-2i})}\right)$ 
for $(n,i,L)=(3,2,v^\perp)$.
Let $h$ be the surjective composition
$v^\perp\otimes v^\perp\rightarrow K(S)\otimes v^\perp
\RightArrowOf{(\ref{eq-surjective-composition-of-m-with-varphi-2i-2})} Z^{4}$.
The image $h(v^\perp\otimes w)$ is saturated, by Lemma
\ref{lemma-torsion-subgroup-is-cyclic}
part \ref{lemma-item-sym-2-plus-w-times-v-perp-is-saturated}.
Let $\gamma:v^\perp\rightarrow \widetilde{Z}^4$ be the restriction of 
$\psi-\widetilde{\sigma}$ to $v^\perp$. It is given by 
$\gamma(x)=\frac{1}{2}c_1(u_x)c_1(u_w)=\frac{1}{2}h(x)$. 
We get
$
\sharp\left(\frac{{\rm Im}(\gamma)+Z^4}{Z^4} \right)
 = 
\sharp\left(\frac{(1/2)h(v^\perp\otimes w)}{h(v^\perp\otimes w)}\right)
= 23. 
$
Now $\Hom_G(v^\perp,Z^4_\RationalNumbers)$ vanishes. 
Lemma \ref{lemma-width}
implies that the width 
$\sharp\left(2\epsilon_{(\ref{eq-extension-of-L-by-Z-2i})}\right)$
is $\geq 23$. The equality 
$\sharp\left(2\epsilon_{(\ref{eq-extension-of-L-by-Z-2i})}\right)=23$
follows, since
$\sharp\left(2\epsilon_{(\ref{eq-extension-of-L-by-Z-2i})}\right)
\leq \sharp(v^\perp)=23$.
}

{\bf Case $i\geq 3$:}
We provide next a lower bound for the rank of $\psi-\widetilde{\sigma}$.
Let $\bar{h}$ be the surjective composition
$
K(S)\otimes v^\perp \ 
\LongRightArrowOf{
(\ref{eq-surjective-composition-of-m-with-varphi-2i-2})} \ Z^{2i} \ 
\longrightarrow \ Z^{2i}/tor.
$
Lemma \ref{lemma-10-in-markman-diagonal} implies the following:

If $i\leq \frac{n}{2}$, then $\bar{h}$ is an isomorphism (and $Z^{2i}$
is torsion free).

If $i= \frac{n+1}{2}$, then  
$\rank(\ker(\bar{h}))+\rank(\ker(\varphi^{n+1}))\leq 1$.

If $i= \frac{n+2}{2}$, then 
$\rank(\ker(\bar{h}))+\rank(\ker(\varphi^{n+2}))\leq 24$.

\noindent
The decomposition of $(v^\perp\otimes v^\perp)_\RationalNumbers$
into irreducible $G$-submodules is 
\begin{equation}
\label{eq-decomposition-of-v-perp-squared}
(v^\perp\otimes v^\perp)_\RationalNumbers^G \ \ \oplus \ \ 
\Wedge{2}{v^\perp_\RationalNumbers} \ \ 
\oplus \ \ v^\perp_\RationalNumbers(2),
\end{equation} 
where the latter is the
submodule of $\Sym^2v^\perp_\RationalNumbers$ 
spanned by squares of isotropic vectors.
We conclude that $\bar{h}$ maps $v^\perp(2)$, 
given in (\ref{eq-decomposition-of-v-perp-squared}), 
injectively into $Z^{2i}/tor$.

Set $L=v^\perp$. 
Let $\gamma$ be the image of $\psi-\widetilde{\sigma}$ 
in $\Hom(v^\perp,\widetilde{Z}^{2i}/tor)$. 
It is given by
$
\gamma(x) 
\stackrel{(\ref{eq-psi-minus-e-sigma-in-terms-of-varphi-2i-2})}{=} 
\left(\frac{i-1}{(v,v)}\right)\bar{h}(x\otimes w),
$
for $x\in v^\perp$.
The kernel of $\gamma$ has rank $\leq 1$, since $\bar{h}$ maps 
$v^\perp(2)$ injectively\footnote{
If $i=\frac{n+3}{2}$, then 
$\rank(\ker(\bar{h}))+\rank(\ker(\varphi^{n+3}))\leq 323$.
This follows from Lemma 
\ref{lemma-10-in-markman-diagonal}, if $n\geq 9$, 
and from the table in the appendix of
\cite{markman-diagonal}, for $n=3, 5, 7$. 
Hence $\bar{h}$ maps $\Wedge{2}v^\perp$ or $v^\perp(2)$ injectively
and $\rank(\ker(\gamma))\leq 1$.
This may help settle some more cases of Conjecture 
\ref{conj-Mon-isomorphic-to-Mon-2}.
}  
into $Z^{2i}/tor$. 
The kernel $A$ of $\bar{h}$ has rank $\leq 24$.
Let $B$ be the saturation of $\bar{h}(v^\perp\otimes w)$ in 
$Z^{2i}/tor$. Then $B/\bar{h}(v^\perp\otimes w)$ is cyclic,
by Lemma \ref{lemma-torsion-subgroup-is-cyclic}
part \ref{lemma-item-torsion-subgroup-is-cyclic}. 
If $\bar{h}(v^\perp\otimes w)$ is saturated, we 
get the isomorphism: \\
${\displaystyle
\frac{{\rm Im}(k\gamma)+Z^{2i}/tor}{Z^{2i}/tor} \ \ \ \cong \ \ \
\left(\left[k\frac{i-1}{(v,v)}\Integers\right]/\Integers\right)\otimes
{\rm Im}(\gamma).
}$\\
In any case, the minimal number $\sharp$ of generators, 
defined in (\ref{eq-sharp}), satisfies
\[
\sharp\left[
\frac{{\rm Im}(k\gamma)+Z^{2i}/tor}{Z^{2i}/tor}
\right]  \ \ \ \geq \ \ \ \rank(\gamma)-1 \ \ \ \geq \ \ \ 21,
\]
for $1\leq k < \frac{(v,v)}{\gcd(i-1,(v,v))}$.

The homomorphism 
$\Hom_G(v^\perp,Z^{2i}/tor)\rightarrow 
\Hom_G(v^\perp,\widetilde{Z}^{2i}/tor+Y)$
is an isomorphism, for any $G$-submodule 
$Y\subset Z^{2i}\otimes\RationalNumbers$ containing $Z^{2i}/tor$, such that 
$\sharp(Y/(Z^{2i}/tor))\leq 21$,
by Lemma \ref{lemma-equivariant-Hom-v-perp-Z-equal-Hom-v-perp-tilde-Z}.
We conclude, that the order of the extension class of 
$\epsilon_{(\ref{eq-a-wide-extension})}$ is 
$\frac{(v,v)}{\gcd(i-1,(v,v))}$
and the width $\sharp(k\epsilon_{(\ref{eq-a-wide-extension})})$
of $k\epsilon_{(\ref{eq-a-wide-extension})}$
is $> 20$, for $1\leq k < \frac{(v,v)}{\gcd(i-1,(v,v))}$, by Lemma
\ref{lemma-width}
(applied with $w=20$, $r=\frac{(v,v)}{\gcd(i-1,(v,v))}$, 
$D=Z^{2i}\otimes\RationalNumbers$, $Z=Z^{2i}/tor$, and 
$Z'=\widetilde{Z}^{2i}/tor$).
Hence, $\frac{(v,v)}{\gcd(i-1,(v,v))}$ divides also the order of 
$\epsilon_{(\ref{eq-extension-of-L-by-Z-2i})_{v^\perp}}$.

\newcounter{step-order-of-extension-of-Q-by-Z-2i}
\setcounter{step-order-of-extension-of-Q-by-Z-2i}{2}
{\bf Step \Roman{step-order-of-extension-of-Q-by-Z-2i}:} 
(The order of $\epsilon_{(\ref{eq-extension-of-Q-by-Z-2i})}$)
Let $d$ be the the order of
$\epsilon_{(\ref{eq-extension-of-L-by-Z-2i})_{v^\perp}}$. 
The class  $\epsilon_{(\ref{eq-extension-of-L-by-Z-2i})_{v^\perp}}$ 
is a pullback of 
$\epsilon_{(\ref{eq-extension-of-Q-by-Z-2i})}$. 
Hence, the order
of $\epsilon_{(\ref{eq-extension-of-Q-by-Z-2i})}$
is divisible by $d$.
We have shown in step \Roman{step-order-of-extension-of-L-by-Z-2i}, that 
$\frac{(v,v)}{\gcd((i-1),(v,v))}$ divides $d$, if $i\geq 3$. 
If $i=2$, we have shown in step 
\Roman{step-order-of-extension-of-L-by-Z-2i}, that $n-1$ divides the order of 
$2\epsilon_{(\ref{eq-extension-of-L-by-Z-2i})_{v^\perp}}$. 
So $n-1$ divides the order of 
$2\epsilon_{(\ref{eq-extension-of-Q-by-Z-2i})}$.
\EndProof

%
\subsection{Proof of Proposition  
\ref{prop-extension-of-Q-bar-by-Z-hat-does-not-split} }
\label{sec-reduction-of-prop-extension-to-prop}

The following Lemmas will be needed for the proof of Proposition 
\ref{prop-extension-of-Q-bar-by-Z-hat-does-not-split}.

\begin{new-lemma}
\label{lemma-class-of-extension-of-K-top-S-mod-L-by-L-mod-L'}
Let $L\subset K(S)$ be $\Integers v$ or $v^\perp$
and $L'$ a $G$-invariant sublattice of $L$, such
that $L/L'$ is torsion. Let $d$ be the exponent of $L/L'$.
The homomorphism $L\rightarrow L/L'$ pushes the extension class of
\begin{equation}
\label{eq-extension-of-K-top-S-mod-L-by-L}
0 \rightarrow L \rightarrow K(S) \rightarrow K(S)/L \rightarrow 0
\end{equation}
to a class in 
$\Ext^1_G(K(S)/L,L/L')$ of order $\gcd(d,(v,v))$.
\end{new-lemma}

\preprint{We postpone the long proof of 
Lemma \ref{lemma-class-of-extension-of-K-top-S-mod-L-by-L-mod-L'}
to section 
\ref{sec-proof-of-lemma-class-of-extension-of-K-top-S-mod-L-by-L-mod-L'}.
}
\journal{Lemma \ref{lemma-class-of-extension-of-K-top-S-mod-L-by-L-mod-L'}
is proven in the preprint version of this paper
\cite{markman-eprint-version}. 
The long proof is elementary and thus omitted.
}

\begin{new-lemma}
\label{lemma-if-varphi-2i-does-not-vanish}
Let $L\subset K(S)$ be $\Integers v$ or $v^\perp$ and $i$ 
an integer $\geq 2$. 
The order of the $G$-module extension class of 
\begin{equation}
\label{eq-extension-of-Q-mod-varphi-L-by-varphi-L}
0 \rightarrow \ \varphi^{2i}(L) \ 
\rightarrow \ Q^{2i} \ 
\rightarrow \ Q^{2i}/\varphi^{2i}(L) \ 
\rightarrow 0
\end{equation}
is $\gcd(d,(v,v))$, if the exponent $d$ of
$\varphi^{2i}(L)$ is finite, and the order is $(v,v)$, if
$\varphi^{2i}(L)$ has positive rank.
\end{new-lemma}

\noindent
{\bf Proof:}
The assumption $i\geq 2$ implies that $\varphi^{2i}$ is $G$-equivariant,
by Theorem \ref{thm-summary-of-monodromy-results} part
\ref{thm-item-markman-monodromy-I-lemma-4.5}.
Set $L':= L\cap\ker(\varphi^{2i})$. 
The push-forward of the extension 
(\ref{eq-extension-of-K-top-S-mod-L-by-L})
via $L\rightarrow L/L'$ 
is isomorphic to the pullback of the extension
(\ref{eq-extension-of-Q-mod-varphi-L-by-varphi-L})
via the homomorphism 
$K(S)/L\rightarrow Q^{2i}/\varphi^{2i}(L)$
induced by $\varphi^{2i}$ 
(use Lemma \ref{lemma-characterization-of-pullback-pushforward}). 
The pushed-forward extension has order $\gcd(d,(v,v))$, by Lemma
\ref{lemma-class-of-extension-of-K-top-S-mod-L-by-L-mod-L'}.
Hence, the extension class 
$\epsilon_{(\ref{eq-extension-of-Q-mod-varphi-L-by-varphi-L})}$ 
pulls back to a class of order $\gcd(d,(v,v))$.
We conclude, that the order of 
$\epsilon_{(\ref{eq-extension-of-Q-mod-varphi-L-by-varphi-L})}$ 
is divisible by $\gcd(d,(v,v))$. 

The order of  $\epsilon_{(\ref{eq-extension-of-Q-mod-varphi-L-by-varphi-L})}$ 
clearly divides $d$. 
It remains to show that the order of 
$\epsilon_{(\ref{eq-extension-of-Q-mod-varphi-L-by-varphi-L})}$ 
divides $(v,v)$.
Let $\nabla$ be the sublattice $\Integers v+v^\perp+\ker(\varphi^{2i})$ of
$K(S)$. The extension 
(\ref{eq-extension-of-Q-mod-varphi-L-by-varphi-L})
is the pullback of 
\[
0\rightarrow L/L'
\rightarrow 
K(S)/\left[L^\perp+\ker(\varphi^{2i})\right]
\rightarrow 
K(S)/\nabla
\rightarrow 0,
\]
via the homomorphism 
$Q^{2i}/\varphi^{2i}(L) \rightarrow K(S)/\nabla$, by Lemma
\ref{lemma-characterization-of-pullback-pushforward}. 
The order of $K(S)/\nabla$ divides the order
$(v,v)$ of $K(S)/[\Integers v + v^\perp]$. 
Hence, the order of 
$\epsilon_{(\ref{eq-extension-of-Q-mod-varphi-L-by-varphi-L})}$ 
divides $(v,v)$.
\EndProof

{\bf Proof of Proposition  
\ref{prop-extension-of-Q-bar-by-Z-hat-does-not-split}:}
Part \ref{prop-item-summary-order-of-extension-is-geq-3})
Part \ref{prop-item-summary-order-of-extension-is-geq-3} 
follows from 
part \ref{prop-item-order-of-extension-class-of-Q-bar-by-Z-hat}.
This is clear for $i=2$, so we check only for $i\geq 3$,
where the order of $\epsilon_{(\ref{eq-extension-of-Q-bar-by-Z-hat})}$
is shown in part \ref{prop-item-order-of-extension-class-of-Q-bar-by-Z-hat} 
to be $r(n,i):=\frac{2n\!-\!2}{\gcd(i\!-\!1,2n\!-\!2)}$.
If $i\leq\frac{n+1}{2}$, then $r(n,i)\geq 4$. 
If $i=\frac{n+2}{2}$ then $n$ is even and we have:
$
r(n,i) \ \ = \ \ 
\left\{\begin{array}{ccc}
2n\!-\!2 & \mbox{if} & 
i=\frac{n+2}{2} \ \mbox{is even}
\\
n\!-\!1 & \mbox{if} & 
i=\frac{n+2}{2} \ \mbox{is odd}.
\end{array}\right.
$\\
Thus, 
$r(n,i)\geq 3$, for $n\geq 3$ and $i$ an integer in the 
range $2\leq i \leq \frac{n+2}{2}$. 

Part
\ref{prop-item-order-of-extension-class-of-Q-bar-by-Z-hat})
If $2\leq i \leq \frac{n}{2}$, then $\rank(Q^{2i})=24$, 
and if $i=\frac{n+1}{2}$ and $n$ is odd, then 
$\rank(Q^{2i})\geq 23$, by Lemma \ref{lemma-10-in-markman-diagonal}.
Hence, we need only consider the following cases.

\newcounter{case-rank-Q-2i-is-24}
\setcounter{case-rank-Q-2i-is-24}{1}
{\bf Case \Roman{case-rank-Q-2i-is-24}:}
If $3\leq i \leq \frac{n}{2}$, or $\rank(Q^{2i})=24$, 
then $Q^{2i}$ is torsion free, by
Proposition \ref{prop-varphi-k-is-independent-of-universal-sheaf}.
Thus, the
exact sequences (\ref{eq-extension-of-Q-bar-by-Z-hat}) and
(\ref{eq-extension-of-Q-by-Z-2i}) coincide. 
The order of the extension class of (\ref{eq-extension-of-Q-by-Z-2i})
is divisible by $\frac{2n-2}{\gcd(i-1,2n-2)}$, by Proposition 
\ref{prop-three-non-trivial-extensions}.

\newcounter{case-rank-Q-2i-is-23}
\setcounter{case-rank-Q-2i-is-23}{2}
{\bf Case \Roman{case-rank-Q-2i-is-23}:} 
The case $\frac{n+1}{2}\leq i \leq \frac{n+2}{2}$, $n\geq 3$, $i\geq 3$, 
and $\rank(Q^{2i})=23$:
The pullback of (\ref{eq-extension-of-Q-bar-by-Z-hat}) 
via $v^\perp\rightarrow \overline{Q}^{2i}$ is an extension
\begin{equation}
\label{eq-extension-of-v-perp-by-Z-hat}
0\rightarrow \widehat{Z}^{2i} \rightarrow
\widehat{E}^{2i}\rightarrow
v^\perp \rightarrow 0.
\end{equation}
It suffices to prove that $\frac{2n\!-\!2}{\gcd(i\!-\!1,2n\!-\!2)}$
divides the order of $\epsilon_{(\ref{eq-extension-of-v-perp-by-Z-hat})}$. 
Note, that the extension (\ref{eq-extension-of-v-perp-by-Z-hat})
is the pushforward of $(\ref{eq-extension-of-L-by-Z-2i})_{v^\perp}$
via the natural inclusion
$Z^{2i}\rightarrow \widehat{Z}^{2i}$, by
Lemma \ref{lemma-characterization-of-pullback-pushforward}. 

The quotient $\widehat{Z}^{2i}/Z^{2i}$ is isomorphic to $Q^{2i}_{tor}$,
and $Q^{2i}_{tor}$ is equal to the image $\varphi^{2i}(\Integers v)$,
since $Q^{2i}$ is assumed to have rank $23$.
The quotient $(\widehat{Z}^{2i}/tor)/(Z^{2i}/tor)$ is itself a quotient of
$\widehat{Z}^{2i}/Z^{2i}$ and is hence cyclic.
The width of $k\epsilon_{(\ref{eq-a-wide-extension})}$
is $>1$, for $1\leq k< \frac{2n\!-\!2}{\gcd(i\!-\!1,2n\!-\!2)}$, 
by part \Alph{prop-item-generic-i}.\ref{prop-item-width} of Proposition 
\ref{prop-three-non-trivial-extensions}.
Thus, the pushforward of $\epsilon_{(\ref{eq-a-wide-extension})}$
via $Z^{2i}/tor\rightarrow \widehat{Z}^{2i}/tor$ has 
order $\frac{2n\!-\!2}{\gcd(i\!-\!1,2n\!-\!2)}$. 
The latter coincides with the
push-forward of $\epsilon_{(\ref{eq-extension-of-v-perp-by-Z-hat})}$
via $\widehat{Z}^{2i}\rightarrow \widehat{Z}^{2i}/tor$. 
We conclude that the order of
$\epsilon_{(\ref{eq-extension-of-v-perp-by-Z-hat})}$, and consequently of
$\epsilon_{(\ref{eq-extension-of-Q-bar-by-Z-hat})}$, is divisible by
$\frac{2n\!-\!2}{\gcd(i\!-\!1,2n\!-\!2)}$.

\newcounter{case-rank-Q-2i-is-1}
\setcounter{case-rank-Q-2i-is-1}{3}
{\bf Case \Roman{case-rank-Q-2i-is-1}:} 
The case $i=\frac{n+2}{2}$, $n$ even and $\geq 4$, and $\rank(Q^{2i})=1$.
The kernel of $\varphi^{2i}$ is contained in $L:=v^\perp$,
since $\varphi^{2i}$ is $G$-equivariant.
Then $Q_{tor}^{2i}=\varphi^{2i}(v^\perp)$.
Consequently, the quotient homomorphism
$\widehat{Z}^{2i}\rightarrow \widehat{Z}^{2i}/Z^{2i}\cong Q^{2i}_{tor}$ 
pushes-forward
the extension (\ref{eq-extension-of-Q-bar-by-Z-hat})
to the extension 
(\ref{eq-extension-of-Q-mod-varphi-L-by-varphi-L}).
We show next, that the order of the extension class of 
(\ref{eq-extension-of-Q-mod-varphi-L-by-varphi-L})
is divisible by $\frac{2n-2}{\gcd(i-1,2n-2)}$.
Hence, so is the order of the $G$-module extension class of 
(\ref{eq-extension-of-Q-bar-by-Z-hat}).

Let $d$ be the greatest common divisor of $(v,v)$ and the order of
$\epsilon_{(\ref{eq-extension-of-L-by-Z-2i})_{v^\perp}}$.
Then $d$
divides the order of $\restricted{\varphi^{2i}}{v^\perp}$, 
which is equal to the exponent of $\varphi^{2i}(v^\perp)$.
Hence, the exponent of $\varphi^{2i}(v^\perp)$
is divisible by  $d$ (or infinite).
The order of 
$\epsilon_{(\ref{eq-extension-of-Q-mod-varphi-L-by-varphi-L})_{v^\perp}}$
is thus divisible by $d$, 
by Lemma \ref{lemma-if-varphi-2i-does-not-vanish}. 
Now $d$ is divisible by
$\frac{(v,v)}{\gcd(i-1,(v,v))}$, by Proposition
\ref{prop-three-non-trivial-extensions}.

\newcounter{case-i-is-2-of-prop-extension-of-Q-bar-by-Z-hat}
\setcounter{case-i-is-2-of-prop-extension-of-Q-bar-by-Z-hat}{4}
{\bf Case \Roman{case-i-is-2-of-prop-extension-of-Q-bar-by-Z-hat}:} 
The case $i=2$ of Proposition 
\ref{prop-extension-of-Q-bar-by-Z-hat-does-not-split}:
When $i=2$ and $n\geq 4$, then $Q^{2i}$ is torsion free, the
exact sequences (\ref{eq-extension-of-Q-bar-by-Z-hat}) and
(\ref{eq-extension-of-Q-by-Z-2i}) coincide, and 
the order of the extension class 
$2\epsilon_{(\ref{eq-extension-of-Q-by-Z-2i})}$
is divisible by $n-1$, by Proposition 
\ref{prop-three-non-trivial-extensions}.
The case $i=2$ and $n=3$ of the Proposition was proven
in Example \ref{example-warm-up-for-non-generic-cases}.
\hide{
When $i=2$ and $n=3$ 
then $\ker(\varphi^{2i})$ is contained in
$\Integers v$, by Lemma \ref{lemma-10-in-markman-diagonal}.
Set $L:=v^\perp$. Then the order of 
$2\epsilon_{(\ref{eq-extension-of-L-by-Z-2i})}$ 
is  $2$, by Proposition 
\ref{prop-three-non-trivial-extensions}.
Thus, the order of $\epsilon_{(\ref{eq-extension-of-L-by-Z-2i})}$ is
$4$. 
The width $\sharp\left(2\epsilon_{(\ref{eq-extension-of-L-by-Z-2i})}\right)$
is $>1$, by Proposition \ref{prop-three-non-trivial-extensions}.
We conclude, that the order of 
$2\epsilon_{(\ref{eq-extension-of-Q-bar-by-Z-hat})}$ is divisible by $2$,
as in case \Roman{case-rank-Q-2i-is-23} 
above (with the simplification that $Z^4$ is torsion free, so 
(\ref{eq-extension-of-v-perp-by-Z-hat}) and (\ref{eq-extension-of-L-by-Z-2i})
coincide). Thus, the order of 
$\epsilon_{(\ref{eq-extension-of-Q-bar-by-Z-hat})}$ is divisible by $4$.
}
%
This completes the proof of
Proposition \ref{prop-extension-of-Q-bar-by-Z-hat-does-not-split}
\EndProof

%

%

%
\preprint{
\subsection{Proof of 
Lemma \ref{lemma-class-of-extension-of-K-top-S-mod-L-by-L-mod-L'}
}
\label{sec-proof-of-lemma-class-of-extension-of-K-top-S-mod-L-by-L-mod-L'}

We keep the notation of section 
\ref{subsec-Monodromy-constaints-via-non-split-extensions}.
In particular, $v$ is a primitive Mukai vector, 
$(v,v)\geq 2$, and $G:=O^+K(S)_v$. 
The following Lemma follows easily from results of Nikulin
\cite{nikulin}.

\begin{new-lemma} (\cite{markman-monodromy-I} Lemma 8.1)
\label{lemma-primitive-embeddings-of-rank-2-lattices}
Let $\Lambda$ be an even unimodular lattice 
of signature $(\ell_+,\ell_-)$, satisfying $\ell_+\geq 3$ and $\ell_-\geq 3$. 
\begin{enumerate}
\item
Let $M=
\left[
\begin{array}{cc}
2a & b \\
b & 2d
\end{array}
\right]$ 
be a symmetric matrix with $a,b, d\in \Integers$, and  
$\lambda_1\in \Lambda$  a primitive element with 
$(\lambda_1,\lambda_1)=2a$.
Then there exists a primitive rank $2$ sublattice $\Sigma\subset \Lambda$, 
containing $\lambda_1$, and an element $\lambda_2\in \Sigma$, such that 
$\{\lambda_1,\lambda_2\}$ is a basis for $\Sigma$, and 
$M$ is the matrix of the bilinear form of 
$\Sigma$ is this basis. 
\item 
\label{lemma-item-transitivity-on-pairs}
Assume that $\rank(M)=2$. Let $\{\lambda_1',\lambda_2'\}$ 
be a basis for another primitive sublattice of $\Lambda$,
having the same matrix $M$. 
Then there exists an isometry $g$ of $\Lambda$ satisfying 
$g(\lambda_i)=\lambda_i'$, $i=1,2$.
\item 
\label{lemma-item-isometry-acts-transitivly-on-primitive-vectors}
The kernel $O^+(\Lambda)$ of the orientation character 
acts transitively on the set of primitive integral classes 
$\lambda\in \Lambda$, 
with fixed ``squared-length'' $(\lambda,\lambda)$. 
\end{enumerate}
\end{new-lemma}

\begin{new-lemma}
\label{lemma-E}
The $OH^2(S,\Integers)$-orbit of every primitive class in $H^2(S,\Integers)$
spans $H^2(S,\Integers)$. Consequently, every $OH^2(S,\Integers)$-invariant 
sublattice of $H^2(S,\Integers)$ is of the form
$dH^2(S,\Integers)$, for some integer $d$. Furthermore, 
$[H^2(S,\Integers)/dH^2(S,\Integers)]^{OH^2(S,\Integers)}=(0).$
\end{new-lemma}

\noindent
{\bf Proof:}
Let $U$ be the hyperbolic plane with the bilinear form 
$((x_1,y_1),(x_2,y_2))=x_1y_2+x_2y_1.$
$U\oplus U$ is an orthogonal direct summand of $H^2(S,\Integers)$. 
For every integer $n$, the lattice $U\oplus U$ is spanned by the vectors
\[
\{[(1,n),(0,0)], \ \ [(1,n),(1,0)], \ \ [(1,n),(0,1)], \ \ [(0,0),(1,n)], \ \ 
[(1,0),(1,n)], \ \ [(0,1),(1,n)]\},
\]
all of which are primitive vectors of self-intersection $2n$. 
All the above vectors belong to a single $OH^2(S,\Integers)$-orbit
in $H^2(S,\Integers)$, by part
\ref{lemma-item-isometry-acts-transitivly-on-primitive-vectors} of
Lemma 
\ref{lemma-primitive-embeddings-of-rank-2-lattices}.
Let $\ell$ be a primitive element of $H^2(S,\Integers)$ and $L$
the sublattice of $H^2(S,\Integers)$ spanned by
its $OH^2(S,\Integers)$-orbit. Set $n:=(\ell,\ell)/2$.
There exists an isometry $g$ of $H^2(S,\Integers)$, satisfying
$g(\ell)=[(1,n),(0,0)]$. Hence, $L$ contains the summand $U\oplus U$.
The $OH^2(S,\Integers)$-orbit of 
every element $u$ of $H^2(S,\Integers)$ intersects $U\oplus U$
non-trivially. Hence, $u$ belongs to $L$ and $L=H^2(S,\Integers)$.
\EndProof

\smallskip
Let $w$ be the Mukai vector $(1,0,n-1)$, $n\geq 2$.
Given integers $d$, $e$, we set 
\begin{equation}
\label{eq-L-d-e}
L_{d,e} \ \ \ := \ \ \ dH^2(S,\Integers)\oplus e\Integers w.
\end{equation}

\begin{new-lemma}
\label{lemma-B}
Assume $v=(1,0,1-n)$, $n\geq 2$. 
Let $L\subset v^\perp$ be a $G$-invariant sub-lattice, such that 
$v^\perp/L$ is finite. 
Then $L=L_{d,e}$, for some positive integers $d$, $e$. 
\end{new-lemma}

\noindent
{\bf Proof:} A $G$-invariant sublattice of $v^\perp$ is also 
$OH^2(S,\Integers)$-invariant. 
$v^\perp$ is the orthogonal direct sum $H^2(S,\Integers)\oplus \Integers w$.
The image of the projection of $L$ to $H^2(S,\Integers)$ 
is $dH^2(S,\Integers)$, for some positive integer $d$, by Lemma
\ref{lemma-E}. The kernel is spanned by $ew$, for some positive integer $e$.
Let $x\in H^2(S,\Integers)$ be a primitive element. 
Choose an integer $y$, such that
$\tilde{x}:=dx+yw$ belongs to $L$.
There exists an isometry $g\in OH^2(S,\Integers)$, such that $x-g(x)$
is primitive. Then $\tilde{x}-g(\tilde{x})=d(x-g(x))$
belongs to $L\cap H^2(S,\Integers)$ and the $OH^2(S,\Integers)$-orbit
of $d(x-g(x))$ spans $dH^2(S,\Integers)$, by Lemma \ref{lemma-E}. 
Thus, $dH^2(S,\Integers)$ is contained in $L$ and $L=L_{d,e}$. 
\EndProof

\begin{new-lemma}
\label{lemma-D}
Let $L\subset v^\perp$ be a $G$-invariant sublattice, such that 
$v^\perp/L$ is torsion of exponent $(v,v)$. The inclusion
$v^\perp\hookrightarrow(v^\perp)^*$, via the integral pairing
on $v^\perp$, maps $[v^\perp/L]^G$ onto 
$[(v^\perp)^*/L]^G$.
\end{new-lemma}

\noindent
{\bf Proof:} We may assume $v=(1,0,1-n)$, $n\geq 2$. Then 
$L=L_{d,e}$, for some positive integers $d$, $e$ satisfying
$\lcm(d,e)=(v,v)$, by Lemma \ref{lemma-B}. 
Let $u$ be the element $\ell a + t (w/(v,v))$ of $(v^\perp)^*$,
where $\ell,t$ are integers, $w=(1,0,n-1)$, 
and $a\in H^2(S,\Integers)$ is a primitive
class. Assume, that the coset $u+L_{d,e}$ is $G$-invariant. 
We need to show that $t\equiv 0$, modulo $(v,v)$. 
It suffices to prove, that both $d$ and $e$ divide $t$.


Let $U$ be the hyperbolic plane. 
Choose a primitive 
embedding $\iota:U\hookrightarrow H^2(S,\Integers)$, 
whose image contains $a$. 
Set $r:= (v,v)+1$. Let $b$ be a primitive class orthogonal to 
$\iota(U)$ in $H^2(S,\Integers)$, 
satisfying $(b,b)=(1-r^2)/(w,w)$.
Then the self-intersection of $(w,w)b+rw$ is equal to $(w,w)$.
The orthogonal complement $\nabla$, of $\iota(U)$ in the Mukai lattice
$H^*(S,\Integers)$, is unimodular. 
Hence, there exists an isometry $g\in O^+(\nabla)$, 
such that $g(v)=v$ and $g(w)=(w,w)b+rw$, 
by Lemma \ref{lemma-primitive-embeddings-of-rank-2-lattices}
part \ref{lemma-item-transitivity-on-pairs} applied to the two pairs 
$\{v,(v-w)/(v,v)\}$ and $\{v,[v-(w,w)b-rw]/(v,v)\}$. 
Extend $g$ to the Mukai lattice via the identity on $\iota(U)$. 
Then $g(u)-u=t[w-b]$. The difference $g(u)-u$ belongs to $L_{d,e}$, by 
the $G$-invariance of the coset $u+L_{d,e}$. Thus, 
both $d$ and $e$ divide $t$.
\EndProof

\medskip
{\bf Proof of Lemma 
\ref{lemma-class-of-extension-of-K-top-S-mod-L-by-L-mod-L'}:}
We reduce first the general case to the case where
the exponent $d$ of $L/L'$ is precisely $(v,v)$. 
Let $\epsilon_{(\ref{eq-extension-of-K-top-S-mod-L-by-L})}$
be the extension class of (\ref{eq-extension-of-K-top-S-mod-L-by-L})
and $\epsilon_d$ its pushforward via $L\rightarrow L/L'$. 
The order of $\epsilon_{(\ref{eq-extension-of-K-top-S-mod-L-by-L})}$ is 
$(v,v)$, by Lemma
\ref{lemma-orders-of-extension-classes-with-Mukai-lattice-in-the-middle}.
Set $\ell:=\lcm(d,(v,v))$, $h:=\frac{(v,v)}{\gcd(d,(v,v))}$, 
and $L_\ell:=hL'$. 
The exponent of $L/L_\ell$ is $\ell$.
Let $L_{(v,v)}$ be the kernel of the composition
$
L   \longrightarrow 
L/L_\ell
\LongRightArrowOf{\frac{d}{\gcd(d,(v,v))}} 
L/L_\ell,
$
where the second homomorphism is multiplication by 
$\frac{d}{\gcd(d,(v,v))}$.
Then $L/L_{(v,v)}$ has exponent $(v,v)$. 
Let $\epsilon_{\ell}$ and $\epsilon_{(v,v)}$ be the extension classes of the
pushforward of (\ref{eq-extension-of-K-top-S-mod-L-by-L})
via $L\rightarrow L/L_\ell$ and $L\rightarrow L/L_{(v,v)}$, respectively.
Assume, that the order of $\epsilon_{(v,v)}$ is $(v,v)$. 
We claim, that the order of $\epsilon_{\ell}$ is $(v,v)$ as well.
Indeed, the order of $\epsilon_\ell$ 
divides the order $(v,v)$ of 
$\epsilon_{(\ref{eq-extension-of-K-top-S-mod-L-by-L})}$ and 
is divisible by the order 
$(v,v)$ of $\epsilon_{(v,v)}$, since the latter is the pushforward of 
$\epsilon_\ell$ under 
$L/L_\ell\rightarrow L/L_{(v,v)}$. 

We have the short exact sequence 
$0\rightarrow L'/hL'\rightarrow L/L_\ell\rightarrow L/L'\rightarrow 0$.
The exponents $h$ of the kernel and $d$ of the co-kernel are relatively 
prime. Hence, the sequence splits $G$-equivariantly, via 
the isomorphism $L/L'\cong hL/L_\ell$. 
Thus, the order of $\epsilon_d$ is equal to that of 
$h\cdot \epsilon_\ell$, which is $\gcd(d,(v,v))$.

We may and thus do assume, that the exponent of $L/L'$ is 
precisely $(v,v)$. 
The homomorphism $L\rightarrow L/L'$ 
pushes the short exact sequence
(\ref{eq-extension-of-K-top-S-mod-L-by-L})
to 
\begin{equation}
\label{eq-extension-of-span-v-over-2n-2-by-v-perp-mod-L}
0\rightarrow L/L'
\rightarrow K(S)/L' \rightarrow K(S)/L 
\rightarrow 0.
\end{equation}
Denote the composition 
$K(S)/L\rightarrow K(S)/[\Integers v+v^\perp]\cong 
L^*/L$ by $f$. 
The extension $(\ref{eq-extension-of-span-v-over-2n-2-by-v-perp-mod-L})$
is the pullback of the short exact sequence 
\[
0\rightarrow L/L'\rightarrow L^*/L'\rightarrow 
L^*/L\rightarrow 0
\]
via $f$. 
Apply the functor $\Hom_G(K(S)/L,\bullet)$ to
the above short exact sequence.
The connecting homomorphism
$\delta:\Hom_G(K(S)/L,L^*/L)\rightarrow
\Ext^1_G(K(S)/L,L/L')$ takes $f$ to 
$\epsilon_{(\ref{eq-extension-of-span-v-over-2n-2-by-v-perp-mod-L})}$.
The homomorphism $f$ is surjective, and is hence an element of order $(v,v)$
in the domain of $\delta$. 
It remains to prove that $\delta$ is injective. 
Equivalently, we need to prove that
the homomorphism $\eta:\Hom_G(K(S)/L,L/L')\rightarrow 
\Hom_G(K(S)/L,L^*/L')$ is surjective. 
When $L=v^\perp$, the integral $G$-module $K(S)/L$ is trivial of 
rank $1$, and $\eta$ is an isomorphism,  by Lemma \ref{lemma-D}.

When $L=\Integers v$, the homomorphism $\eta$ factors as the composition
\[
[v^\perp/(v,v)v^\perp]^G \ \ \ \LongRightArrowOf{\alpha} \ \ \ 
[(v^\perp)^*/(v,v)v^\perp]^G \ \ \ \LongRightArrowOf{\beta} \ \ \ 
[\frac{1}{(v,v)}v^\perp/(v,v)v^\perp]^G.
\]
The homomorphism $\alpha$ is an isomorphism,  by Lemma \ref{lemma-D}.
We prove that $\beta$ is an isomorphism as well. The cokernel of
$\beta$ is a subgroup of $[\frac{1}{(v,v)}v^\perp/(v^\perp)^*]^G$,
which is isomorphic to 
$[v^\perp/(v,v)(v^\perp)^*]^G$. It suffices to prove
its vanishing when $v=(1,0,1-n)$. 
Then $(v,v)(v^\perp)^*$ is the sublattice $L_{(v,v),1}$ of $v^\perp$, given in
(\ref{eq-L-d-e}). 
The quotient $v^\perp/L_{(v,v),1}$ is isomorphic to
$H^2(S,\Integers)/(v,v)H^2(S,\Integers)$ and 
$[H^2(S,\Integers)/(v,v)H^2(S,\Integers)]^G$ is contained in 
$[H^2(S,\Integers)/(v,v)H^2(S,\Integers)]^{O^+H^2(S,\Integers)}$, 
which vanishes, by Lemma 
\ref{lemma-E}.
\EndProof

}

\journal{\begin{footnotesize}}

\journal{\end{footnotesize}}

\end{document}